\documentclass[11pt]{article}

\usepackage{amsmath}
\usepackage{amsfonts}
\usepackage{hyperref}
\usepackage{enumitem}
\usepackage[square,sort,comma,numbers]{natbib}
\usepackage{graphicx}
\usepackage{multirow}
\usepackage[table]{xcolor}
\usepackage{tikz}
\usetikzlibrary{shapes.geometric,arrows,positioning}
\tikzstyle{startstop} = [rectangle, rounded corners, minimum width=2cm, minimum height=0.6cm,text centered, draw=black, fill=red!20]
\tikzstyle{process}   = [rectangle, minimum width=2cm, minimum height=0.6cm, text centered, draw=black, fill=orange!20, align=center]
\tikzstyle{decision}  = [rectangle, rounded corners, aspect=2, minimum width=2cm, minimum height=0.6cm, text centered, draw=black, fill=blue!20, inner sep=5pt, align=center]
\tikzstyle{arrow}     = [thick,->,>=stealth]
\usepackage[ruled,vlined,linesnumbered]{algorithm2e}
\SetKwRepeat{Do}{do}{while}

\definecolor{lightgray}{gray}{0.9}

\newcommand{\R}{\mathbb{R}}
\DeclareMathOperator*{\argmin}{\mathrm{argmin}}

\begin{document}

\title{Local search and trajectory metaheuristics for the flexible job
  shop scheduling problem with sequencing flexibility and
  position-based learning effect\thanks{This work has been partially
    supported by the Brazilian agencies FAPESP (grants 2013/07375-0,
    2018/24293-0, and 2022/05803-3) and CNPq (grants 311536/2020-4 and
    302073/2022-1).}}

\author{
    K. A. G. Araujo\thanks{Department of Applied Mathematics, Institute of
    Mathematics and Statistics, University of S\~ao Paulo, Rua do
    Mat\~ao, 1010, Cidade Universit\'aria, 05508-090, S\~ao Paulo, SP,
    Brazil. e-mail: kennedy94@ime.usp.br}
 \and 
    E. G. Birgin\thanks{Department of Computer Science, Institute of
    Mathematics and Statistics, University of S\~ao Paulo, Rua do
    Mat\~ao, 1010, Cidade Universit\'aria, 05508-090, S\~ao Paulo, SP,
    Brazil. e-mail: egbirgin@ime.usp.br}
 \and
    D. P. Ronconi\thanks{Department of Production Engineering, Polytechnic 
    School, University of S\~ao Paulo, Av. Prof. Luciano Gualberto, 1380, 
    Cidade Universit\'aria, 05508-010 São Paulo, SP, Brazil. e-mail: 
    dronconi@usp.br}
}

\date{March 22, 2024}

\maketitle

\begin{abstract}
The flexible job shop scheduling problem with sequencing flexibility and position-based learning effect is considered in the present work. In [K. A. G. Araujo, E. G. Birgin, and D. P. Ronconi, Models, constructive heuristics, and benchmark instances for the flexible job shop scheduling problem with sequencing flexibility and position-based learning effect, Technical Report MCDO02022024, 2024], models, constructive heuristics, and benchmark instances for the same problem were introduced. In the present work, we are concerned with the development of effective and efficient methods for its resolution. For this purpose, a local search method and four trajectory metaheuristics are considered. In the local search, we show that the classical strategy of only reallocating operations that are part of the critical path can miss better quality neighbors, as opposed to what happens in the case where there is no learning effect. Consequently, we analyze an alternative type of neighborhood reduction that eliminates only neighbors that are not better than the current solution. In addition, we also suggest a neighborhood cut and experimentally verify that this significantly reduces the neighborhood size, bringing efficiency, with minimal loss in effectiveness. Extensive numerical experiments with the local search and the metaheuristics are carried on. The experiments show that tabu search, built on the reduced neighborhood, when applied to large-sized instances, stands out in relation to other the other three metaheuristics, namely, iterated local search, greedy randomized adaptive search procedure, and simulating annealing. Experiments with classical instances without sequencing flexibility show that the introduced methods also stand out in relation to methods from the literature. All the methods introduced, as well as the instances and solutions found, are freely available. As a whole, we build a test suite that can be used in future work.\\

\noindent
\textbf{Keywords:} flexible job shop scheduling problem, routing
flexibility, sequencing flexibility, learning effect, local search,
reduced neighbourhood, trajectory metaheuristics.\\
\end{abstract}

\section{Introduction}

The flexible job shop (FJS) with sequencing flexibility is a production environment with a wide range of relevant practical applications, especially in the nowadays on-demand printing industry~\cite{Lunardi2020,Lunardi2021}. Today, companies in the print-on-demand business must deal with customized production and prioritize on-time delivery in an effort to meet their customers' needs. In this context, production activities are organized in flexible workshops to better manage the execution of the wide range of tasks demanded. Other businesses that fit into this production environment include the glass industry~\cite{AlvarezValdes2005}, the mold industry~\cite{Gan2002}, the scheduling of aircraft support operations in flight decks~\cite{Yu2017}, the scheduling of repairing orders in automobile collision repair shops~\cite{AndradePineda2020}, and the construction of production programs for steelmaking~\cite{DeMoerloose2023}. It is therefore important that the solving methods are prepared to cope with the most diverse characteristics of the actual problems encountered in this production environment. One such factor is the learning effect, i.e.\ how the processing time of an operation varies with the number of times it is executed. It is clear that using processing times that are not entirely consistent with reality can lead to inaccurate schedules and result in significant economic losses.

The FJS scheduling problem is an extension of the classical job shop (JS) scheduling problem in which each operation can be processed by one within a set of machines instead of a single machine. This characteristic is known as routing flexibility. Two additional features are considered in the present work: sequencing flexibility and learning effect. In the FJS without sequencing flexibility, there exists the concept of a task, which consists of a set of operations that must respect a sequential order of execution (first the first operation, then the second, then the third, etc). The sequencing flexibility consists in considering that the precedences between the operations of a same task are given by an arbitrary directed acyclic graph (DAG). In a classical scheduling problem, given an operation and a machine that can process this operation, a fixed processing time is given that corresponds to the time demanded by the machine to process the operation. The learning effect corresponds to the real-life ingredient that consists in the fact that a person learns through the execution of a repetitive task and, the more times they execute it, the faster they do it. In this work we consider a learning function that depends on the position that an operation occupies within the list of operations executed by a machine, i.e. a position-based learning effect function.

A recent literature review on the FJS scheduling problem was done in~\cite{DauzrePrs2024}, while a recent review of the FJS with sequencing flexibility was included in~\cite{arbro2023}. For applications in real-life problems see~\cite{AndradePineda2020,AlvarezValdes2005,birgin2014milp,DeMoerloose2023,Gan2002,Lunardi2020,Lunardi2021,Yu2017} and for development of methods for the FJS with sequencing flexibility and a variety of additional features see~\cite{Birgin2015,Cao2021,Gao2006,Kasapidis2023,Kasapidis2021,Kim2003,VitalSoto2020}.

The impact on a worker's qualification, by repeated processing of an operation, and the resulting reduction in the operation's processing time was investigated in~\cite{DeJong1957} in 1957. Since the publications of the pioneer works~\cite{Biskup1999,Cheng2000,Gupta1988} that introduced the concept of learning effect in scheduling problems near two decades ago, a large number of papers has been published devoted to this subject. Surveys and classifications can be found in~\cite{Azzouz2018,Biskup2008,Janiak2011,Pei2022}. However, few papers, which are reviewed below, address the FJS scheduling problem with learning effect and, to the best of our knowledge, none of them consider sequencing flexibility.

In \cite{TayebiAraghi2014} the FJS scheduling problem with sequence-dependent setup times as well as position-dependent learning and time-dependent deterioration effects is considered. The authors propose a hybrid metaheuristic that combines the genetic algorithm and variable neighborhood search for the purpose of minimizing the makespan. The same problem is considered in \cite{Azzouz2020}, where the problem is modeled as a bilevel optimization problem in which both levels have the same objective, namely, minimizing the makespan. The proposed method, called evolutionary bilevel optimization, constructs feasible solutions in a two-stage hierarchical approach that assigns operations to machines in the first stage and schedules operations in the second stage. In \cite{Wu2018} a dual resource FJS is considered in which a machine and a worker are required to process an operation. It is also assumed that there is flexibility in the choice of both, and therefore, for each operation there is a set of machines and a set of workers capable of processing it. The problem is described using a mathematical model and a hybrid method that combines genetic algorithms and variable neighborhood search is developed. In \cite{Wu2019} an FJS scheduling problem with operation processing time deterioration effect is considered with the objective of minimizing a function that combines makespan and energy consumption. The proposed methodology hybridizes pigeon-inspired optimization and simulated annealing. In \cite{Zhu2020}, the FJS scheduling problem with time-dependent learning effect is considered. The problem is multi-objective and the goal is to minimize the makespan, total carbon emissions, and total workers' cost. A memetic algorithm, based on NSGA-II, associated with a variable neighborhood search is proposed. Four different neighborhoods are proposed and special attention is given to the development of constructive heuristics to develop the initial population. In \cite{Peng2021} a dual resource (machine and worker) FJS with position-dependent learning effect is considered. Additionally, a transportation time between machines is also taken into account. The goal is to minimize makespan, energy consumption, and noise. Therefore, the problem is multi-objective. All of these ingredients come from a real-life sand casting problem. To tackle the problem, a discrete multi-objective imperial competition algorithm is developed in which a local search based on simulated annealing is used. In \cite{Li2023}, a multi-objective FJS scheduling problem with position-dependent learning effect, in which makespan and total carbon emissions are minimized, is considered. The problem is modeled as a mixed integer linear multiobjective optimization problem and an improved multiobjective sparrow search algorithm is developed.

The problem considered in the present work was recently considered in~\cite{arbro2023}, where integer linear programming and constraint programming models were introduced and compared. To improve the performance of exact methods applied to small-sized instances of the problem, constructive heuristics were also developed. In the present work we continue that work by developing local search strategies and metaheuristics that can compute good quality solutions to large-sized instances. (It should be noted that the problem is NP-hard since it contains the JS scheduling problem, knowingly NP-hard~\cite{Garey1976}, as a particular case). In~\cite{Mastrolilli2000} it was introduced, for the FJS scheduling problem, a neighborhood reduction that eliminates only neighbors that do not improve the current solution. The reduction is based on the fact that, given a feasible solution, removing and relocating an operation that is not in the ``critical path'' has no chance of leading to a better solution. This idea was extended to the FJS scheduling problem with sequencing flexibility in~\cite{Lunardi2020,Lunardi2021}. In the present paper, we begin by showing that the fundamental principle on which this neighborhood reduction is based does not hold when we consider learning effect. Consequently, a new neighborhood reduction is introduced in the present work. In addition, we also introduce a neighborhood cutoff that can leave better quality neighbors behind. However, numerical experiments will show that if on the one hand the neighborhood cutoff is about 90\%, the solution found gets worse by about 1\% on average. Based on the introduced local search, we then implement four trajectory metaheuristics. With them, we manage to find good quality solutions for large-sized instances and to find solutions very close to the optimal solutions in small-sized instances with known optimal solution.

The rest of this paper is organized as follows. In Section~\ref{ls} we elaborate on how a feasible solution of the problem can be represented by a DAG. This representation allows us to introduce in a simple way the concept of neighborhood and, in the sequel, a local search. Still in Section~\ref{ls}, we analyze different possibilities to reduce the number of neighbors of a current solution to be inspected. In Section~\ref{meta} we describe the considered metaheuristics, namely, iterated local search, greedy randomized adaptive search procedure, tabu search, and simulated annealing. Section~\ref{experiments} is devoted to extensive numerical experiments. Conclusions and lines of future work are stated in the final section.

\section{Local search} \label{ls}

The data of an instance of the FJS with sequencing flexibility and position-based learning effect consists of (a) a set of operations $\mathcal{O}$ and a set of machines $\mathcal{F}$; (b) for each operation $i \in \mathcal{O}$, a subset $\mathcal{F}_i \subseteq \mathcal{F}$ containing the machines that can process $i$; (c) for each operation-machine pair $(i,k)$ with $i \in \mathcal{O}$ and $k \in \mathcal{F}_i$, a standard processing time $p_{ik}$; and (d) a set of arcs $\widehat A \subseteq \mathcal{O} \times \mathcal{O}$ representing the precedence relations between the operations. The learning effect is given by a function $\psi_{\alpha}(p,r)$ which, given a standard processing time~$p$ and a position~$r$, returns the actual processing time of an operation with standard processing time~$p$ when processed at the $r$-th position of a machine. The parameter $\alpha>0$ represents the learning effect rate. In the current work, we consider $\psi_{\alpha}(p,r) = \lfloor 100 \, p / r^{\alpha} + 1/2 \rfloor$. A simple example of an instance is shown in Figure~\ref{fig1}.

A feasible solution to an instance of the FJS with sequencing flexibility and position-based learning effect can be represented by a DAG $G=(V,A)$ as depicted in Figure~\ref{fig2}a. This graph is sometimes referred to as a solution graph in the literature. The vertices of~$G$, represented by the set~$V$, correspond to the operations plus the dummy vertices~$s$ and~$t$, i.e.\ $V=\mathcal{O} \cup \{s,t\}$. The arcs, represented by the set~$A$, correspond to the arcs in $\widehat A$ representing the precedence relationships between operations (in black in the picture), arcs going from $s$ to operations that have no predecessors and arcs going to $t$ from operations that do not precede any other operation (in purple in the picture). Arcs that start from $s$ and arcs that reach $t$ are named dummy arcs. Furthermore, dashed arcs represent the assignment of operations to machines and the order in which the operations are processed by each machine. These arcs are named machine arcs. Each node $i \in V \setminus \{s,t\} = \mathcal{O}$, that is, each operation, has a value $w_i$ associated with it that represents its actual processing time, which is calculated with the learning function using the operation's standard processing time and the position that the operation occupies in the machine to which it was assigned. Nodes $s$ and $t$ are associated with the value zero. The longest path between nodes $s$ and $t$ is called critical path (highlighted in yellow in the figure), and its length corresponds to the makespan of the represented solution.

\begin{figure}[ht!]
\begin{center}
\begin{tabular}{cc}
\begin{minipage}{6cm}
\begin{tikzpicture}[node distance={12mm}, thick, main/.style = {draw, circle}] 
\node[main] (1) {1};
\node[main] (2) [right =of 1] {2};
\node[main] (3) [right =of 2] {3};

\node[main] (4) [below =of 1] {4};
\node[main] (5) [right =of 4] {5};

\draw[->] (1)--(2);
\draw[->] (2)--(3);
\draw[->] (4)--(5);
\end{tikzpicture}
\end{minipage}
&
\begin{minipage}{4cm}
\begin{tabular}{|cc|cc|}
\cline{3-4}
\multicolumn{2}{c|}{} & \multicolumn{2}{c|}{Machines}\\
\multicolumn{2}{c|}{} & 1 & 2 \\
\hline
\multirow{5}{*}{\rotatebox{90}{Operations$\phantom{a}$}}
&  1 & 1  & 1  \\
&  2 & 1  & 1  \\
&  3 & 1  & 1  \\
&  4 & 10 & 10 \\
&  5 & 1  & 1  \\
\hline
\end{tabular}
\end{minipage}
\end{tabular}
\end{center}
\caption{On the left, representation of the operations' precedence constraints by a DAG $D=(\mathcal{O},\widehat A)$, where $\mathcal{O} = \{ 1, 2, \dots, 5 \}$ represents the set of operations and $\widehat A = \{(1,2), (2,3), (4,5)\}$ is the set of arcs that represents the precedence constraints. In this simple example, precedence constraints are given by a linear order, i.e.\ there is no sequencing flexibility. This instance has two machines and each of the five operations can be processed in any of the two machines, i.e.\ $\mathcal{F}=\{1,2\}$ and $\mathcal{F}_i=\mathcal{F}$ for all $i \in \mathcal{O}$. This means that there is full routing flexibility. The table on the right shows the standard processing times $p_{ik}$ of the five operations on each of the two machines.}
\label{fig1}
\end{figure}
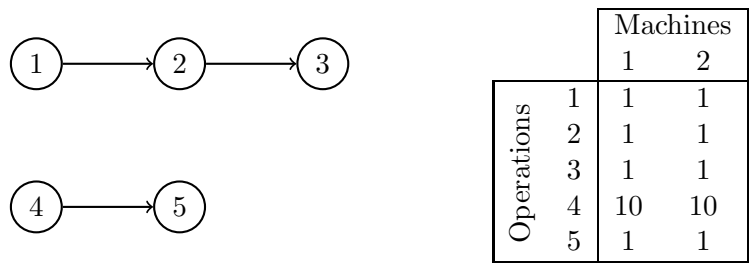

\begin{figure}[ht!]
\centering
\begin{tabular}{ccc}
\begin{tikzpicture}[node distance={14mm}, thick, main/.style = {draw, circle}] 
\node[main] (1) [label={[orange]90:$w_1=100$}] {1};
\node[main] (2) [right =of 1, label={[cyan]90:$w_2=100$}] {2};
\node[main] (3) [right =of 2, label={[orange]90:$w_3=25$}] {3};
\node[main] (6) [below left =of 1] {$s$};
\node[main] (4) [right =of 6, label={[orange]-90:$w_4=500$}] {4};
\node[main] (5) [right =of 4, label={[orange]-90:$w_5=33$}] {5};
\node[main] (7) [below right =of 3] {$t$};
\draw[->,yellow,ultra thick] (6) -- (1) -- (4) to [bend right] (5);
\draw[->,yellow,ultra thick] (5) -- (3) -- (7);
\draw[->] (1)--(2);
\draw[->] (2)--(3);
\draw[->] (4)--(5);
\draw[->,purple] (6) -- (1);
\draw[->,purple] (6) -- (4);
\draw[->,purple] (3) -- (7);
\draw[->,purple] (5) -- (7);
\draw[->,orange,dashed] (1) to (4);
\draw[->,orange,dashed] (4) to [bend right] (5);
\draw[->,orange,dashed] (5) to (3);
\end{tikzpicture}
&&
\begin{tikzpicture}[node distance={12mm}, thick, main/.style = {draw, circle}] 
\node[main] (1) [label={[orange]90:$w_1=100$}] {1};
\node[main] (2) [right =of 1, label={[orange]90:$w_2=50$}] {2};
\node[main] (3) [right =of 2, label={[orange]90:$w_3=20$}] {3};
\node[main] (6) [below left =of 1] {$s$};
\node[main] (4) [right =of 6, label={[orange]-90:$w_4=333$}] {4};
\node[main] (5) [right =of 4, label={[orange]-90:$w_5=25$}] {5};
\node[main] (7) [below right =of 3] {$t$};
\draw[->,yellow,ultra thick] (6) -- (1) to [bend left] (2);
\draw[->,yellow,ultra thick] (2)-- (4) to [bend right] (5);
\draw[->,yellow,ultra thick] (5) -- (3) -- (7);
\draw[->] (1)--(2);
\draw[->] (2)--(3);
\draw[->] (4)--(5);
\draw[->,purple] (6) -- (1);
\draw[->,purple] (6) -- (4);
\draw[->,purple] (3) -- (7);
\draw[->,purple] (5) -- (7);
\draw[->,orange,dashed] (1) to [bend left] (2);
\draw[->,orange,dashed] (2) to (4);
\draw[->,orange,dashed] (4) to [bend right] (5);
\draw[->,orange,dashed] (5) to (3);
\end{tikzpicture}\\
(a)&&(b)
\end{tabular}
\caption{In this figure we consider the instance in Figure~\ref{fig1} with learning rate $\alpha=1$. The digraph on the left (Figure~\ref{fig2}a) represents a feasible solution in which machine 1 (associated with the cyan color) processes operation 2 only, while machine 2 (associated with the orange color) processes operations 1, 4, 5, and 3 in that order. The colored numbers represent the actual processing time of the operations, with the influence of the learning effect. The critical path, which length corresponds to the makespan, is given by the path $s, 1, 4, 5, 3, t$ (highlighted in yellow in the picture). The digraph on the right (Figure~\ref{fig2}b) represents the feasible solution obtained by reallocating operation~2, that was not in the critical path, from machine~1 to machine~2 between operations~1 and~4. The constructed feasible solution, with critical path given by $s, 1, 2, 4, 5, 3, t$, has a makespan smaller than the original one (528 versus 658).}
\label{fig2}
\end{figure}
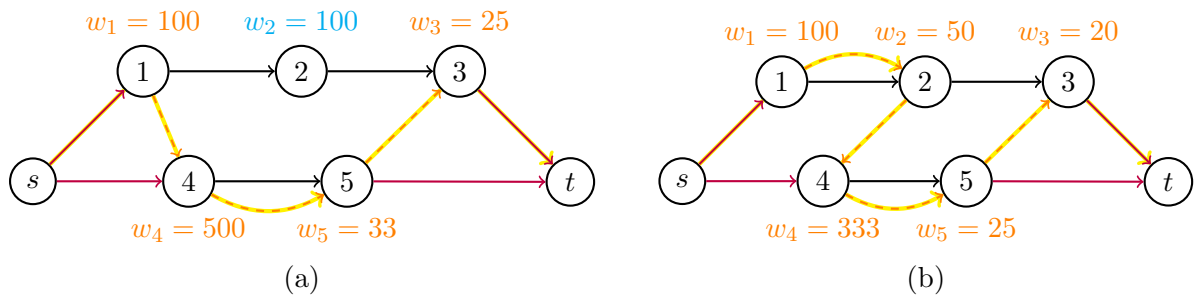

Given a feasible solution and a DAG $G=(V,A)$ representing it, a new feasible solution can be constructed by removing an operation from the machine to which it was assigned and reinserting it in the same machine but in another position or in another machine. When an operation is removed, the machine arcs adjacent to it must be removed and a new arc going from the operation before the removed one to the one after the removed one (if both exist) must be created. When the operation is reinserted, a similar reverse operation must also be done. When reinserting the operation, it is important to verify that a cycle is not produced in the digraph. Only reinsertions that do not create cycles build a digraph that corresponds to a feasible solution. When there is no learning effect, it is known~\cite{Mastrolilli2000} that there are chances to build a better feasible solution by removing and relocating operations that are part of the critical path only. If an operation is not part of the critical path, its removal and reinsertion cannot decrease the length of the critical path. It can increase it or create another even longer path. This is false when considering learning effect. An example is shown in Figure~\ref{fig2}b.

Given a feasible solution, we can define its neighborhood as the set of all feasible solutions that can be obtained by removing and reinserting a single operation. When there is no learning effect, just the removal and reinsertion of operations from the critical path may lead to better neighbors, i.e. neighbors with lower makespan. This fact is widely used to generate just promising neighbors. The observation in the previous paragraph shows that this reduction cannot be used in the problem we are considering in the present work. This leads us to analyze whether every removal and reinsertion that does not generate cycles has the potential to generate a neighbor with lower makespan or whether any neighborhood reduction is possible. The main point is to note that when an operation is removed from a machine, the operations that were scheduled to be processed later on that machine have their position decreased by one unit and, consequently, their actual processing time increased. Similarly, in the machine where the operation is inserted, the operations scheduled to be processed after the inserted operation have their position increased by one and, therefore, their processing time is decreased. These modifications may change the makespan whether improving or worsening.

Consider a feasible solution represented by a DAG $G=(V,A)$. Let $\mathcal{P}$ be a critical path in~$G$, with length $C_{\max}$. Let $i \in \mathcal{O}$ be an arbitrary operation. We name~$f_i$ the machine to which~$i$ is assigned. Let $k \in \mathcal{F}$ be an arbitrary machine. We name $Q_k = i_1, \dots, i_{|Q_k|}$ the ordered list of operations assigned to machine~$k$. If an operation $i$ is assigned to machine $f_i$ and it is in the $\gamma$-th position of $Q_{f_i}$, then its actual processing time is given by $w_i = \psi_{\alpha}(p_{i,f_i},\gamma)$. We intend to compute all neighbors of the solution represented by $G$, $f$, $Q$, and $w$. The neighbors will be constructed, for all $v \in \mathcal{O}$, by removing $v$ and reinserting $v$ at every possible place that does not generate a cycle. We want to determine if there are insertions that can be ignored because we know \textit{a priori} that they will not lead to a makespan reduction.

Let $v \in \mathcal{O}$ be an arbitrary operation. The computation of the neighbors of the current solution (associated with the remotion and reinsertion of~$v$) begins by computing a digraph $G_v^-=(V^-,A^-)$ in which the operation $v$ is removed. Such a graph is sometimes referred to as a reduced graph in the literature. The quantities $f^-$, $Q^-$, and $w^-$ associated with $G_v^-$ are also calculated. This digraph with its associated information is an intermediate structure necessary for the computation of the neighbors and, since the operation $v$ is not assigned to any machine, it does not represent a feasible solution. This task is described in Algorithm~\ref{algo1}. In line 2, the assignment of operations to machines is copied for all operations other than~$v$. In line~3, the lists of operations of all machines other than $f_v$ are copied, and, in the list of machine~$f_v$, operation~$v$ is eliminated. In line~4, the actual processing times of most of the operations are copied, except for the operations succeeding~$v$ in the list~$Q_{f_v}$. The actual processing times of these operations, whose position within the machine was reduced by one, need to be recalculated. The value of $w_v^-$ is set to zero. In line~5, machine arcs adjacent to~$v$ are removed (if they exist) and a new machine arc between its predecessor and successor is inserted, unless $v$ is the first or last element in~$Q_{f_v}$. In lines~6 and~7, the set $\mathcal{R}_v^{\leftarrow}$ of vertices that reaches $v$ and the set of vertices $\mathcal{R}_v^{\rightarrow}$ that are reached from $v$ in $G_v^-$ are calculated, which will be useful for detecting cycles in future possible reinsertions of~$v$. In line~8, the longest path $\mathcal{P}^-$ in the digraph $G_v^-$ is calculated, which will be useful to determine whether a reinsertion has chances to decrease the makespan or not. We call $\xi$ the length of $\mathcal{P}^-$. (We do not call it $C_{\max}^-$ because as $G_v^-$ does not represent a feasible solution, then the length of the path $\mathcal{P}^-$ does not represent a makespan.) Along with the computation of $\mathcal{P}^-$, the algorithm also computes, for each machine~$k$, the smallest position~$\tau_k$ such that, for all $\gamma>\tau_k$, the $\gamma$-th operation processed by machine~$k$ is not in~$\mathcal{P}^-$. (If machine~$k$ does not process any operation in $\mathcal{P}^-$, then $\tau_k=0$.) The tasks in lines~6, 7, and~8 use auxiliary well-known algorithms for topological sorting, depth-first search, and an adaptation~\cite[\S 22.2]{Cormen:2009:IA:1614191} of the Bellman-Ford algorithm described in Algorithms~\ref{algo2}, \ref{algo3}, and~\ref{algo4}, respectively, for completeness.

Let $G$ be the digraph, with associated quantities $f$, $Q$, and $w$, representing the current feasible solution. Let $\mathcal{P}$ be the critical path in $G$, with length~$C_{\max}$. Let $v$ be the operation we removed and wish to reinsert. Let $G_v^-$ be the digraph with $v$ removed and let $f^-$, $Q^-$, and $w^-$ be the quantities associated with $G_v^-$. Let $\mathcal{P}^-$ be the critical path in $G_v^-$, with length $\xi$, and, for each machine~$k$, let $\tau_k$ be the smallest position in $Q_k$ such that every operation in a position after~$\tau_k$ is not in~$\mathcal{P}^-$. Let~$\kappa$ be a machine and~$\gamma$ be a position in the list~$Q_\kappa^-$ such that inserting~$v$ at position~$\gamma$ of~$Q_\kappa^-$ does not generate a cycle. Does such an insertion have a chance of generating a new digraph whose associated feasible solution has a makespan smaller than~$C_{\max}$? If $\xi \geq C_{\max}$ and $\gamma > \tau_k$, then the answer is ``no''. This is because the path $\mathcal{P}^-$ with length~$\xi$ not smaller than~$C_{\max}$ already exists and the insertion of~$v$ in machine~$\kappa$, at a position~$\gamma$ later than~$\tau_k$, will not modify the actual processing time of any operation in $\mathcal{P}^-$. If $\xi < C_{\max}$ or $\xi \geq C_{\max}$ but $\gamma \leq \tau_k$, then chances exist. 

It should be noted that, strictly speaking, the fact of~$v$ being in~$\mathcal{P}$ or not is not related to the answer to the question above. But already anticipating something that will come later, as the neighborhood reduction driven by the answer to the question may be quite small, we will consider in the experiments, heuristically, $v \in \mathcal{P}$ as equivalent or strongly correlated to $\xi<C_{\max}$. That is, we will consider that removing an operation from the critical path will most likely imply that $\xi < C_{\max}$. This is very plausible for moderate values of the learning factor~$\alpha$, with which a possible reduction of one in the machine position of some critical path operations does not cancel out the benefit of removing an operation from the critical path.

\begin{algorithm}[!ht]
\caption{Computes $G^-_v=(V^-_v,A^-_v)$, $f^-$, $Q^-$, and $w^-$ by removing operation $v$ from the solution graph $G = (V, A)$. Then, in $G^-_v$, it computes the set $\mathcal{R}_v^{\leftarrow}$ of vertices that reaches $v$ and the set of vertices $\mathcal{R}_v^{\rightarrow}$ that are reached from $v$, the largest path $\mathcal{P^-}$ from $s$ to $t$ and its length $\xi$, and the position $\tau_k$ of the last critical operation in each machine $k \in \mathcal{F}$.}
\label{algo1}
\KwIn{$\mathcal{O}$, $\mathcal{F}$, $p$, $v$, $f$, $Q$, $w$, $G=(V,A)$}
\KwOut{$f^-$, $Q^-$, $w^-$, $G^-_v=(V^-_v,A^-_v)$, $\mathcal{P}^-$, $\xi$, $\mathcal{R}_v^{\leftarrow}$, $\mathcal{R}_v^{\rightarrow}$, $\tau$}
\SetKwBlock{Begin}{function}{end function}
\DontPrintSemicolon
\Begin($\text{RemoveOp}{(}\mathcal{O}$, $p$, $v$, $f$, $Q$, $w$, $G$, $f^-$, $Q^-$, $w^-$, $G^-_v$, $\mathcal{P}^-$, $\xi$, $\mathcal{R}_v^{\leftarrow}$, $\mathcal{R}_v^{\rightarrow}$, $\tau${)})
{
Define $f^-_i := f_i$ for all $i \in \mathcal{O} \setminus \{v\}$.\;
Let $Q_{f_v}$ be given by the sequence $i_1, \dots, i_{\gamma-1}, i_{\gamma}, i_{\gamma+1}, \dots, i_{|Q_{f_v}|}$ with $i_{\gamma}=v$. Define $Q^-_k := Q_k$ for all $k \in \mathcal{F} \setminus \{f_v\}$ and $Q^-_{f_v} := i_1, \dots, i_{\gamma-1}, i_{\gamma+1}, \dots, i_{|Q_{f_v}|}$.\;
Define $w^-_i := w_i$ for all $i \in \mathcal{O} \setminus \{ v, i_{\gamma+1}, \dots, i_{|Q_{f_v}|}\}$, $w_v^-:=0$, and $w_{i_\ell}^- := \psi_{\alpha}(p_{i_\ell,f_{i_\ell}},\ell-1)$ for $\ell = \gamma+1, \dots, |Q_{f_v}|$.\;
Define the graph $G^-_v = (V^-_v,A^-_v)$ with $V^-_v := V$ and $A^-_v := ( A \setminus \{ (i_{\gamma-1},v), (v,i_{\gamma+1}) \} ) \cup \{ (i_{\gamma-1},i_{\gamma+1}) \}$.\;
Initialize $\mathcal{V} \gets \emptyset$, $\mathcal{U}$ as an empty list, and $\mathcal{R}_v^{\leftarrow} \gets \{v\}$, and compute in $\mathcal{U}$ and $\mathcal{R}_v^{\leftarrow}$ a topological sort of the vertices in $V^ -_v$ and the set of vertices $i \in V^-_v$ such that a path from $i$ to $v$ exists, respectively, by calling TopologicalSort+($G^-_v$, $\mathcal{U}$, $s$, $\mathcal{V}$, $\mathcal{R}_v^{\leftarrow}$).\;
Initialize $\mathcal{R}_v^{\rightarrow} \gets \emptyset$ and compute in $\mathcal{R}_v^{\rightarrow}$ the set of vertices $i \in V^-_v$ such that a path from $v$ to $i$ exists, by calling DFS($G^-_v$, $v$, $\mathcal{R}_v^{\rightarrow}$).\;
Compute the largest path $\mathcal{P^-}$ from $s$ to $t$, its length $\xi$, and determine the position~$\tau_k$ of the last critical operation at each machine $k$ for all $k \in \mathcal{F}$ by calling CriticalPath($\mathcal{F}$, $f^-$, $w^-$, $Q^-$, $G^-_v = (V^-_v,A^-_v)$, $\mathcal{U}$, $\mathcal{P}^-$, $\xi$, $\tau$).\;
}
\end{algorithm}

\begin{algorithm}[!ht]
\caption{Computes a topological sort $\mathcal{U}$ of the vertices of $G = (V, A)$. In addition, if~$\mathcal{R}_v^{\leftarrow}$ is present as an input parameter, computes the set $\mathcal{R}_v^{\leftarrow}$ of vertices that reaches $v$ in $G = (V, A)$.} 
\label{algo2}
\KwIn{$G = (V,A)$, $\mathcal{U}$, $v$, $\mathcal{V}$, $\mathcal{R}_v^{\leftarrow}$}
\KwOut{$\mathcal{U}$, $\mathcal{V}$, $\mathcal{R}_v^{\leftarrow}$}
\DontPrintSemicolon
\SetKwBlock{Begin}{function}{end function}
\Begin(TopologicalSort+{(}$G$, $\mathcal{U}$, $v$, $\mathcal{V}$, $\mathcal{R}_v^{\leftarrow}${)})
{
Set $\mathcal{V} \gets \mathcal{V} \cup \{v\}$.\;
\For{$j \; \mathrm{such} \; \mathrm{that} \; (v,j) \in A$}{
    \If{$j \notin \mathcal{V}$}{
        TopologicalSort+($G$, $\mathcal{U}$,  $j$, $\mathcal{V}$, $\mathcal{R}_v^{\leftarrow}$)\;
        }
    \If{$i \not\in \mathcal{R}_v^{\leftarrow} \;\mathrm{and} \;  j \in \mathcal{R}_v^{\leftarrow}$}{
        set $\mathcal{R}_v^{\leftarrow} \gets \mathcal{R}_v^{\leftarrow} \cup \{ v\}$.\;}
    }
    Insert $i$ at the beginning of $\mathcal{U}$.\;
}
\end{algorithm}

\begin{algorithm}[!ht]
\caption{Computes the set of vertices $\mathcal{R}_v^{\rightarrow}$ as the set of vertices that can be reached by $v$ in $G = (V, A)$.} 
\label{algo3}
\KwIn{$G = (V,A)$, $v$, $\mathcal{R}_v^{\rightarrow}$}
\KwOut{$\mathcal{R}_v^{\rightarrow}$}
\DontPrintSemicolon
\SetKwBlock{Begin}{function}{end function}
\Begin(DFS{(}$G$, $v$, $\mathcal{R}_v^{\rightarrow}${)})
{
Set $\mathcal{R}_v^{\rightarrow} \gets \mathcal{R}_v^{\rightarrow} \cup \{v\}$.\;
\For{$j\; \mathrm{such} \; \mathrm{that} \;(v,j) \in A$}{
    \If{$j \notin \mathcal{R}_v^{\rightarrow}$}{
        DFS($G$, $j$, $\mathcal{R}_v^{\rightarrow}$)\;
        }
    }
}
\end{algorithm}

\begin{algorithm}[ht!]
\caption{Computes a critical path $\mathcal{P}$ and its length $\xi$ for a given graph $G = (V, A)$. In addition, if $\tau$ is present as an input parameter, determines the last critical operation in each machine.}
\label{algo4}
\KwIn{$\mathcal{F}$, $f$, $w$, $Q$, $G = (V,A)$, $\tau$}
\KwOut{$\mathcal{U}$, $\mathcal{P}$, $\xi$, $\tau$}
\SetKwBlock{Begin}{function}{end function}
\DontPrintSemicolon
\Begin($\text{CriticalPath}{(}\mathcal{F}$, $f$, $w$, $Q$, $G$, $\mathcal{U}$, $\mathcal{P}$, $\xi$, $\tau${)}){
    
Initialize $d_i \leftarrow -\infty$ for all $i \in V\setminus\{s\}$ and define $d_s := 0$ and $\pi_s := 0$.\;

Initialize $\mathcal{V} \leftarrow \emptyset$ and $\mathcal{U}$ as an empty list and compute in $\mathcal{U}$ a topological sort of the vertices in $V$, by calling TopologicalSort+($G$, $\mathcal{U}$, $s$, $\mathcal{V}$).\;

\For{$\ell=1,\dots,|V|$}{
    Let $i$ be the $\ell$-th operation in the topological order given by $\mathcal{U}$.\;
    \For{$j\; \mathrm{such} \; \mathrm{that} \;(i,j) \in A$}{
        \If{$d_j < d_i + w_i$}{
            $d_j \gets d_i + w_i$ and $\pi_j \gets i$.\;}}}

$\xi := d_t$\;

Initialize $i \gets \pi_t$, $\mathcal{P} \leftarrow \emptyset$,
and $\tau_k \leftarrow 0$ for all $k \in \mathcal{F}$.\;

\Do{$i \neq s$}{
  \If{$\tau_{f_i} = 0$}{
      Let $Q_{f_i}$ be given by the sequence $i_1,
      \dots,i_{\ell-1},i,i_{\ell+1},\dots,i_{|Q_{f_i}|}$. Define
      $\tau_{f_i} := \ell$.}
  
  $\mathcal{P} \gets \mathcal{P} \cup \{i\}$ and $i \gets \pi_i$.}
}
\end{algorithm}

The task of reinserting $v$ in $G_v^-$ at position $\gamma$ of machine $\kappa$ generates a DAG that we name $G_v^+$. This task is similar to the removing task. The construction of $G_v^+$, its associated quantities $f^+$, $Q^+$, and $w^+$, and its critical path~$\mathcal{P}^+$ with length~$C_{\max}^+$ is described in Algorithm~\ref{algo5}. The only relevant detail that remains to be explained is how to determine whether an insertion generates a cycle or not. A cycle will only be created in $G_v^+$ if $v$ is inserted in a position that leaves some $u \in \mathcal{R}_v^{\leftarrow}$ to be processed after $v$ in machine $\kappa$ or some $u \in \mathcal{R}_v^{\rightarrow}$ to be processed before $v$ in machine $\kappa$. The limits $\underline{\gamma}$ and $\bar \gamma$ such that $\underline{\gamma} + 1 \leq \gamma \leq \bar \gamma$ avoids cycles are calculates in lines~7 and~8. A possible reduction of this interval is computed in lines~9 and~10, eliminating the possibility of making insertions after $\tau_\kappa$ if $\xi \geq C_{\max}$, as already discussed. Algorithm~\ref{algo6} implements a best neighbor local search with the neighbor reduction already discussed. It corresponds to a classical local search with a best neighbor strategy.

\begin{algorithm}[!ht]
\caption{Computes $G^+_v=(V^+_v,A^+_v)$, $f^+$, $Q^+$, and $w^+$ by inserting operation $v$ at position $\gamma$ of machine $\kappa$ in the reduced graph $G^-_v$. Then, in $G^+_v$, it computes the largest path $\mathcal{P^+}$ from $s$ to $t$ and its length $C_{\max}^+$.} 
\label{algo5}
\KwIn{$\mathcal{O}$, $\mathcal{F}$, $p$, $v$, $\gamma$, $\kappa$, $f^-$, $Q^-$, $w^-$, $G^-_v = (V^-_v,A^-_v)$}
\KwOut{$f^+$, $Q^+$, $w^+$, $G^+_v=(V^+_v,A^+_v)$, $\mathcal{P}^+$, $C_{\max}^+$}
\SetKwBlock{Begin}{function}{end function}
\DontPrintSemicolon
\Begin(InsertOp{(}$\mathcal{O}$, $p$, $v$, $\gamma$, $\kappa$, $f^-$, $Q^-$, $w^-$, $G^-_v$, $f^+$, $Q^+$, $w^+$, $G^+_v$, $\mathcal{P}^+$, $C_{\max}^+${)})
{
Define $f^+_i := f^-_i$ for all $i \in \mathcal{O} \setminus \{v\}$ and $f^+_v := \kappa$.\;
Let $Q^-_{\kappa}$ be given by the sequence $i_1, i_2, \dots, i_{|Q^-_{\kappa}|}$.\;
Define $Q^+_k := Q^-_k$ for all $k \in \mathcal{F} \setminus \{\kappa\}$ and $Q^+_{\kappa} := i_1, \dots, i_{\gamma-1}, v, i_{\gamma}, \dots, i_{|Q^-_{\kappa}|}$.\;
Define $w^+_i := w^-_i$ for all $i \in \mathcal{O} \setminus \{ v, i_{\gamma}, \dots, i_{|Q_{\kappa}|}\}$, and $w_{i_\ell}^+ := \psi(p_{i_\ell,\kappa},\ell)$ for $\ell = \gamma, \dots, |Q^+_{\kappa}|$.\;
Define the graph $G^+_v = (V^+_v , A^+_v)$ with $V^+_v := V^-_v$ and $A^+_v := ( A^-_v \setminus \{(i_{\gamma-1},i_{\gamma}))\} \cup \{ (i_{\gamma-1},v),(v,i_{\gamma}) \}$. \;
Initialize $\mathcal{V} \gets \emptyset$ and $\mathcal{U}$ as an empty list and compute in $\mathcal{U}$ a topological sort of the vertices in $V^+_v$, by calling TopologicalSort+($G^+_v,$ $\mathcal{U}$, $s$, $\mathcal{V}$).\;
Compute the critical path $\mathcal{P^+}$ from $s$ to $t$ and its length $C_{\max}$ by calling CriticalPath($\mathcal{F}$, $f^+$, $w^+$, $Q^+$, $G^+_v = (V^+_v,A^+_v)$, $\mathcal{U}$, $\mathcal{P}^+$, $C_{\max}^+$).\;
}
\end{algorithm}

\begin{algorithm}[!ht]

\caption{Best neighbor local search with reduced neighbor.}
\label{algo6}
\KwIn{$\mathcal{O}$, $\mathcal{F}$, $p$, $G = (V,A)$, $f$, $w$, $Q$, $\mathcal{P}$, $C_{\max}$}

\KwOut{$G^{\star} = (V^\star,A^\star)$, $f^{\star}$, $w^{\star}$, $Q^{\star}$, $\mathcal{P}^{\star}$, $C^\star_{\max}$}
\SetKwBlock{Begin}{function}{end function}
\DontPrintSemicolon
\Begin(LocalSearch{(}$\mathcal{O}$, $\mathcal{F}$, $p$, $G$, $f$, $w$, $Q$, $\mathcal{P}$, $C_{\max}$, $G^\star$, $f^\star$, $w^\star$, $Q^\star$, $\mathcal{P}^\star$, $C^\star_{\max}${)})
{
\Do{$\delta>0$}{
    $C_{\max}^{\mathrm{bn}} \gets +\infty $\;
    \For{$v \in \mathcal{O}$}{
        RemoveOp($\mathcal{O}$, $p$, $v$, $f$, $Q$, $w$, $G$, $f^-$, $Q^-$, $w^-$, $G^-_v$, $\mathcal{P}^-$, $\xi$, $\mathcal{R}_v^{\leftarrow}$, $\mathcal{R}_v^{\rightarrow}$, $\tau$)\;
	    \For{$k \in \mathcal{F}_v$}{
	        Let $\underline{\gamma}$ be the position of the last operation in $Q^-_k = i_1, \dots, i_{|Q^-_k|}$ such that $i_{\underline{\gamma}} \in \mathcal{R}_v^{\leftarrow}$ and $\underline{\gamma}=0$ if $i_{\ell} \not\in \mathcal{R}_v^{\leftarrow}$ for $\ell=1,\dots,|Q^-_k|$.\;
	        Let $\bar{\gamma}$ be the position of the first operation in $Q^-_k = i_1, \dots, i_{|Q^-_k|}$ such that $i_{\bar{\gamma}} \in \mathcal{R}_v^{\rightarrow}$ and $\bar{\gamma}=|Q^-_k| + 1$ if $i_{\ell} \not\in \mathcal{R}^{\rightarrow}_v$ for $\ell=1,\dots,|Q^-_k|$.\;
	        
	        \If{$\xi \geq  C_{\max}$}{
	            %\Comment{$O(N)$}
	            $\bar \gamma \gets \min\{ \bar \gamma, \tau_k\}$, where $\tau_k$ is such that there is no critical operation after $\tau_k$ in $Q^-_k$ ($\tau_k=0$ if there is  no critical operation in $Q^-_k$).\;
	            %\Comment{$O(N)$}
                }  
	        \For{$\gamma = \underline{\gamma}+1,\dots,\bar \gamma$}{
    	       InsertOp($\mathcal{O}$, $p$, $v$, $\gamma$, $k$, $f^-$, $Q^-$, $w^-$, $G^-_v, f^+, Q^+, w^+, G^+_v$, $\mathcal{P}^+$, $C^+_{\max}$)\;
    	   
    	       \If{$C_{\max}^+ < C_{\max}^{\mathrm{bn}}$}{%\Comment{$O(N)$}
                    $G^{\mathrm{bn}}, f^{\mathrm{bn}}, w^{\mathrm{bn}}, Q^{\mathrm{bn}}, \mathcal{P}^{\mathrm{bn}}, C_{\max}^{\mathrm{bn}} \gets G^+_v, f^+, w^+, Q^+, \mathcal{P}^+, C_{\max}^+$
                    }  
	           }
            }
        }
    $\delta \gets C_{\max} - C_{\max}^{\mathrm{bn}}$\;
    \If{$\delta > 0$}{
        $G, f, w, Q, \mathcal{P}, C_{\max} \gets G^{\mathrm{bn}}, f^{\mathrm{bn}}, w^{\mathrm{bn}}, Q^{\mathrm{bn}}, \mathcal{P}^{\mathrm{bn}}, C_{\max}^{\mathrm{bn}}$
        }
    }
$G^\star, f^\star, w^\star, Q^\star, \mathcal{P}^\star, C_{\max}^\star \gets G, f, w, Q, \mathcal{P}, C_{\max}$    
}
\end{algorithm}

\section{Metaheuristics} \label{meta}

It is known that the direct application of a local search as the one introduced in the previous section suffers from premature convergence to local solutions. Therefore, it is natural to think of using it in connection with metaheuristics. In particular, trajectory metaheuristics are a natural choice since they make direct use of local search strategies. In this paper we consider the well-known metaheuristics iterated local search (ILS), greedy randomized adaptive search procedure (GRASP), tabu search (TS), and simulating annealing (SA). The basic components of each of them are presented below.

\subsection{Iterated local search}

The ILS~\cite{Lourenco2003} strategy consists of iteratively running a local search until it converges to a local solution and perturbing the solution found to be used as the initial solution for the next run of the local search. For the first run of the local search, it uses as initial solution the best among the solutions constructed by the constructive heuristics ECT and EST introduced in~\cite{arbro2023}. The perturbation consists of a sequence of $\ell^p$ modifications, where $\ell^p$ is a random number between $\ell^p_{\min}$ and $\ell^p_{\max}$. ($\ell^p_{\min}$ and $\ell^p_{\max}$ are the only two parameters of the ILS.) Each modification consists of removing a random operation and relocating it to a random place so that the solution thus constructed is feasible. Algorithms~\ref{algo7} and~\ref{algo8} describe the ILS completely. They make use of Algorithm~\ref{algo6} for the local search and Algorithms~\ref{algo1} and~\ref{algo5} for removing and reinserting an operation, respectively. The random choices at line~12 of Algorithm~\ref{algo7} and lines~2, 4, and 7 of Algorithm~\ref{algo8} follow a discrete uniform distribution within the prescribed range. The same applies to all the methods described below, which make random choices following discrete or continuous (depending on the situation) uniform distributions within prescribed ranges.

\begin{algorithm}[!ht]
\caption{Iterated local search.} 
\label{algo7}
\KwIn{$\mathcal O$, $\mathcal F$, $p$, $\widehat A$, $\ell^p_{\min}$, $\ell^p_{\max}$}
\KwOut{$f^\star$, $w^\star$, $Q^\star$, $G^\star$, $C_{\max}^\star$}
\SetKwBlock{Begin}{function}{end function}
\DontPrintSemicolon
\Begin(ILS{(}$\mathcal{O}$, $\mathcal{F}$, $p$, $\widehat A$, $\ell_{\min}p$, $\ell_{\max}p$, $f^\star$, $w^\star$, $Q^\star$, $G^\star$, $C_{\max}^\star${)})
{
ECT($\mathcal{O}$, $\mathcal{F}$, $p$, $\widehat A$, $f^1$, $w^1$, $Q^1$, $G^1$, $\mathcal{U}^1$, $\mathcal{P}^1$, $C_{\max}^1$, $\tau^1$)\;

EST($\mathcal{O}$, $\mathcal{F}$, $p$, $\widehat A$, $f^2$, $w^2$, $Q^2$, $G^2$, $\mathcal{U}^2$, $\mathcal{P}^2$, $C_{\max}^2$, $\tau^2$)\;

\uIf{$C_{\max}^1 < C_{\max}^2$}{
  $f^\star \gets f^1$, $w^\star \gets w^1$, $Q^\star \gets Q^1$, $G^\star \gets G^1$ and $C_{\max}^\star \gets C_{\max}^1$\;}
\Else{
  $f^\star \gets f^2$, $w^\star \gets w^2$, $Q^\star \gets Q^2$, $G^\star \gets G^2$ and $C_{\max}^\star \gets C_{\max}^2$\;}
     
\While{the stopping criterion is not satisfied}{
  LocalSearch($\mathcal{O}$, $\mathcal{F}$, $p$, $G$, $f$, $w$, $Q$, $\mathcal{P}$, $C_{\max}$, $G'$, $f'$, $w'$, $Q'$, $\mathcal{P}'$, $C'_{\max}$)\;
            
  \If{$C'_{\max} < C^\star_{\max}$}{
    $f^\star \gets f'$, $w^\star \gets w'$, $Q^\star \gets Q'$, $G^\star \gets G'$ and $C_{\max}^\star \gets C'_{\max}$ \;
  }

  Set $\ell^p$ a random number between $\ell^p_{\min}$ and $\ell^p_{\max}$.
        
  \For{$i = 1, \dots, \ell^p$}{
    Perturb($\mathcal{O}$, $\mathcal{F}$, $p$, $f'$, $w'$, $Q'$, $G'$, $f$, $w$, $Q$, $G$, $\mathcal{P}$, $C_{\max}$)\;
  }
}
}
\end{algorithm}

\begin{algorithm}[!ht]
\caption{Perturbs a solution.} 
\label{algo8}
\KwIn{${\cal O}$, ${\cal F}$, $p$, $f$, $w$, $Q$, $G = (V, A)$}
\KwOut{$f'$, $w'$, $Q'$, $G' = (V', A')$, $\mathcal{P}'$, $C'_{\max}$}
\SetKwBlock{Begin}{function}{end function}
\DontPrintSemicolon
\Begin(Perturb{(}${\cal O}$, ${\cal F}$, $p$, $f$, $w$, $Q$, $G$, $f'$, $w'$, $Q'$, $G'$, $\mathcal{P}'$, $C'_{\max}${)})
{
Let $v \in V \setminus \{s,t\}$ be a random operation.\;

RemoveOp($\mathcal{O}$, $p$, $v$, $f$, $Q$, $w$, $G$, $f^-$, $Q^-$, $w^-$, $G^-_v$, $\mathcal{P}^-$,
$\xi$, $\mathcal R^{\leftarrow}_v$, $\mathcal R^{\rightarrow}_v$, $\tau$)\;

Let $k \in \mathcal{F}_{v}$ be a random machine and let $Q^-_{k} = i_1, \dots, i_{|Q^-_{k}|}$.\;
    
Let $\underline{\gamma}$ be the position of the last operation in $Q^-_{k}$ such that
$i_{\underline{\gamma}} \in \mathcal{R}_v^{\leftarrow}$ and let $\underline{\gamma}=0$ if
$i_{\ell} \not\in \mathcal{R}_v^{\leftarrow}$ for $\ell=1,\dots,|Q^-_k|$.\;
    
Let $\bar{\gamma}$ be the position of the first operation in $Q^-_{k}$ such that
$i_{\bar{\gamma}} \in \bar{R}_v^{\rightarrow}$ and let $\bar{\gamma}=|Q^-_k| + 1$ if
$i_{\ell} \not\in \mathcal{R}^{\rightarrow}_v$ for $\ell=1,\dots,|Q^-_k|$.\;
         
Let $\gamma \in \{ \underline{\gamma} + 1, \dots, \bar{\gamma}\}$ be a random feasible position.\;
    
InsertOp{(}$\mathcal{O}, p, v, \gamma, k,f^-,Q^-,w^-,G^-_v, f',Q', w', G'$, $\mathcal{P}'$, $C'_{\max}${)}
}
\end{algorithm}

\subsection{Greedy randomized adaptive adaptive procedure}

The GRASP~\cite{Feo1995} consists of iteratively generating an initial solution and running a local search starting from the initial solution just generated. The initial solutions are generated as the best among the solutions generated by randomized versions of the constructive heuristics ECT and EST introduced in~\cite{arbro2023}. 

The two constructive heuristics introduced in~\cite{arbro2023} are based on the earliest starting time (EST) rule~\cite{birgin2014milp} and the earliest completion time (ECT) rule~\cite{Leung2005}. The heuristics schedule one operation at a time. In the EST-based constructive heuristic, the instant~$r_{\min}$ which is the earliest instant at which an unscheduled operation could be initiated is computed first. All operation/machine pairs that could start at that instant are considered and the pair with the shortest processing time is selected. Since they all would start at instant~$r_{\min}$, selecting the pair with the shortest processing time is the same as selecting the pair that ends earliest. This idea is taken to the extreme in the constructive heuristic based on the ECT rule: without limiting the choice to the operation/machine pairs that could start as soon as possible, the operation/machine pair that will finish earliest is chosen, even if the processing of the operation does not start as soon as possible.

In the EST-based constructive heuristic, the randomization is done as follows. First, the earliest time $r_{\min}$ at which an unscheduled operation can be initiated is calculated. Then we calculate the shortest and the longest processing times $a$ and $b$ associated with the operation/machine pairs that can start at $r_{\min}$. Given $\alpha \in (0,1]$, a random pair is chosen from those whose processing time is between $a$ and $a + \alpha (b-a)$. In the ECT-based constructive heuristic, the randomization consists of, from among the unscheduled operations, calculating the lowest and highest completion instants $a$ and $b$ and drawing an operation/machine pair from among those that would terminate between $a$ and $a + \alpha (b-a)$. It is worth noting that the calculation of the effective processing time always takes into consideration the learning effect. The only parameter of the method is $\alpha$. The GRASP strategy is described in Algorithm~\ref{algo9} and the randomized versions of the constructive heuristics based on the EST and ECT dispatching rules are described in Algorithms~\ref{algo10} and~\ref{algo11}, respectively. Algorithm~\ref{algo9} uses Algorithm~\ref{algo6} to perform the local search and Algorithms~\ref{algo10} and~\ref{algo11} use Algorithm~\ref{algo4} for the calculation of the critical path associated with the initial solution, which is a necessary input for the local search. For more details on the deterministic versions of the constructive heuristics based on the EST and ECT dispatching rules see \cite[Algorithms~1 and~4]{arbro2023}.
  
\begin{algorithm}[!ht]
\caption{Greedy randomized adaptive procedure.} 
\label{algo9}
\KwIn{$\mathcal{O}$, $\mathcal{F}$, $p$, $\widehat A$, $\alpha$}
\KwOut{$f^\star$, $w^\star$, $Q^\star$, $G^\star$, $C_{\max}^\star$}
\SetKwBlock{Begin}{function}{end function}
\DontPrintSemicolon
\Begin(GRASP{(}$\mathcal{O}$, $\mathcal{F}$, $p$, $\widehat A$, $\alpha$, $f^\star$, $w^\star$, $Q^\star$, $G^\star$, $C_{\max}^\star${)})
{
Initialize $C_{\max}^\star \gets +\infty$\;
     
\While{the stopping criterion is not satisfied}{
  RandomizedECT($\mathcal{O}$, $\mathcal{F}$, $p$, $\widehat A$, $\alpha$, $f^1$, $w^1$, $Q^1$, $G^1$, $\mathcal{P}^1$, $C_{\max}^1$) \;
        
  RandomizedEST($\mathcal{O}$, $\mathcal{F}$, $p$, $\widehat A$, $\alpha$, $f^2$, $w^2$, $Q^2$, $G^2$, $\mathcal{P}^2$, $C_{\max}^2$) \;

  \uIf{$C_{\max}^1 < C_{\max}^2$}{
    $f \gets f^1$, $w \gets w^1$, $Q \gets Q^1$, $\mathcal{P} \gets \mathcal{P}^1$, $G \gets G^1$, and $C_{\max} \gets C^1_{\max}$ \;
  }
  \Else{
    $f \gets f^2$, $w \gets w^2$, $Q \gets Q^2$, $\mathcal{P} \gets \mathcal{P}^2$, $G \gets G^2$, and $C_{\max} \gets C^2_{\max}$ \;
  }
  
  LocalSearch($\mathcal{O}$, $\mathcal{F}$, $p$, $G$, $f$, $w$, $Q$, $\mathcal{P}$, $C_{\max}$, $G'$, $f'$, $w'$, $Q'$, $\mathcal{P}'$, $C'_{\max}$) \;
  
  \If{$C'_{\max} < C^\star_{\max}$}{
    $f^\star \gets f'$, $w^\star \gets w'$, $Q^\star \gets Q'$, $G^\star \gets G'$ and $C_{\max}^\star \gets C'_{\max}$ \;
  }
}
}
\end{algorithm}

\begin{algorithm}[!ht]
\caption{Computes a feasible solution using a randomization of the earliest starting time (EST) rule.} 
\label{algo10}
\KwIn{$\mathcal{O}$, $\mathcal{F}$, $p$, $\widehat A$, $\alpha$}
\KwOut{$f$, $w$, $Q$, $G$, $\mathcal{P}$, $C_{\max}$}
\SetKwBlock{Begin}{function}{end function}
\DontPrintSemicolon
\Begin(RandomizedEST{(}$\mathcal{O}$, $\mathcal{F}$, $p$, $\widehat A$, $\alpha$, $f$, $w$, $Q$, $G$, $\mathcal{P}$, $C_{\max}${)})
{
Set $A \gets \widehat A \cup \{ (s,j) \mid (\cdot,j) \not\in \widehat A \} \cup \{ (i,t) \mid (i,\cdot) \not\in \widehat A \}$ and define $V := \mathcal{O} \cup \{s,t\}$ and $G = (V,A)$.\;
    
Set $r^{\mathrm{op}}_v \gets +\infty$ and define $r^{\mathrm{op}}_s := 0$, $w_s := w_t := 0$, and $c_s:=0$. \;
    
Set $r^{\mathrm{mac}}_k \gets 0$ and $g_k \gets 1$ for all $k \in \mathcal{F}$.\;
    
Initialize $\Pi \gets V \setminus\{s,t\}$ as the set of non-scheduled operations, and $Q_k$ as an empty list for all $k \in \mathcal{F}$.\; 
    
\While{$\Pi \neq \emptyset$}{
  \For{$v \in \Pi$}{
    \If{$\Pi \cap \{i \;|\:(i,v) \in A \}  = \emptyset$}{
      $r^{\mathrm{op}}_v \gets \max \{ c_i \mid i \in V\setminus \Pi$ such that $(i,v) \in A \}$\;
    }
  }
  
  Set $r_{\min} = {\displaystyle \min_{v \in \Pi, k \in \mathcal{F}_v} \{ \max( r^{\mathrm{op}}_v,r^{\mathrm{mac}}_k) \}}$ and
  let $E$ be the set of pairs $(v,k)$ with $v \in \Pi$ and $k \in \mathcal F_v$ such that $\max( r^{\mathrm{op}}_v,r^{\mathrm{mac}}_k) = r_{\min}$.\;

  Compute $ \underline{\psi} \gets {\displaystyle \min_{(v,k) \in E} \{ \psi(p_{v, k}, g_k) \}}$ and
  $\overline{\psi} \gets {\displaystyle \max_{(v,k) \in E} \{ \psi(p_{v, k}, g_k) \}}$, set
  $\mathcal{L} \gets \{(v,k)  \in E \mid \psi(p_{v, k}, g_k) \leq \underline{\psi} + \alpha (\overline{\psi} - \underline{\psi}) \}$
  and choose $(\hat v, \hat k) \in \mathcal{L}$ randomly.\;
    
  Define $w_{\hat v} := \psi(p_{\hat v,\hat k}, g_{\hat k})$, $f_{\hat v} := \hat k$, $c_{\hat v} := c_{\hat v, \hat k}$ and
  set $r^{\mathrm{mac}}_{\hat k} \gets c_{\hat v}$ and $g_{\hat k} \gets g_{\hat k} + 1 $.\;
        
  \If{$|Q_{\hat k}| \neq 0$}{
    Set $A \gets A \cup \{ (i_{|Q_{\hat k}|},\hat v) \} $, where $Q_{\hat k} = i_1, \dots, i_{|Q_{\hat k}|}$.\;
  }
  
  Insert $\hat v$ at the end of $Q_{\hat k}$ and set $\Pi \gets \Pi \setminus \{\hat v\}$. \;
}

CriticalPath($\mathcal{F}$, $f$, $w$, $Q$, $G$, $\mathcal{U}$, $\mathcal{P}$, $C_{\max}$, $\tau$).
}
\end{algorithm}

\begin{algorithm}[!ht]
  \caption{Computes a feasible solution using a randomization of the earliest completion time (ECT) rule
    .} 
\label{algo11}
\KwIn{$\mathcal{O}$, $\mathcal{F}$, $p$, $\widehat A$, $\alpha$}
\KwOut{$f$, $w$, $Q$, $G$, $\mathcal{P}$, $C_{\max}$}
\SetKwBlock{Begin}{function}{end function}
\DontPrintSemicolon
\Begin(RandomizedECT{(}$\mathcal{O}$, $\mathcal{F}$, $p$, $\widehat A$, $\alpha$, $f$, $w$, $Q$, $G$, $\mathcal{P}$, $C_{\max}${)})
{
Set $A \gets \widehat A \cup \{ (s,j) \mid (\cdot,j) \not\in \widehat A \} \cup \{ (i,t) \mid (i,\cdot) \not\in \widehat A \}$
and define $V := \mathcal{O} \cup \{s,t\}$ and $G = (V,A)$.\;
    
Set $r^{\mathrm{op}}_v \gets +\infty$ and define $r^{\mathrm{op}}_s := 0$, $w_s := w_t := 0$, and $c_s:=0$. \;
    
Set $r^{\mathrm{mac}}_k \gets 0$ and $g_k \gets 1$ for all $k \in \mathcal{F}$.\;
    
Initialize $\Pi \gets V \setminus\{s,t\}$ as the set of non-scheduled operations, and $Q_k$ as an empty list for all $k \in \mathcal{F}$.\; 
    
\While{$\Pi \neq \emptyset$}{
  \For{$v \in \Pi$}{
    \If{$\Pi \cap \{i \;|\:(i,v) \in A \}  = \emptyset$}{
      $r^{\mathrm{op}}_v \gets \max \{ c_i \mid i \in V\setminus \Pi$ such that $(i,v) \in A \}$\;
    }
  }
  \For{$v \in \Pi$}{
    \For{$k \in \mathcal{F}_v$}{
      $c_{v, k} \gets \max ( r^{\mathrm{op}}_v,r^{\mathrm{mac}}_k) + \psi(p_{v, k}, g_k)$
    }
  }
  Set $c^{\min} \gets \min \{ c_{v, k} \mid v \in \Pi, k \in \mathcal F_v\}$,
  $c^{\max} \gets \max \{ c_{v, k} \mid v \in \Pi, k \in \mathcal F_v\}$ and
  $\mathrm{RCL} \gets \{(v, k) \mid c_{v, k} \leq c^{\min} + \alpha(c^{\max} - c^{\min}), v \in \Pi, k \in \mathcal F_v \}$.\;
  
  Choose $(\hat v, \hat k)$ randomly in RCL.\;
    
  Define $w_{\hat v} := \psi(p_{\hat v,\hat k}, g_{\hat k})$, $f_{\hat v} := \hat k$,
  $c_{\hat v} := c_{\hat v, \hat k}$ and set $r^{\mathrm{mac}}_{\hat k} \gets c_{\hat v}$ and $g_{\hat k} \gets g_{\hat k} + 1 $.\;
        
  \If{$|Q_{\hat k}| \neq 0$}{
    Let $Q_{\hat k} = i_1, \dots, i_{|Q_{\hat k}|}$. Set $A \gets A \cup \{ (i_{|Q_{\hat k}|},\hat v) \} $.\;
  }
  
  Insert $\hat v$ at the end of $Q_{\hat k}$ and set $\Pi \gets \Pi \setminus \{\hat v\}$. \;
}
    
CriticalPath($\mathcal{F}$, $f$, $w$, $Q$, $G$, $\mathcal{U}$, $\mathcal{P}$, $C_{\max}$, $\tau$).
}
\end{algorithm}

\subsection{Tabu search}

The TS~\cite{Glover1989,Glover1990} considered consists basically of a local search with a modification in the acceptance criteria of a neighbor as a new solution. The acceptance depends on a list of tabu moves. The size $t_{\max}$ of this list is the only parameter of the method. In a given iteration, the best neighbor is computed, disregarding neighbors constructed with a tabu move unless they correspond to the best solution already constructed. The way of constructing the neighborhood follows exactly the same scheme as the local search described in Algorithm~\ref{algo6}. The move that transforms the solution of the current iteration into the solution of the next iteration consists in relocating a certain operation $v$ at the $\gamma$ position of a machine $k \in \mathcal{F}_v$. We keep in the tabu list the operation/machine pair $(v,k)$. The complete method is described in Algorithm~\ref{algo12}.

\begin{algorithm}[!ht]
\caption{Tabu Search.}
\label{algo12}
\KwIn{$ \mathcal{O}$, $\mathcal{F}$, $p$, $\widehat A$, $t_{\max}$}
\KwOut{$f^{\star}$, $w^{\star}$, $Q^{\star}$, $G^\star$, $C_{\max}^\star$}
\SetKwBlock{Begin}{function}{end function}
\DontPrintSemicolon
\Begin(TS{(}$\mathcal{O}$, $\mathcal{F}$, $p$, $\widehat A$, $t_{\max}$, $f^\star$, $w^\star$, $Q^\star$, $G^\star$, $C_{\max}^\star${)})
{
ECT($\mathcal{O}$, $\mathcal{F}$, $p$, $\widehat A$, $f^1$, $w^1$, $Q^1$, $G^1$, $\mathcal{U}^1$, $\mathcal{P}^1$, $C_{\max}^1$, $\tau^1$)\;

EST($\mathcal{O}$, $\mathcal{F}$, $p$, $\widehat A$, $f^2$, $w^2$, $Q^2$, $G^2$, $\mathcal{U}^2$, $\mathcal{P}^2$, $C_{\max}^2$, $\tau^2$)\;

\uIf{$C_{\max}^1 < C_{\max}^2$}{
  $f^\star \gets f^1$, $w^\star \gets w^1$, $Q^\star \gets Q^1$, $G^\star \gets G^1$ and $C_{\max}^\star \gets C_{\max}^1$\;}
\Else{
  $f^\star \gets f^2$, $w^\star \gets w^2$, $Q^\star \gets Q^2$, $G^\star \gets G^2$ and $C_{\max}^\star \gets C_{\max}^2$\;}
     
Initialize $\mathcal T$ as an empty list.\;

\While{the stopping criterion is not satisfied}{
  $C'_{\max} \gets +\infty $\;
  \For{$v \in \mathcal{O}$}{
    RemoveOp$(\mathcal{O}$, $p$, $v$, $f$, $Q$, $w$, $G$, $f^-$, $Q^-$, $w^-$, $G^-_v$, $\mathcal{P}^-$, $\xi$,
    $\mathcal{R}_v^{\leftarrow}$, $\mathcal{R}_v^{\rightarrow}$, $\tau$)\;
        
    \For{$k \in \mathcal{F}_v$}{
      Let $\underline{\gamma}$ be the position of the last operation in $Q^-_k = i_1, \dots, i_{|Q^-_k|}$
      such that $i_{\underline{\gamma}} \in \mathcal{R}_v^{\leftarrow}$ and $\underline{\gamma}=0$
      if $i_{\ell} \not\in \mathcal{R}_v^{\leftarrow}$ for $\ell=1,\dots,|Q^-_k|$.\;
      
      Let $\bar{\gamma}$ be the position of the first operation in $Q^-_k = i_1, \dots, i_{|Q^-_k|}$ such that
      $i_{\bar{\gamma}} \in \bar{R}_v^{\rightarrow}$ and $\bar{\gamma}=|Q^-_k| + 1$
      if $i_{\ell} \not\in \mathcal{R}^{\rightarrow}_v$ for $\ell=1,\dots,|Q^-_k|$.\;
	        
      \If{$\xi \geq  C_{\max}$}{
	%\Comment{$O(N)$}
	$\bar \gamma \gets \min\{ \bar \gamma, \tau_k\}$, where $\tau_k$ is such that there is no critical
        operation after $\tau_k$ in $Q^-_k$ ($\tau_k=0$ if there is  no critical operation in $Q^-_k$).\;
	%\Comment{$O(N)$}
      }  

      \For{$\gamma = \underline{\gamma}+1,\dots,\bar \gamma$}{
    	InsertOp($\mathcal{O}$, $p$, $v$, $\gamma$, $\kappa$, $f^-$, $Q^-$, $w^-$, $G^-_v, f^+, Q^+, w^+, G^+_v$, $\mathcal{P}^+$, $C^+_{\max}$)\;

    	\If{$( C^+_{\max} < C'_{\max} \, \mathrm{ and } \, (v,k) \notin \mathcal{T} ) \, \mathrm{or} \,
          C^+_{\max} < \min\{C'_{\max}, C^\star_{\max} \}$}{%\Comment{$O(N)$}
          $G' \gets G^+_v$, $f' \gets f^+$, $w' \gets w^+$, $\mathcal{P}' \gets \mathcal{P}^+$ and $C'_{\max} \gets C_{\max}^+$ \;
          $v' \gets v,\; k' \gets k$\;
        }  
      }
    }
  }
  
  If $(v', k') \in \mathcal{T}$ then remove it from $\mathcal{T}$. Anyhow, insert $(v', k')$ at the end of $\mathcal T$.\;

  \If{$|\mathcal T| > t_{\max}$}{
    Remove the first element of $\mathcal T$
  }
  
  $G \gets G'_v$,  $f \gets f'$, $w \gets w'$,$\mathcal{P} \gets \mathcal{P}'$ and $C_{\max} \gets C'_{\max}$\;
    
  \If{$C_{\max} < C^\star_{\max}$}{
    $G^\star \gets G_v$,  $f^\star \gets f$, $w^\star \gets w'$, $\mathcal{P}^\star \gets \mathcal{P}$ and $C^\star_{\max} \gets C_{\max}$\;
  }
}
}
\end{algorithm}

\subsection{Simulated annealing}

In the SA~\cite{Kirkpatrick1983}, as in the other metaheuristics, the initial solution is given by the best solution among those provided by the constructive heuristics based on the EST and ECT dispatch rules. The simulating annealing advances from one solution to another by making perturbations of the solution consisting of a random movement of removal and insertion (described in Algorithm~\ref{algo8}). The main feature is that, for a new solution to be accepted, it need not be better than the current solution. On the other hand, the probability of accepting a solution that is worse than the current solution decreases as the method progresses, i.e.\ as the temperature decreases. The method has as parameters the initial and final temperatures $T_0$ and $T_f$. The temperature decreases as the method progresses being multiplied by a value $\delta$ that varies between given parameters $0 < \delta_{\min} < \delta_{\max} < 1$. The amount of perturbations for a fixed temperature varies between the parameters $\bar \ell_{\min}$ and $\bar \ell_{\max}$ and increases as the temperature decreases. The method as a whole is described in Algorithm~\ref{algo13}.

\begin{algorithm}[!ht]
\caption{Simulated Annealing.}
\label{algo13}
\KwIn{$\mathcal{O}$, $\mathcal{F}$, $p$, $\widehat A$, $\bar{\ell}$, $\delta$, $T_0$, $T_f$}
\KwOut{$f^\star$, $w^\star$, $Q^\star$, $G^\star$, $C_{\max}^\star$}
\SetKwBlock{Begin}{function}{end function}
\DontPrintSemicolon
\Begin(SA{(}$\mathcal{O}$, $\mathcal{F}$, $p$, $\widehat A$, $\bar{\ell}$, $\delta$, $T_0$, $T_f$, $f^\star$, $w^\star$, $Q^\star$, $G^\star$, $C_{\max}^\star${)})
{
ECT($\mathcal{O}$, $\mathcal{F}$, $p$, $\widehat A$, $f^1$, $w^1$, $Q^1$, $G^1$, $\mathcal{U}^1$, $\mathcal{P}^1$, $C_{\max}^1$, $\tau^1$)\;

EST($\mathcal{O}$, $\mathcal{F}$, $p$, $\widehat A$, $f^2$, $w^2$, $Q^2$, $G^2$, $\mathcal{U}^2$, $\mathcal{P}^2$, $C_{\max}^2$, $\tau^2$)\;

\uIf{$C_{\max}^1 < C_{\max}^2$}{
  $f^\star \gets f^1$, $w^\star \gets w^1$, $Q^\star \gets Q^1$, $G^\star \gets G^1$ and $C_{\max}^\star \gets C_{\max}^1$\;}
\Else{
  $f^\star \gets f^2$, $w^\star \gets w^2$, $Q^\star \gets Q^2$, $G^\star \gets G^2$ and $C_{\max}^\star \gets C_{\max}^2$\;}
     
Set $T \gets T_0$.\;

\While{the stopping criterion is not satisfied}{
        
  \For {$\ell = 1, \dots, \bar \ell$}{
    Perturb{(}${\cal O}$, ${\cal F}$, $p$, $f$, $w$, $Q$, $G$, $f'$, $w'$, $Q'$, $G'$, $\mathcal{P}'$, $C'_{\max}${)}
            
    Set $\Delta C \gets ( C'_{\max}- C_{\max} ) / C_{\max}$ and let $r \in [0,1]$ be a random number.
            
    \If{$e^{-\Delta C /T} \geq r$}{
      Set $f \gets f'$, $w \gets w'$, $Q \gets Q'$, $G \gets G'$ and $C_{\max} \gets C'_{\max}$. \;
                
      \If{$C'_{\max} < C^\star_{\max}$}{
        Set $f^\star \gets f'$, $w^\star \gets w'$, $Q^\star \gets Q'$, $G^\star \gets G'$ and $C_{\max}^\star \gets C'_{\max}$. \;
      }
    }
  }
  
  $T \gets \max(\delta \, T, T_f)$
}
}
\end{algorithm}

\section{Numerical experiments} \label{experiments}

In this section, we present numerical experiments with the proposed local search and the considered metaheuristics. In all experiments, we consider the~60 small-sized instances introduced in~\cite{arbro2023} and the~50 large-sized instances introduced in~\cite{birgin2014milp}, with learning rate $\alpha \in \{ 0.1, 0.2, 0.3 \}$, totaling~330 instances. Small-sized instances have between 5 and 7 machines and between 9 and 24 operations. The operations are divided into up to 6 tasks, which correspond to connected components of the DAG of precedence relations, and the DAGs have up to 21 edges. The MILP models~\cite[Eq.(1--4, p.7)]{arbro2023} of the small-sized instances have up to almost 1{,}000 binary variables and 12{,}000 constraints. Large-sized instances have between 5 and 26 machines and between 25 and 289 operations. The operations are divided into up to 17 tasks and the DAGs representing the precedence relations have up to 272 edges. Their MILP models have up to 73{,}000 binary variables and 3{,}800{,}000 constraints. It is worth recalling that the considered problem has two types of flexibility: routing flexibility and sequencing flexibility. The first refers to the fact that an operation can be processed by a machine within a set of machines instead of a single machine and corresponds to the ``flexible'' of the FJS. The second corresponds to the fact that the operations of the same task have their precedences represented by an arbitrary DAG instead of obeying a linear order. The sum of these two flexibilities causes the problem to have a large search space and makes it difficult to find a proven optimal solution even in instances that may initially appear to be simple. For a more detailed description of the instances and their characteristics, see~\cite{arbro2023}.

The local search and the metaheuristics were implemented in C++ programming language. The code was compiled using g++ 10.2.1. To ensure full reproducibility of the results presented in the present work, as well as future comparisons, the instances, code and solutions found are available at \url{https://github.com/kennedy94/FJS}. The experiments were carried out in an Intel i9-12900K (12th Gen) processor operating at 5.200GHz and 128 GB of RAM. 

\subsection{Experiments with local search variations}

In this section we evaluate variations of the local search described in Algorithm~\ref{algo6}. In~\cite{arbro2023} two constructive heuristics are introduced for the same problem being considered in this paper. The constructive heuristics are different and, while one constructs, in general, better solutions for instances of type Y, the other constructs better solutions for instances of type DA; see~\cite{arbro2023} for details. Either way, the two constructive heuristics take negligible time and, for that reason, in the present work we use, as the initial solution for the local search, the best among the solutions constructed by the two constructive heuristics introduced in~\cite{arbro2023}.

As described in Algorithm~\ref{algo6}, the local search uses the best neighbor strategy and makes use of the neighborhood reduction described in Section~\ref{ls}. Therefore, we call this version of ``local search with the best neighbor strategy and reduced neighborhood''. The neighborhood reduction is implemented in lines~9 and~10. If we remove those two lines we obtain a version that we call ``local search with the best neighbor and full neighborhood''. The version with reduced neighborhood does not consider neighbors that are guaranteed to be no better than the current solution. Therefore, the solution obtained with the reduced neighborhood must be identical to the solution obtained with the full neighborhood. (In fact all iterates of the two versions must be identical and not just the final solution). Only a reduction of CPU time is expected. As already mentioned in Section~\ref{ls}, we decided to consider yet another version that would exhibit a more drastic reduction in CPU time, albeit with a possible loss of quality in the solution. We call this version of ``local search with the best neighbor and cropped neighborhood''. This version consists of changing $v \in \mathcal{O}$ to $v \in \mathcal{P}$ in line~4 of Algorithm~\ref{algo6}. That is, only operations on the critical path are reallocated, since there is a greater tendency for these reallocations to generate better quality neighbors. We have then three different versions of the local search with the best neighbor strategy that are distinguished by the neighborhood used: full neighborhood, reduced neighborhood and cropped neighborhood. Each of them corresponds to minimal variations of Algorithm~\ref{algo6} as already described. Furthermore, we consider the same three versions but using the strategy of interrupting the inspection of the neighborhood when finding the first neighbor that improves the current solution, i.e.\ the first-improvement strategy. This change corresponds, in Algorithm~\ref{algo6}, to interrupting the loop of line~4 the first time line~16 is executed.

The results of the six variations of the local search applied to the~330 considered instances are shown in Tables~\ref{tab1}--\ref{tab4}. Tables~\ref{tab1} and~\ref{tab2} correspond to the first-improvement strategy while Tables~\ref{tab3} and~\ref{tab4} correspond to the best-improvement strategy. In the tables, we show the makespan of the obtained solution, the number of iterations that the local search made until finding an iterate that is better than all its neighbors (this is the stopping criterion as described in Algorithm~\ref{algo6}), and the CPU time in seconds. In the case of the small-sized instances introduced in~\cite{arbro2023}, the tables do not show the CPU time because it was always less than~1 millisecond. The tables also do not show anything related to the full neighborhood. What should be said about the use of the full neighborhood is that, in all instances, as expected, the solution obtained was identical to the solution obtained with the reduced neighborhood, the number of iterations was also the same, and the reduced neighborhood promoted a reduction of 52.51\% in the CPU time.

When we compare the first-improvement and best-improvement strategies (Tables~\ref{tab1} and~\ref{tab2} versus Tables~\ref{tab3} and~\ref{tab4}), the results are quite similar, but the best-improvement strategy always finds better quality solutions using less CPU time. Specifically, the best improvement strategy returns solutions that are on average 1.02\% and 0.70\% better than the solutions returned by the first-improvement strategy, when we consider the reduced and cropped neighborhoods, respectively. Therefore, from now on, we focus on evaluating the reduced neighborhood and the cropped neighborhood associated with the best-improvement strategy.

The cropped neighborhood eliminates, on average, 90.34\% of the neighbors of the reduced neighborhood, promoting a proportional reduction in CPU time. However, adopting the cropped neighborhood can lead to a loss of quality in the final solution obtained by the local search method. On average, when compared to the local search with the reduced neighborhood, the local search with the cropped neighborhood finds solution with a makespan 0.69\% worse. When we compare the final solution with the initial solution, the local search using the reduced neighborhood improves the initial solution by, on average, 6.88\%, while the local search using the cropped neighborhood improves the initial solution by 6.11\%. In conclusion, the local search with the cropped neighborhood is significantly faster than the local search with the reduced neighborhood and finds solutions that are only slightly worse than the solutions found by the latter.

%Acredito que esse dado também deveria estar. O gap em relação ao ótimo apenas de instâncias com ótimo comprovado, no total de 180 contando com mini e médias, para reduced e cropped respectivamente.
%XX8% = 6.65%
%XX9% = 7.48%

%Esses dados também seriam bons de ser incluídos, do total de 330 instâncias.
%XX10 = 265 -- Número de vezes em que Best improvement achou um makespan melhor do que First Improvement. 
%XX11 = 245 -- Número de vezes em que Best Improvement foi mais rápido do que First Improvement. 
%XX12 = 193 -- Número de vezes que o Best Improvement foi mais rápido e achou makespan melhor do que First Improvement

\begin{table}[ht!]
\begin{center}  
\resizebox{!}{0.45\textheight}{
\begin{tabular}{|c|cc|cc|cc|cc|cc|cc|} 
\hline
\multirow{3}{*}{instance} & \multicolumn{4}{c|}{$ \alpha = 0.1$} & \multicolumn{4}{c|}{$ \alpha = 0.2$} & \multicolumn{4}{c|}{$ \alpha = 0.3$} \\ 
\cline{2-13}
& \multicolumn{2}{c|}{Cropped} & \multicolumn{2}{c|}{Reduced} & \multicolumn{2}{c|}{Cropped} & \multicolumn{2}{c|}{Reduced} & \multicolumn{2}{c|}{Cropped} & \multicolumn{2}{c|}{Reduced} \\ 
\cline{2-13}
& $C_{\max}$ & \#it & $C_{\max}$ & \#it & $C_{\max}$ & \#it & $C_{\max}$ & \#it & $C_{\max}$ & \#it & $C_{\max}$ & \#it \\ 
\hline
\hline
miniDAFJS01 & 23,264 &  1 & 23,019 &  2 & 22,099 & 3 & 21,741 & 4 & 20,007 & 4 & 19,767 & 7 \\
miniDAFJS02 & 23,242 &  1 & 22,927 &  4 & 22,161 & 1 & 22,046 & 2 & 21,136 & 2 & 20,987 & 2 \\
miniDAFJS03 & 18,363 &  1 & 18,363 &  1 & 17,972 & 1 & 17,972 & 1 & 17,619 & 2 & 17,619 & 2 \\
miniDAFJS04 & 21,690 &  1 & 21,151 &  4 & 20,757 & 2 & 19,602 & 5 & 19,867 & 2 & 18,800 & 5 \\
miniDAFJS05 & 22,598 &  1 & 22,418 &  2 & 20,826 & 1 & 20,523 & 2 & 19,253 & 1 & 18,878 & 2 \\
miniDAFJS06 & 23,370 &  1 & 23,370 &  1 & 20,783 & 2 & 20,783 & 2 & 18,712 & 3 & 18,712 & 3 \\
miniDAFJS07 & 28,644 &  1 & 28,644 &  1 & 25,088 & 1 & 24,715 & 2 & 24,636 & 2 & 24,256 & 4 \\
miniDAFJS08 & 19,878 &  1 & 19,878 &  1 & 18,857 & 1 & 18,857 & 1 & 17,900 & 2 & 17,900 & 2 \\
miniDAFJS09 & 25,425 &  1 & 25,425 &  1 & 23,830 & 1 & 23,715 & 2 & 22,398 & 1 & 22,258 & 2 \\
miniDAFJS10 & 21,563 &  1 & 20,359 &  4 & 21,021 & 1 & 18,823 & 6 & 19,723 & 2 & 17,466 & 6 \\
miniDAFJS11 & 33,689 &  2 & 33,686 &  3 & 30,941 & 2 & 30,550 & 5 & 28,510 & 2 & 27,907 & 6 \\
miniDAFJS12 & 20,342 &  1 & 20,342 &  1 & 19,624 & 1 & 19,624 & 1 & 19,094 & 1 & 19,094 & 1 \\
miniDAFJS13 & 17,091 &  2 & 16,313 &  3 & 16,616 & 2 & 15,143 & 3 & 15,278 & 2 & 15,053 & 2 \\
miniDAFJS14 & 23,248 &  4 & 23,140 &  5 & 21,990 & 5 & 21,817 & 5 & 20,711 & 5 & 20,620 & 5 \\
miniDAFJS15 & 22,472 &  1 & 22,472 &  1 & 20,653 & 4 & 20,620 & 5 & 18,628 & 3 & 20,328 & 2 \\
miniDAFJS16 & 25,691 &  1 & 25,426 &  2 & 24,593 & 1 & 24,114 & 2 & 22,939 & 3 & 22,939 & 3 \\
miniDAFJS17 & 21,070 &  1 & 20,155 &  3 & 20,183 & 2 & 20,183 & 2 & 19,139 & 2 & 19,139 & 2 \\
miniDAFJS18 & 18,512 &  3 & 18,135 &  7 & 18,829 & 1 & 18,201 & 4 & 17,784 & 1 & 17,295 & 2 \\
miniDAFJS19 & 21,293 &  1 & 20,945 &  3 & 20,107 & 1 & 19,642 & 3 & 19,030 & 1 & 18,474 & 3 \\
miniDAFJS20 & 23,443 &  2 & 23,212 &  3 & 21,286 & 1 & 21,286 & 1 & 19,587 & 1 & 19,587 & 1 \\
miniDAFJS21 & 24,404 &  3 & 24,274 &  4 & 22,215 & 2 & 22,151 & 6 & 21,182 & 3 & 20,368 & 5 \\
miniDAFJS22 & 25,923 &  4 & 25,923 &  4 & 24,273 & 4 & 24,273 & 4 & 22,767 & 3 & 22,767 & 3 \\
miniDAFJS23 & 25,839 &  2 & 25,788 &  3 & 24,253 & 2 & 24,164 & 3 & 21,772 & 5 & 21,561 & 6 \\
miniDAFJS24 & 26,932 &  1 & 26,610 &  4 & 24,500 & 1 & 23,709 & 3 & 22,524 & 1 & 21,714 & 4 \\
miniDAFJS25 & 23,370 &  1 & 22,982 &  2 & 22,613 & 1 & 21,946 & 2 & 19,900 & 3 & 19,900 & 3 \\
miniDAFJS26 & 22,306 &  6 & 22,230 &  9 & 20,724 & 5 & 20,724 & 7 & 19,333 & 5 & 19,333 & 7 \\
miniDAFJS27 & 27,086 &  3 & 27,086 &  3 & 23,865 & 6 & 23,145 & 9 & 22,513 & 3 & 22,334 & 4 \\
miniDAFJS28 & 23,841 & 10 & 23,841 & 11 & 22,745 & 7 & 22,745 & 7 & 20,020 & 2 & 20,020 & 2 \\
miniDAFJS29 & 21,151 &  1 & 21,057 &  2 & 19,789 & 1 & 19,789 & 1 & 18,560 & 1 & 18,414 & 3 \\
miniDAFJS30 & 26,428 &  2 & 26,126 &  5 & 23,736 & 2 & 22,973 & 5 & 21,356 & 2 & 20,760 & 5 \\ 
\hline
mean & 23,405.60 & 2.03 & 23,176.57 & 3.30 & 21,897.63 & 2.17 & 21,519.20 & 3.50 & 20,395.93 & 2.33 & 20,141.67 & 3.47 \\ 
wins & 12 & & 30 & & 10 & & 30 & & 12 & & 19 & \\
\hline
\multicolumn{13}{c}{}\\
\hline
\multirow{3}{*}{instance} & \multicolumn{4}{c|}{$ \alpha = 0.1$} & \multicolumn{4}{c|}{$ \alpha = 0.2$} & \multicolumn{4}{c|}{$ \alpha = 0.3$} \\ 
\cline{2-13}
& \multicolumn{2}{c|}{Cropped} & \multicolumn{2}{c|}{Reduced} & \multicolumn{2}{c|}{Cropped} & \multicolumn{2}{c|}{Reduced} & \multicolumn{2}{c|}{Cropped} & \multicolumn{2}{c|}{Reduced} \\ 
\cline{2-13}
& $C_{\max}$ & \#it & $C_{\max}$ & \#it & $C_{\max}$ & \#it & $C_{\max}$ & \#it & $C_{\max}$ & \#it & $C_{\max}$ & \#it \\ 
\hline
\hline
miniYFJS01 & 35,243 & 1 & 35,046 &  2 & 34,443 & 1 & 33,132 &  4 & 33,697 &  1 & 32,199 &  2 \\
miniYFJS02 & 28,688 & 1 & 28,688 &  1 & 27,557 & 1 & 27,557 &  1 & 25,969 &  1 & 25,969 &  1 \\
miniYFJS03 & 48,141 & 3 & 47,391 &  6 & 44,195 & 3 & 42,896 &  6 & 39,880 &  8 & 38,935 & 11 \\
miniYFJS04 & 25,394 & 1 & 25,394 &  1 & 24,669 & 1 & 24,485 &  2 & 24,017 &  1 & 23,774 &  2 \\
miniYFJS05 & 26,371 & 2 & 25,750 &  3 & 25,635 & 3 & 24,235 &  3 & 24,743 &  1 & 24,743 &  1 \\
miniYFJS06 & 30,952 & 1 & 30,952 &  1 & 29,080 & 1 & 29,080 &  1 & 27,366 &  1 & 27,366 &  1 \\
miniYFJS07 & 46,782 & 2 & 45,705 &  6 & 43,886 & 3 & 42,578 &  5 & 41,487 &  3 & 39,694 &  5 \\
miniYFJS08 & 33,954 & 2 & 33,883 &  7 & 31,992 & 2 & 31,773 &  7 & 30,276 &  2 & 29,990 &  6 \\
miniYFJS09 & 37,049 & 2 & 37,049 &  2 & 36,098 & 2 & 36,098 &  2 & 34,357 &  1 & 34,357 &  1 \\
miniYFJS10 & 29,416 & 6 & 29,416 &  6 & 30,290 & 3 & 29,858 &  2 & 27,547 &  2 & 27,425 &  2 \\
miniYFJS11 & 51,212 & 1 & 51,129 &  2 & 47,079 & 1 & 46,939 &  2 & 43,254 &  1 & 43,254 &  1 \\
miniYFJS12 & 36,343 & 4 & 36,343 &  6 & 33,404 & 5 & 33,404 &  5 & 29,989 &  5 & 29,989 &  5 \\
miniYFJS13 & 32,915 & 7 & 32,792 &  8 & 30,219 & 5 & 30,219 &  7 & 25,881 &  3 & 25,619 &  4 \\
miniYFJS14 & 31,826 & 7 & 31,826 &  9 & 31,284 & 4 & 31,284 &  4 & 28,131 &  6 & 28,131 &  6 \\
miniYFJS15 & 45,901 & 4 & 45,442 &  6 & 42,828 & 4 & 45,903 &  5 & 39,862 &  4 & 46,387 &  2 \\
miniYFJS16 & 33,965 & 3 & 33,791 &  7 & 32,481 & 3 & 32,165 &  7 & 30,483 &  2 & 30,483 &  2 \\
miniYFJS17 & 44,181 & 5 & 52,936 &  3 & 41,316 & 5 & 48,475 &  4 & 38,739 &  5 & 44,650 &  3 \\
miniYFJS18 & 34,133 & 4 & 34,044 &  9 & 29,883 & 3 & 30,201 &  5 & 27,944 &  3 & 29,698 &  2 \\
miniYFJS19 & 39,165 & 4 & 36,706 & 14 & 33,965 & 6 & 33,805 & 12 & 31,743 &  6 & 31,743 &  9 \\
miniYFJS20 & 36,071 & 9 & 36,071 & 11 & 35,100 & 4 & 34,689 &  8 & 32,603 &  4 & 32,718 &  5 \\
miniYFJS21 & 40,994 & 7 & 39,978 & 10 & 37,414 & 5 & 35,570 & 16 & 33,355 &  5 & 33,355 &  7 \\
miniYFJS22 & 35,319 & 7 & 34,282 & 15 & 32,485 & 1 & 32,337 &  2 & 30,212 &  1 & 29,512 &  4 \\
miniYFJS23 & 45,713 & 9 & 45,725 & 10 & 42,016 & 3 & 42,016 &  7 & 38,777 &  4 & 37,963 &  9 \\
miniYFJS24 & 40,291 & 6 & 36,367 & 14 & 38,144 & 2 & 38,144 &  4 & 30,173 & 11 & 31,054 & 14 \\
miniYFJS25 & 43,592 & 3 & 43,555 &  4 & 40,933 & 1 & 40,425 &  3 & 34,388 &  6 & 34,388 &  5 \\
miniYFJS26 & 54,739 & 2 & 54,675 &  3 & 49,412 & 1 & 49,412 &  1 & 44,766 &  8 & 43,452 & 18 \\
miniYFJS27 & 36,899 & 6 & 36,765 & 18 & 34,359 & 7 & 34,943 &  7 & 32,539 &  3 & 31,749 &  7 \\
miniYFJS28 & 41,951 & 2 & 41,838 &  3 & 38,640 & 2 & 33,096 &  8 & 35,729 &  2 & 35,729 &  7 \\
miniYFJS29 & 44,930 & 9 & 47,807 &  4 & 38,993 & 3 & 36,561 &  9 & 35,726 &  3 & 38,122 &  5 \\
miniYFJS30 & 44,038 & 3 & 43,725 &  8 & 41,290 & 3 & 41,640 &  5 & 32,238 &  6 & 33,215 & 13 \\ 
\hline
mean & 38,538.93 & 4.10 & 38,502.37 & 6.63 & 35,969.67 & 2.93 & 35,764.00 & 5.13 & 32,862.37 & 3.63 & 33,188.77 & 5.33 \\
wins & 11 & & 27 & & 14 & & 25 & & 19 & & 13 & \\
\hline
\end{tabular}}
\end{center}
\caption{Local search results with the first-improvement strategy
  using the reduced neighborhood and the cropped neighborhood, applied
  to the small-sized instances with learning rate $\alpha \in \{0.1,
  0.2, 0.3\}$.}
\label{tab1}
\end{table}

\begin{table}[ht!]
\begin{center}  
\resizebox{\textwidth}{!}{
\begin{tabular}{|cc|ccc|ccc|ccc|ccc|ccc|ccc|} 
\hline
\multirow{3}{*}{instance} &  & \multicolumn{6}{c|}{$ \alpha = 0.1$} & \multicolumn{6}{c|}{$ \alpha = 0.2$} & \multicolumn{6}{c|}{$ \alpha = 0.3$} \\ 
\cline{3-20}
&  & \multicolumn{3}{c|}{Cropped} & \multicolumn{3}{c|}{Reduced} & \multicolumn{3}{c|}{Cropped} & \multicolumn{3}{c|}{Reduced} & \multicolumn{3}{c|}{Cropped} & \multicolumn{3}{c|}{Reduced} \\ 
\cline{3-20}
&  & $C_{\max}$ & \#it & Time & $C_{\max}$ & \#it & Time & $C_{\max}$ & \#it & Time & $C_{\max}$ & \#it & Time & $C_{\max}$ & \#it & Time & $C_{\max}$ & \#it & Time \\ 
\hline
\hline
DAFJS01 &  & 28,036 & 4 & 0.001 & 26,794 & 14 & 0.004 & 26,371 & 3 & 0.000 & 25,970 & 19 & 0.001 & 21,685 & 3 & 0.000 & 21,439 & 5 & 0.000 \\
DAFJS02 &  & 28,098 & 9 & 0.002 & 28,098 & 10 & 0.001 & 28,141 & 2 & 0.000 & 27,860 & 5 & 0.000 & 23,672 & 8 & 0.000 & 23,808 & 7 & 0.000 \\
DAFJS03 &  & 53,688 & 1 & 0.001 & 51,579 & 15 & 0.005 & 47,879 & 3 & 0.000 & 45,458 & 37 & 0.011 & 43,509 & 3 & 0.000 & 42,867 & 13 & 0.005 \\
DAFJS04 &  & 54,082 & 1 & 0.001 & 52,563 & 15 & 0.001 & 47,778 & 5 & 0.000 & 46,577 & 8 & 0.001 & 41,899 & 2 & 0.000 & 41,697 & 3 & 0.000 \\
DAFJS05 &  & 40,195 & 25 & 0.004 & 41,604 & 24 & 0.002 & 39,771 & 16 & 0.001 & 37,797 & 36 & 0.006 & 35,610 & 9 & 0.000 & 31,254 & 47 & 0.009 \\
DAFJS06 &  & 41,271 & 37 & 0.006 & 50,586 & 10 & 0.002 & 38,087 & 20 & 0.001 & 36,252 & 69 & 0.014 & 32,168 & 29 & 0.002 & 31,122 & 36 & 0.006 \\
DAFJS07 &  & 54,989 & 7 & 0.004 & 54,802 & 36 & 0.026 & 45,737 & 34 & 0.006 & 48,959 & 29 & 0.021 & 40,830 & 6 & 0.001 & 39,662 & 24 & 0.024 \\
DAFJS08 &  & 58,930 & 13 & 0.006 & 60,359 & 27 & 0.026 & 52,816 & 5 & 0.001 & 50,928 & 28 & 0.022 & 44,590 & 4 & 0.001 & 42,777 & 43 & 0.038 \\
DAFJS09 &  & 45,139 & 11 & 0.001 & 44,917 & 17 & 0.003 & 36,871 & 29 & 0.002 & 36,735 & 39 & 0.009 & 34,185 & 19 & 0.002 & 34,025 & 39 & 0.011 \\
DAFJS10 &  & 53,977 & 6 & 0.001 & 51,565 & 50 & 0.011 & 45,971 & 4 & 0.000 & 45,560 & 12 & 0.003 & 36,279 & 8 & 0.001 & 34,289 & 42 & 0.010 \\
DAFJS11 &  & 66,872 & 6 & 0.003 & 63,878 & 55 & 0.052 & 51,746 & 37 & 0.016 & 54,661 & 30 & 0.049 & 45,657 & 17 & 0.006 & 44,809 & 46 & 0.068 \\
DAFJS12 &  & 62,198 & 31 & 0.014 & 59,612 & 157 & 0.155 & 53,738 & 9 & 0.002 & 52,467 & 38 & 0.040 & 44,490 & 8 & 0.004 & 43,717 & 28 & 0.043 \\
DAFJS13 &  & 57,874 & 16 & 0.005 & 57,310 & 36 & 0.021 & 48,847 & 22 & 0.004 & 48,816 & 29 & 0.016 & 41,717 & 4 & 0.001 & 41,960 & 11 & 0.006 \\
DAFJS14 &  & 70,169 & 39 & 0.008 & 65,321 & 59 & 0.029 & 56,639 & 35 & 0.006 & 55,881 & 44 & 0.027 & 46,392 & 26 & 0.006 & 44,707 & 54 & 0.035 \\
DAFJS15 &  & 61,882 & 62 & 0.032 & 63,408 & 124 & 0.166 & 53,807 & 6 & 0.003 & 53,686 & 14 & 0.031 & 47,692 & 10 & 0.004 & 47,474 & 33 & 0.067 \\
DAFJS16 &  & 70,727 & 40 & 0.017 & 70,838 & 84 & 0.125 & 63,277 & 11 & 0.005 & 58,436 & 107 & 0.183 & 55,112 & 12 & 0.006 & 53,095 & 78 & 0.172 \\
DAFJS17 &  & 68,446 & 50 & 0.018 & 70,306 & 26 & 0.016 & 57,893 & 41 & 0.009 & 57,518 & 54 & 0.031 & 46,340 & 80 & 0.026 & 47,331 & 96 & 0.087 \\
DAFJS18 &  & 71,871 & 31 & 0.006 & 74,631 & 21 & 0.015 & 58,533 & 19 & 0.004 & 58,907 & 29 & 0.019 & 50,818 & 7 & 0.002 & 48,672 & 26 & 0.015 \\
DAFJS19 &  & 63,999 & 12 & 0.002 & 63,314 & 26 & 0.019 & 52,824 & 24 & 0.006 & 57,082 & 10 & 0.014 & 41,426 & 28 & 0.004 & 41,865 & 56 & 0.029 \\
DAFJS20 &  & 64,005 & 75 & 0.020 & 63,841 & 225 & 0.207 & 61,871 & 5 & 0.002 & 57,378 & 99 & 0.115 & 46,868 & 53 & 0.018 & 48,600 & 40 & 0.051 \\
DAFJS21 &  & 72,142 & 57 & 0.020 & 76,433 & 31 & 0.032 & 61,483 & 12 & 0.004 & 62,677 & 9 & 0.009 & 50,318 & 15 & 0.005 & 49,422 & 39 & 0.040 \\
DAFJS22 &  & 60,081 & 132 & 0.043 & 60,143 & 315 & 0.448 & 53,248 & 44 & 0.014 & 54,318 & 42 & 0.097 & 43,566 & 17 & 0.011 & 41,615 & 94 & 0.235 \\
DAFJS23 &  & 48,338 & 9 & 0.001 & 48,054 & 18 & 0.012 & 41,433 & 17 & 0.003 & 44,473 & 31 & 0.018 & 39,105 & 3 & 0.000 & 38,334 & 12 & 0.006 \\
DAFJS24 &  & 55,507 & 9 & 0.004 & 55,353 & 33 & 0.032 & 47,277 & 14 & 0.004 & 47,656 & 9 & 0.007 & 40,921 & 17 & 0.004 & 42,473 & 27 & 0.019 \\
DAFJS25 &  & 79,286 & 26 & 0.011 & 76,253 & 92 & 0.227 & 59,672 & 68 & 0.033 & 61,805 & 92 & 0.184 & 48,469 & 25 & 0.013 & 49,390 & 39 & 0.098 \\
DAFJS26 &  & 73,969 & 54 & 0.020 & 75,329 & 96 & 0.190 & 70,211 & 8 & 0.003 & 65,524 & 109 & 0.373 & 53,944 & 39 & 0.023 & 55,269 & 39 & 0.117 \\
DAFJS27 &  & 79,234 & 13 & 0.004 & 79,064 & 23 & 0.053 & 62,437 & 21 & 0.010 & 62,596 & 29 & 0.079 & 53,719 & 46 & 0.023 & 55,569 & 46 & 0.134 \\
DAFJS28 &  & 59,201 & 18 & 0.003 & 54,866 & 109 & 0.103 & 48,414 & 21 & 0.004 & 47,884 & 46 & 0.039 & 41,408 & 11 & 0.003 & 40,097 & 34 & 0.050 \\
DAFJS29 &  & 68,684 & 21 & 0.006 & 66,929 & 80 & 0.106 & 57,356 & 16 & 0.005 & 57,179 & 34 & 0.044 & 49,819 & 7 & 0.005 & 49,519 & 24 & 0.050 \\
DAFJS30 &  & 57,350 & 17 & 0.005 & 60,016 & 13 & 0.013 & 51,283 & 14 & 0.004 & 50,961 & 45 & 0.052 & 43,426 & 3 & 0.001 & 42,592 & 6 & 0.008 \\ 
\hline
mean &  & 59,008.00 & 27.73 & 0.009 & 58,925.53 & 61.37 & 0.070 & 50,713.70 & 18.83 & 0.005 & 50,400.03 & 39.40 & 0.051 & 42,854.47 & 17.30 & 0.006 & 42,315.00 & 36.23 & 0.048 \\ 
wins &  & 12 & & & 19 & & & 10 & & & 20 & & & 9 & & & 12 & & \\
\hline
\multicolumn{20}{c}{}\\
\hline
\multirow{3}{*}{instance} &  & \multicolumn{6}{c|}{$ \alpha = 0.1$} & \multicolumn{6}{c|}{$ \alpha = 0.2$} & \multicolumn{6}{c|}{$ \alpha = 0.3$} \\ 
\cline{3-20}
&  & \multicolumn{3}{c|}{Cropped} & \multicolumn{3}{c|}{Reduced} & \multicolumn{3}{c|}{Cropped} & \multicolumn{3}{c|}{Reduced} & \multicolumn{3}{c|}{Cropped} & \multicolumn{3}{c|}{Reduced} \\ 
\cline{3-20}
&  & $C_{\max}$ & \#it & Time & $C_{\max}$ & \#it & Time & $C_{\max}$ & \#it & Time & $C_{\max}$ & \#it & Time & $C_{\max}$ & \#it & Time & $C_{\max}$ & \#it & Time \\ 
\hline
\hline
YFJS01 &  & 79,662 & 10 & 0.000 & 79,737 & 21 & 0.003 & 74,547 & 9 & 0.000 & 70,624 & 31 & 0.004 & 66,246 & 9 & 0.000 & 67,588 & 10 & 0.001 \\
YFJS02 &  & 76,529 & 3 & 0.000 & 76,813 & 4 & 0.000 & 66,853 & 1 & 0.000 & 66,853 & 1 & 0.000 & 59,211 & 3 & 0.000 & 59,001 & 4 & 0.001 \\
YFJS03 &  & 38,424 & 3 & 0.000 & 35,936 & 14 & 0.000 & 33,273 & 5 & 0.000 & 33,273 & 5 & 0.000 & 30,995 & 3 & 0.000 & 31,407 & 2 & 0.000 \\
YFJS04 &  & 45,683 & 3 & 0.000 & 45,533 & 4 & 0.000 & 42,553 & 4 & 0.000 & 42,553 & 4 & 0.000 & 38,092 & 5 & 0.000 & 38,092 & 5 & 0.000 \\
YFJS05 &  & 45,421 & 4 & 0.000 & 45,421 & 4 & 0.000 & 40,738 & 2 & 0.000 & 40,738 & 2 & 0.000 & 38,324 & 5 & 0.000 & 38,324 & 9 & 0.000 \\
YFJS06 &  & 44,848 & 13 & 0.000 & 47,384 & 12 & 0.001 & 43,384 & 15 & 0.000 & 42,207 & 21 & 0.002 & 45,370 & 2 & 0.000 & 45,370 & 2 & 0.000 \\
YFJS07 &  & 48,980 & 7 & 0.000 & 45,722 & 21 & 0.001 & 43,925 & 12 & 0.000 & 45,971 & 11 & 0.000 & 38,453 & 5 & 0.000 & 38,453 & 5 & 0.000 \\
YFJS08 &  & 43,073 & 7 & 0.000 & 41,658 & 19 & 0.001 & 42,901 & 3 & 0.000 & 41,573 & 10 & 0.000 & 35,429 & 2 & 0.000 & 33,432 & 17 & 0.001 \\
YFJS09 &  & 25,067 & 3 & 0.000 & 24,350 & 4 & 0.000 & 23,059 & 3 & 0.000 & 22,211 & 5 & 0.000 & 20,916 & 7 & 0.000 & 20,916 & 7 & 0.000 \\
YFJS10 &  & 42,583 & 1 & 0.000 & 42,278 & 2 & 0.000 & 39,126 & 5 & 0.000 & 38,683 & 7 & 0.000 & 34,487 & 8 & 0.000 & 36,415 & 7 & 0.000 \\
YFJS11 &  & 55,765 & 4 & 0.000 & 55,416 & 7 & 0.001 & 49,781 & 2 & 0.000 & 49,160 & 6 & 0.001 & 45,051 & 1 & 0.000 & 44,455 & 4 & 0.000 \\
YFJS12 &  & 56,674 & 24 & 0.001 & 54,852 & 59 & 0.007 & 55,174 & 11 & 0.000 & 53,329 & 21 & 0.003 & 47,092 & 3 & 0.000 & 44,873 & 12 & 0.001 \\
YFJS13 &  & 44,083 & 16 & 0.000 & 45,844 & 16 & 0.001 & 42,442 & 10 & 0.000 & 39,720 & 27 & 0.003 & 35,776 & 16 & 0.000 & 35,524 & 19 & 0.003 \\
YFJS14 &  & 119,977 & 27 & 0.016 & 124,518 & 60 & 0.141 & 104,506 & 19 & 0.010 & 105,770 & 58 & 0.181 & 89,687 & 12 & 0.007 & 87,368 & 61 & 0.182 \\
YFJS15 &  & 126,713 & 33 & 0.016 & 127,073 & 171 & 0.481 & 106,511 & 28 & 0.016 & 108,412 & 58 & 0.192 & 93,075 & 14 & 0.007 & 86,141 & 166 & 0.529 \\
YFJS16 &  & 116,201 & 37 & 0.018 & 121,141 & 73 & 0.196 & 101,164 & 13 & 0.005 & 99,801 & 65 & 0.129 & 87,351 & 10 & 0.005 & 87,648 & 48 & 0.136 \\
YFJS17 &  & 104,688 & 21 & 0.026 & 102,614 & 159 & 1.134 & 83,355 & 18 & 0.014 & 82,211 & 87 & 0.621 & 71,201 & 26 & 0.036 & 71,075 & 65 & 0.493 \\
YFJS18 &  & 126,647 & 34 & 0.032 & 128,529 & 87 & 0.752 & 97,157 & 15 & 0.022 & 94,689 & 179 & 2.000 & 84,518 & 9 & 0.015 & 82,967 & 77 & 0.834 \\
YFJS19 &  & 93,306 & 96 & 0.100 & 100,855 & 209 & 1.353 & 86,398 & 21 & 0.026 & 85,028 & 128 & 1.115 & 68,792 & 53 & 0.054 & 68,955 & 135 & 0.984 \\
YFJS20 &  & 96,248 & 17 & 0.019 & 92,240 & 196 & 1.205 & 79,408 & 85 & 0.088 & 80,875 & 294 & 1.923 & 71,676 & 8 & 0.009 & 68,562 & 123 & 0.993 \\ 
\hline
mean &  & 71,528.60 & 18.15 & 0.011 & 71,895.70 & 57.10 & 0.264 & 62,812.75 & 14.05 & 0.009 & 62,184.05 & 51.00 & 0.309 & 55,087.10 & 10.05 & 0.007 & 54,328.30 & 38.90 & 0.208 \\
wins &  & 10 & & & 11 & & & 8 & & & 16 & & & 10 & & & 15 & & \\
\hline
\end{tabular}}
\end{center}
\caption{Local search results with the first-improvement strategy
  using the reduced neighborhood and the cropped neighborhood, applied
  to the instances proposed in~\cite{birgin2014milp} with learning
  rate $\alpha \in \{0.1, 0.2, 0.3\}$.}
\label{tab2}
\end{table}

\begin{table}[ht!]
\begin{center}  
\resizebox{!}{0.45\textheight}{
\begin{tabular}{|c|cc|cc|cc|cc|cc|cc|} 
\hline
\multirow{3}{*}{instance} & \multicolumn{4}{c|}{$ \alpha = 0.1$} & \multicolumn{4}{c|}{$ \alpha = 0.2$} & \multicolumn{4}{c|}{$ \alpha = 0.3$} \\ 
\cline{2-13}
& \multicolumn{2}{c|}{Cropped} & \multicolumn{2}{c|}{Reduced} & \multicolumn{2}{c|}{Cropped} & \multicolumn{2}{c|}{Reduced} & \multicolumn{2}{c|}{Cropped} & \multicolumn{2}{c|}{Reduced} \\ 
\cline{2-13}
& $C_{\max}$ & \#it & $C_{\max}$ & \#it & $C_{\max}$ & \#it & $C_{\max}$ & \#it & $C_{\max}$ & \#it & $C_{\max}$ & \#it \\ 
\hline
\hline
miniDAFJS01 & 23,264 & 1 & 23,019 & 2 & 21,588 & 3 & 21,546 & 4 & 20,007 & 4 & 19,848 & 4 \\
miniDAFJS02 & 23,242 & 1 & 22,927 & 3 & 22,161 & 1 & 22,046 & 2 & 21,136 & 2 & 20,987 & 2 \\
miniDAFJS03 & 18,363 & 1 & 18,363 & 1 & 17,972 & 1 & 17,972 & 1 & 17,619 & 2 & 17,619 & 2 \\
miniDAFJS04 & 21,690 & 1 & 21,151 & 4 & 20,691 & 2 & 19,602 & 4 & 19,867 & 2 & 18,800 & 4 \\
miniDAFJS05 & 22,598 & 1 & 22,418 & 2 & 20,826 & 1 & 20,523 & 2 & 19,253 & 1 & 18,878 & 2 \\
miniDAFJS06 & 23,370 & 1 & 23,370 & 1 & 20,783 & 2 & 20,783 & 2 & 19,605 & 2 & 19,605 & 2 \\
miniDAFJS07 & 28,644 & 1 & 28,644 & 1 & 25,088 & 1 & 24,715 & 2 & 24,636 & 2 & 24,256 & 2 \\
miniDAFJS08 & 19,878 & 1 & 19,878 & 1 & 18,857 & 1 & 18,857 & 1 & 17,900 & 2 & 17,900 & 2 \\
miniDAFJS09 & 25,425 & 1 & 25,425 & 1 & 23,830 & 1 & 23,715 & 2 & 22,398 & 1 & 22,258 & 2 \\
miniDAFJS10 & 21,563 & 1 & 20,359 & 4 & 21,021 & 1 & 18,823 & 5 & 19,723 & 2 & 17,466 & 5 \\
miniDAFJS11 & 33,689 & 2 & 33,585 & 4 & 30,941 & 2 & 30,550 & 5 & 28,510 & 2 & 28,342 & 2 \\
miniDAFJS12 & 20,342 & 1 & 20,342 & 1 & 19,624 & 1 & 19,624 & 1 & 19,094 & 1 & 19,094 & 1 \\
miniDAFJS13 & 17,091 & 2 & 16,313 & 3 & 16,616 & 2 & 15,143 & 3 & 15,278 & 2 & 15,053 & 2 \\
miniDAFJS14 & 23,248 & 4 & 23,140 & 5 & 21,990 & 5 & 21,817 & 5 & 20,711 & 5 & 20,620 & 5 \\
miniDAFJS15 & 22,472 & 1 & 22,472 & 1 & 20,620 & 3 & 20,620 & 3 & 18,628 & 3 & 18,628 & 3 \\
miniDAFJS16 & 25,691 & 1 & 25,426 & 2 & 24,593 & 1 & 24,114 & 2 & 22,939 & 3 & 22,939 & 3 \\
miniDAFJS17 & 21,070 & 1 & 20,155 & 3 & 20,183 & 2 & 20,183 & 2 & 19,139 & 2 & 19,139 & 2 \\
miniDAFJS18 & 18,445 & 2 & 18,135 & 4 & 18,829 & 1 & 18,201 & 4 & 17,784 & 1 & 17,295 & 2 \\
miniDAFJS19 & 21,293 & 1 & 20,945 & 3 & 20,107 & 1 & 19,642 & 3 & 19,030 & 1 & 18,474 & 3 \\
miniDAFJS20 & 23,443 & 2 & 23,212 & 3 & 21,286 & 1 & 21,286 & 1 & 19,587 & 1 & 19,587 & 1 \\
miniDAFJS21 & 24,404 & 3 & 24,274 & 4 & 22,215 & 2 & 22,151 & 4 & 21,182 & 3 & 20,368 & 4 \\
miniDAFJS22 & 25,923 & 3 & 25,923 & 3 & 24,273 & 3 & 24,273 & 3 & 22,767 & 3 & 22,767 & 3 \\
miniDAFJS23 & 25,839 & 2 & 25,788 & 3 & 24,253 & 2 & 24,164 & 3 & 21,772 & 5 & 22,710 & 3 \\
miniDAFJS24 & 26,932 & 1 & 26,610 & 4 & 24,500 & 1 & 23,709 & 3 & 22,524 & 1 & 21,714 & 4 \\
miniDAFJS25 & 23,370 & 1 & 22,982 & 2 & 22,613 & 1 & 21,946 & 2 & 18,472 & 6 & 20,834 & 3 \\
miniDAFJS26 & 22,306 & 3 & 22,230 & 4 & 20,273 & 4 & 20,273 & 4 & 18,921 & 4 & 18,921 & 4 \\
miniDAFJS27 & 27,086 & 3 & 27,086 & 3 & 23,145 & 5 & 23,145 & 5 & 22,513 & 3 & 22,334 & 4 \\
miniDAFJS28 & 23,841 & 6 & 23,841 & 6 & 21,927 & 7 & 21,927 & 7 & 20,020 & 2 & 20,020 & 2 \\
miniDAFJS29 & 21,151 & 1 & 21,057 & 2 & 19,789 & 1 & 19,789 & 1 & 18,560 & 1 & 18,414 & 2 \\
miniDAFJS30 & 26,428 & 2 & 26,126 & 4 & 23,736 & 2 & 22,973 & 4 & 21,356 & 2 & 20,760 & 3 \\ 
\hline
mean & 23,403.37 & 1.73 & 23,173.20 & 2.80 & 21,811.00 & 2.03 & 21,470.40 & 3.00 & 20,364.37 & 2.37 & 20,187.67 & 2.77 \\ 
wins & 12 &  & 30 &  & 12 &  & 30 &  & 13 &  & 18 & \\
\hline
\multicolumn{13}{c}{}\\
\hline
\multirow{3}{*}{instance} & \multicolumn{4}{c|}{$ \alpha = 0.1$} & \multicolumn{4}{c|}{$ \alpha = 0.2$} & \multicolumn{4}{c|}{$ \alpha = 0.3$} \\ 
\cline{2-13}
& \multicolumn{2}{c|}{Cropped} & \multicolumn{2}{c|}{Reduced} & \multicolumn{2}{c|}{Cropped} & \multicolumn{2}{c|}{Reduced} & \multicolumn{2}{c|}{Cropped} & \multicolumn{2}{c|}{Reduced} \\ 
\cline{2-13}
& $C_{\max}$ & \#it & $C_{\max}$ & \#it & $C_{\max}$ & \#it & $C_{\max}$ & \#it & $C_{\max}$ & \#it & $C_{\max}$ & \#it \\ 
\hline
\hline
miniYFJS01 & 35,243 & 1 & 35,046 & 2 & 34,443 & 1 & 33,132 & 4 & 33,697 & 1 & 32,199 & 2 \\
miniYFJS02 & 28,688 & 1 & 28,688 & 1 & 27,557 & 1 & 27,557 & 1 & 25,969 & 1 & 25,969 & 1 \\
miniYFJS03 & 53,381 & 2 & 52,098 & 3 & 49,111 & 2 & 46,806 & 3 & 42,787 & 3 & 42,787 & 3 \\
miniYFJS04 & 25,394 & 1 & 25,394 & 1 & 24,669 & 1 & 24,485 & 2 & 24,017 & 1 & 23,774 & 2 \\
miniYFJS05 & 26,371 & 2 & 25,750 & 3 & 25,635 & 3 & 24,235 & 3 & 24,743 & 1 & 24,743 & 1 \\
miniYFJS06 & 30,952 & 1 & 30,952 & 1 & 29,080 & 1 & 29,080 & 1 & 27,366 & 1 & 27,366 & 1 \\
miniYFJS07 & 46,782 & 2 & 46,089 & 4 & 44,137 & 2 & 42,578 & 5 & 41,737 & 2 & 39,229 & 5 \\
miniYFJS08 & 33,954 & 2 & 33,883 & 3 & 31,992 & 2 & 31,597 & 4 & 30,276 & 2 & 29,631 & 4 \\
miniYFJS09 & 37,049 & 2 & 37,049 & 2 & 36,098 & 2 & 36,098 & 2 & 34,357 & 1 & 34,357 & 1 \\
miniYFJS10 & 29,416 & 5 & 29,416 & 5 & 27,604 & 4 & 29,858 & 2 & 27,547 & 2 & 27,425 & 2 \\
miniYFJS11 & 51,212 & 1 & 51,129 & 2 & 47,079 & 1 & 46,939 & 2 & 43,254 & 1 & 43,254 & 1 \\
miniYFJS12 & 36,343 & 4 & 36,343 & 4 & 33,404 & 5 & 33,404 & 5 & 29,989 & 5 & 29,989 & 5 \\
miniYFJS13 & 32,915 & 5 & 32,792 & 5 & 30,219 & 4 & 30,219 & 4 & 25,881 & 3 & 25,619 & 4 \\
miniYFJS14 & 31,826 & 6 & 31,826 & 6 & 31,284 & 4 & 31,284 & 4 & 28,131 & 6 & 28,131 & 6 \\
miniYFJS15 & 45,597 & 3 & 45,442 & 4 & 42,306 & 3 & 41,893 & 4 & 39,187 & 3 & 38,619 & 4 \\
miniYFJS16 & 33,965 & 3 & 33,791 & 7 & 32,481 & 3 & 32,165 & 6 & 30,483 & 2 & 30,483 & 2 \\
miniYFJS17 & 44,181 & 5 & 44,044 & 7 & 41,316 & 5 & 41,199 & 5 & 38,739 & 5 & 38,584 & 5 \\
miniYFJS18 & 34,133 & 3 & 34,044 & 4 & 29,883 & 2 & 29,883 & 2 & 27,944 & 2 & 27,944 & 2 \\
miniYFJS19 & 39,165 & 3 & 39,046 & 4 & 33,965 & 6 & 33,805 & 7 & 31,743 & 6 & 31,082 & 8 \\
miniYFJS20 & 30,837 & 6 & 30,837 & 6 & 31,304 & 6 & 31,304 & 6 & 29,599 & 4 & 29,599 & 4 \\
miniYFJS21 & 40,523 & 6 & 40,523 & 6 & 35,570 & 9 & 35,570 & 9 & 33,355 & 5 & 33,355 & 5 \\
miniYFJS22 & 35,019 & 4 & 35,019 & 4 & 32,485 & 1 & 31,946 & 2 & 30,212 & 1 & 29,498 & 2 \\
miniYFJS23 & 45,725 & 4 & 45,725 & 4 & 42,446 & 2 & 42,016 & 3 & 38,777 & 4 & 37,963 & 5 \\
miniYFJS24 & 40,291 & 5 & 40,291 & 6 & 38,144 & 2 & 38,144 & 2 & 30,173 & 7 & 31,572 & 9 \\
miniYFJS25 & 43,592 & 3 & 43,555 & 3 & 40,933 & 1 & 40,425 & 3 & 34,388 & 3 & 34,388 & 3 \\
miniYFJS26 & 54,739 & 2 & 54,675 & 3 & 49,412 & 1 & 49,412 & 1 & 44,633 & 3 & 44,207 & 6 \\
miniYFJS27 & 36,765 & 4 & 36,765 & 4 & 34,125 & 4 & 34,225 & 6 & 31,749 & 4 & 31,571 & 6 \\
miniYFJS28 & 41,951 & 2 & 41,838 & 3 & 38,640 & 2 & 38,450 & 3 & 35,729 & 2 & 35,729 & 2 \\
miniYFJS29 & 44,930 & 5 & 44,930 & 5 & 38,993 & 3 & 38,586 & 4 & 34,330 & 4 & 34,604 & 3 \\
miniYFJS30 & 44,038 & 3 & 43,725 & 4 & 41,290 & 3 & 40,744 & 4 & 31,093 & 7 & 31,093 & 7 \\ 
\hline
mean & 38,499.23 & 3.20 & 38,356.83 & 3.87 & 35,853.50 & 2.87 & 35,567.97 & 3.63 & 32,729.50 & 3.07 & 32,492.13 & 3.70 \\
wins & 14 & & 30 & & 13 & & 28 & & 17 & & 18 & \\
\hline
\end{tabular}}
\end{center}
\caption{Local search results with the best-improvement strategy
  using the reduced neighborhood and the cropped neighborhood, applied
  to the small-sized instances with learning rate $\alpha \in \{0.1,
  0.2, 0.3\}$.}
\label{tab3}
\end{table}

\begin{table}[ht!]
\begin{center}  
\resizebox{\textwidth}{!}{
\begin{tabular}{|cc|ccc|ccc|ccc|ccc|ccc|ccc|} 
\hline
\multirow{3}{*}{instance} &  & \multicolumn{6}{c|}{$ \alpha = 0.1$} & \multicolumn{6}{c|}{$ \alpha = 0.2$} & \multicolumn{6}{c|}{$ \alpha = 0.3$} \\ 
\cline{3-20}
&  & \multicolumn{3}{c|}{Cropped} & \multicolumn{3}{c|}{Reduced} & \multicolumn{3}{c|}{Cropped} & \multicolumn{3}{c|}{Reduced} & \multicolumn{3}{c|}{Cropped} & \multicolumn{3}{c|}{Reduced} \\ 
\cline{3-20}
&  & $C_{\max}$ & \#it & Time & $C_{\max}$ & \#it & Time & $C_{\max}$ & \#it & Time & $C_{\max}$ & \#it & Time & $C_{\max}$ & \#it & Time & $C_{\max}$ & \#it & Time \\ 
\hline
\hline
DAFJS01 &  & 26,465 & 4 & 0.000 & 26,203 & 7 & 0.004 & 26,371 & 3 & 0.000 & 24,993 & 8 & 0.001 & 21,685 & 3 & 0.000 & 21,685 & 3 & 0.000 \\
DAFJS02 &  & 28,098 & 3 & 0.000 & 28,098 & 3 & 0.001 & 27,949 & 2 & 0.000 & 27,727 & 4 & 0.000 & 23,021 & 3 & 0.000 & 23,021 & 3 & 0.000 \\
DAFJS03 &  & 53,688 & 1 & 0.000 & 51,580 & 14 & 0.009 & 47,879 & 2 & 0.000 & 45,916 & 15 & 0.011 & 43,509 & 2 & 0.000 & 42,258 & 8 & 0.006 \\
DAFJS04 &  & 54,082 & 1 & 0.000 & 52,494 & 11 & 0.004 & 47,571 & 4 & 0.000 & 46,502 & 7 & 0.003 & 41,899 & 2 & 0.000 & 41,533 & 4 & 0.001 \\
DAFJS05 &  & 40,723 & 10 & 0.006 & 40,723 & 12 & 0.005 & 37,118 & 15 & 0.002 & 38,484 & 12 & 0.005 & 33,331 & 12 & 0.001 & 33,813 & 7 & 0.003 \\
DAFJS06 &  & 45,234 & 11 & 0.007 & 44,877 & 15 & 0.009 & 41,995 & 6 & 0.001 & 41,909 & 9 & 0.005 & 33,956 & 5 & 0.001 & 30,556 & 20 & 0.011 \\
DAFJS07 &  & 54,937 & 5 & 0.006 & 53,651 & 15 & 0.032 & 44,414 & 16 & 0.009 & 44,147 & 21 & 0.046 & 40,442 & 5 & 0.002 & 40,057 & 13 & 0.027 \\
DAFJS08 &  & 58,941 & 10 & 0.009 & 56,561 & 14 & 0.031 & 52,803 & 3 & 0.001 & 51,128 & 16 & 0.038 & 44,590 & 4 & 0.001 & 43,737 & 16 & 0.034 \\
DAFJS09 &  & 43,032 & 9 & 0.003 & 45,796 & 5 & 0.002 & 38,062 & 12 & 0.002 & 37,876 & 12 & 0.007 & 36,488 & 5 & 0.001 & 36,303 & 10 & 0.008 \\
DAFJS10 &  & 54,858 & 5 & 0.002 & 51,277 & 24 & 0.025 & 45,971 & 4 & 0.001 & 45,868 & 3 & 0.002 & 35,273 & 8 & 0.002 & 35,023 & 9 & 0.009 \\
DAFJS11 &  & 63,795 & 8 & 0.010 & 64,796 & 20 & 0.066 & 53,770 & 7 & 0.007 & 52,845 & 17 & 0.061 & 44,701 & 7 & 0.005 & 44,599 & 11 & 0.035 \\
DAFJS12 &  & 62,433 & 15 & 0.024 & 61,533 & 28 & 0.128 & 54,073 & 5 & 0.007 & 52,355 & 17 & 0.080 & 44,051 & 5 & 0.004 & 43,449 & 11 & 0.045 \\
DAFJS13 &  & 58,902 & 6 & 0.004 & 56,453 & 22 & 0.027 & 48,838 & 13 & 0.007 & 49,179 & 10 & 0.011 & 41,717 & 4 & 0.001 & 41,717 & 4 & 0.005 \\
DAFJS14 &  & 62,637 & 42 & 0.029 & 70,357 & 8 & 0.014 & 55,503 & 19 & 0.012 & 54,609 & 25 & 0.040 & 44,777 & 16 & 0.008 & 44,646 & 23 & 0.042 \\
DAFJS15 &  & 64,527 & 21 & 0.032 & 62,030 & 39 & 0.205 & 52,115 & 9 & 0.011 & 53,677 & 7 & 0.035 & 46,921 & 12 & 0.012 & 47,268 & 15 & 0.071 \\
DAFJS16 &  & 69,300 & 21 & 0.029 & 66,153 & 71 & 0.358 & 61,676 & 8 & 0.011 & 60,010 & 20 & 0.093 & 54,072 & 15 & 0.015 & 49,031 & 46 & 0.260 \\
DAFJS17 &  & 70,803 & 9 & 0.011 & 69,975 & 22 & 0.054 & 55,294 & 6 & 0.006 & 54,456 & 18 & 0.045 & 45,667 & 17 & 0.016 & 47,550 & 28 & 0.103 \\
DAFJS18 &  & 67,984 & 27 & 0.022 & 71,052 & 10 & 0.031 & 57,965 & 16 & 0.012 & 57,274 & 25 & 0.042 & 48,351 & 7 & 0.004 & 49,753 & 5 & 0.013 \\
DAFJS19 &  & 57,786 & 16 & 0.014 & 57,471 & 22 & 0.046 & 51,065 & 8 & 0.004 & 51,065 & 8 & 0.019 & 43,998 & 9 & 0.003 & 40,678 & 26 & 0.053 \\
DAFJS20 &  & 67,045 & 22 & 0.021 & 66,674 & 30 & 0.115 & 57,271 & 27 & 0.027 & 55,355 & 47 & 0.180 & 45,160 & 23 & 0.024 & 49,481 & 14 & 0.055 \\
DAFJS21 &  & 72,979 & 17 & 0.025 & 71,913 & 23 & 0.103 & 58,044 & 14 & 0.018 & 57,787 & 20 & 0.070 & 48,693 & 16 & 0.023 & 48,912 & 16 & 0.064 \\
DAFJS22 &  & 60,238 & 46 & 0.068 & 62,189 & 29 & 0.196 & 54,147 & 11 & 0.017 & 53,728 & 15 & 0.088 & 41,480 & 28 & 0.040 & 42,035 & 25 & 0.156 \\
DAFJS23 &  & 47,719 & 9 & 0.004 & 47,520 & 15 & 0.022 & 42,480 & 8 & 0.003 & 42,943 & 11 & 0.014 & 39,105 & 3 & 0.001 & 38,185 & 8 & 0.011 \\
DAFJS24 &  & 54,427 & 13 & 0.009 & 55,115 & 12 & 0.033 & 46,797 & 15 & 0.012 & 46,985 & 12 & 0.024 & 41,832 & 4 & 0.005 & 39,979 & 21 & 0.063 \\
DAFJS25 &  & 74,359 & 37 & 0.049 & 73,353 & 48 & 0.301 & 58,767 & 31 & 0.043 & 61,465 & 29 & 0.165 & 49,919 & 9 & 0.013 & 49,527 & 14 & 0.059 \\
DAFJS26 &  & 68,577 & 50 & 0.076 & 76,328 & 12 & 0.080 & 63,243 & 36 & 0.049 & 68,469 & 18 & 0.142 & 51,544 & 16 & 0.021 & 50,829 & 31 & 0.219 \\
DAFJS27 &  & 79,254 & 2 & 0.004 & 78,967 & 6 & 0.036 & 62,520 & 12 & 0.020 & 62,359 & 18 & 0.120 & 51,382 & 29 & 0.041 & 53,010 & 23 & 0.168 \\
DAFJS28 &  & 58,564 & 4 & 0.004 & 53,505 & 41 & 0.100 & 48,593 & 12 & 0.012 & 48,941 & 9 & 0.029 & 41,206 & 4 & 0.004 & 40,923 & 7 & 0.019 \\
DAFJS29 &  & 70,450 & 6 & 0.005 & 65,758 & 36 & 0.120 & 56,387 & 20 & 0.014 & 57,708 & 14 & 0.054 & 49,906 & 5 & 0.005 & 48,412 & 18 & 0.074 \\
DAFJS30 &  & 57,326 & 6 & 0.004 & 56,839 & 13 & 0.032 & 49,649 & 17 & 0.012 & 49,312 & 30 & 0.092 & 42,819 & 2 & 0.002 & 42,592 & 4 & 0.011 \\ 
\hline
mean &  & 58,305.43 & 14.53 & 0.016 & 57,974.70 & 20.90 & 0.073 & 49,944.33 & 12.03 & 0.011 & 49,835.73 & 15.90 & 0.051 & 42,383.17 & 9.33 & 0.009 & 42,018.73 & 14.77 & 0.054 \\ 
wins &  & 9 & & & 23 & & & 10 & & & 21 & & & 11 & & & 13 & & \\
\hline
\multicolumn{20}{c}{}\\
\hline
\multirow{3}{*}{instance} &  & \multicolumn{6}{c|}{$ \alpha = 0.1$} & \multicolumn{6}{c|}{$ \alpha = 0.2$} & \multicolumn{6}{c|}{$ \alpha = 0.3$} \\ 
\cline{3-20}
&  & \multicolumn{3}{c|}{Cropped} & \multicolumn{3}{c|}{Reduced} & \multicolumn{3}{c|}{Cropped} & \multicolumn{3}{c|}{Reduced} & \multicolumn{3}{c|}{Cropped} & \multicolumn{3}{c|}{Reduced} \\ 
\cline{3-20}
&  & $C_{\max}$ & \#it & Time & $C_{\max}$ & \#it & Time & $C_{\max}$ & \#it & Time & $C_{\max}$ & \#it & Time & $C_{\max}$ & \#it & Time & $C_{\max}$ & \#it & Time \\ 
\hline
\hline
YFJS01 &  & 80,167 & 6 & 0.000 & 79,657 & 8 & 0.002 & 71,592 & 9 & 0.001 & 68,822 & 14 & 0.004 & 66,246 & 6 & 0.000 & 60,208 & 15 & 0.005 \\
YFJS02 &  & 76,529 & 3 & 0.000 & 76,529 & 3 & 0.000 & 66,853 & 1 & 0.000 & 66,853 & 1 & 0.000 & 59,001 & 2 & 0.000 & 59,001 & 2 & 0.000 \\
YFJS03 &  & 34,166 & 5 & 0.000 & 36,409 & 4 & 0.000 & 32,864 & 5 & 0.000 & 32,864 & 5 & 0.000 & 30,995 & 3 & 0.000 & 31,407 & 2 & 0.000 \\
YFJS04 &  & 45,683 & 3 & 0.000 & 45,533 & 4 & 0.000 & 42,553 & 4 & 0.000 & 42,553 & 4 & 0.000 & 38,092 & 3 & 0.000 & 38,092 & 3 & 0.000 \\
YFJS05 &  & 45,421 & 3 & 0.000 & 45,421 & 3 & 0.000 & 40,738 & 2 & 0.000 & 40,738 & 2 & 0.000 & 38,324 & 4 & 0.000 & 38,324 & 8 & 0.001 \\
YFJS06 &  & 44,936 & 10 & 0.001 & 44,936 & 10 & 0.002 & 40,272 & 9 & 0.000 & 43,744 & 6 & 0.001 & 43,123 & 6 & 0.000 & 43,099 & 7 & 0.001 \\
YFJS07 &  & 48,491 & 9 & 0.000 & 48,355 & 10 & 0.002 & 40,941 & 12 & 0.001 & 40,941 & 12 & 0.002 & 38,453 & 2 & 0.000 & 38,453 & 2 & 0.000 \\
YFJS08 &  & 47,429 & 3 & 0.000 & 41,754 & 9 & 0.001 & 41,994 & 3 & 0.000 & 37,681 & 10 & 0.001 & 35,429 & 2 & 0.000 & 33,392 & 5 & 0.000 \\
YFJS09 &  & 24,413 & 2 & 0.000 & 23,944 & 4 & 0.000 & 23,059 & 2 & 0.000 & 22,211 & 4 & 0.000 & 20,916 & 7 & 0.000 & 20,916 & 7 & 0.001 \\
YFJS10 &  & 42,583 & 1 & 0.000 & 42,278 & 2 & 0.000 & 39,126 & 4 & 0.000 & 39,361 & 4 & 0.000 & 34,487 & 7 & 0.000 & 34,487 & 10 & 0.001 \\
YFJS11 &  & 52,815 & 4 & 0.000 & 52,196 & 8 & 0.002 & 49,781 & 2 & 0.000 & 49,206 & 4 & 0.001 & 45,051 & 1 & 0.000 & 44,729 & 2 & 0.000 \\
YFJS12 &  & 53,094 & 11 & 0.001 & 52,650 & 12 & 0.004 & 56,022 & 7 & 0.000 & 54,186 & 12 & 0.004 & 46,471 & 5 & 0.000 & 43,150 & 9 & 0.003 \\
YFJS13 &  & 47,064 & 4 & 0.000 & 46,578 & 9 & 0.003 & 41,618 & 5 & 0.000 & 41,498 & 6 & 0.002 & 36,962 & 4 & 0.000 & 35,800 & 11 & 0.003 \\
YFJS14 &  & 119,792 & 13 & 0.016 & 118,731 & 27 & 0.202 & 97,201 & 18 & 0.022 & 100,177 & 23 & 0.161 & 89,371 & 6 & 0.008 & 87,518 & 21 & 0.162 \\
YFJS15 &  & 121,960 & 14 & 0.020 & 125,754 & 66 & 0.553 & 110,858 & 8 & 0.011 & 109,314 & 33 & 0.284 & 90,041 & 15 & 0.019 & 88,664 & 28 & 0.220 \\
YFJS16 &  & 116,033 & 19 & 0.024 & 115,312 & 31 & 0.204 & 101,073 & 7 & 0.008 & 98,976 & 34 & 0.238 & 88,521 & 5 & 0.006 & 85,252 & 27 & 0.176 \\
YFJS17 &  & 104,862 & 9 & 0.025 & 104,466 & 18 & 0.396 & 83,240 & 14 & 0.042 & 82,167 & 42 & 0.895 & 71,114 & 9 & 0.027 & 70,851 & 46 & 1.004 \\
YFJS18 &  & 124,738 & 13 & 0.043 & 119,224 & 84 & 1.872 & 97,366 & 9 & 0.027 & 94,472 & 45 & 1.123 & 83,681 & 9 & 0.030 & 82,795 & 27 & 0.639 \\
YFJS19 &  & 98,795 & 24 & 0.072 & 100,158 & 38 & 0.847 & 86,289 & 9 & 0.031 & 83,530 & 48 & 1.159 & 70,454 & 22 & 0.061 & 68,592 & 69 & 1.435 \\
YFJS20 &  & 95,573 & 10 & 0.030 & 94,201 & 42 & 0.764 & 79,653 & 24 & 0.062 & 79,638 & 35 & 0.630 & 70,928 & 9 & 0.028 & 70,253 & 19 & 0.365 \\ 
\hline
mean &  & 71,227.20 & 8.30 & 0.012 & 70,704.30 & 19.60 & 0.243 & 62,154.65 & 7.70 & 0.010 & 61,446.60 & 17.20 & 0.225 & 54,883.00 & 6.35 & 0.009 & 53,749.15 & 16.00 & 0.201 \\
wins &  & 6 & & & 17 & & & 8 & & & 17 & & & 7 & & & 19 & & \\
\hline
\end{tabular}}
\end{center}
\caption{Local search results with the best-improvement strategy
  using the reduced neighborhood and the cropped neighborhood, applied
  to the large-sized instances proposed in~\cite{birgin2014milp} with learning
  rate $\alpha \in \{0.1, 0.2, 0.3\}$.}
\label{tab4}
\end{table}

\subsection{Experiments with the metaheuristics in the large-sized instances}

In this section we report experiments to evaluate the behavior of the metaheuristics on the~50 large-sized instances, totaling 150 instances as we vary $\alpha \in \{0.1, 0.2, 0.3\}$. The ILS and GRASP metaheuristics (Algorithms~\ref{algo7} and~\ref{algo9}, respectively) make use of the local search (Algorithm~\ref{algo6}). We analyzed these two metaheuristics in connection with the local search using the reduced neighborhood (RN) and in connection with the local search using the cropped neighborhood (CN). We call these versions of ILS-RN, ILS-CN, GRASP-RN and GRASP-CN, respectively, hereafter. The tabu search metaheuristic was already conceived using the reduced neighborhood (see lines 16 and 17 of Algorithm~\ref{algo12}), and applying to it the concept of cropped neighborhood corresponds simply to changing in line~11 of Algorithm~\ref{algo12}, $v \in \mathcal{O}$ by $v \in \mathcal{P}$. Hereafter we refer to these two versions as TS-RN and TS-CN, respectively. The simulated annealing metaheuristic does not use local search and therefore we considered only a single version of it, which we denote SA hereafter. In all metaheuristics we used the total CPU time as a stopping criterion and considered a maximum limit of~5 minutes. Additionally, and only during the parameters calibration phase, we placed an additional stopping criterion, also based on CPU time, which is satisfied if~5 seconds elapse without improvement of the incumbent solution.

The ILS metaheuristic has as parameters $\ell^p_{\min}$ and $\ell^p_{\max} \in \mathbb{N}$ which must satisfy $1 \leq \ell^p_{\min} \leq \ell^p_{\max}$. In the calibration experiments, we evaluated $\ell^p_{\min}, \ell^p_{\max} \in \{1,2,3,4,5\}$. The GRASP metaheuristic has as a single parameter $\alpha \in \mathbb{R}$ which must satisfy $0 \leq \alpha \leq 1$. In the calibration experiments we evaluated $\alpha \in \{ 0.10, 0.11, \dots, 0.99 \}$. The TS metaheuristic has as a single parameter $t_{\max} \in \mathbb{N}$. In the calibration experiments, we evaluated $t_{\max} = \lceil (|\mathcal{O}| + |\mathcal{F}|) \times t \rceil$ with $t \in \{ 0.1, 0.2, \dots, 5.0 \}$. In the SA metaheuristic, the parameters are $\bar{\ell} \in \mathbb{N}$ and $T_0, T_f, \delta \in \mathbb{R}$ and must satisfy $1 \leq \bar \ell$, $T_0 \geq T_f > 0$ and $\delta \in (0,1)$. In the calibration experiments, we evaluated $\bar \ell \in \{ 1{,}000, 2{,}000, \dots, 5{,}000 \}$, $T_0 = - T_0^p / \ln(T_0^m)$ with $T_0^p \in \{0.10, 0.11, \dots, 0.99 \}$ and $T_0^m \in \{0.10, 0.11, \dots, 0.99 \}$, $T_f \in \{ 10^{-1}, 10^{-3}, 10^{-5}\}$, and $\delta \in \{ 0.80, 0.81, \dots, 0.99 \}$.
(See~\cite{Johnson1989} for details on how to choose the initial temperature in SA). We calibrated the seven methods independently using the irace package~\cite{LpezIbez2016}. In irace all default parameters were considered, except for $\mathrm{maxExperiments} = 10{,}000$. We considered in the calibration only the~50 large-sized instances with learning rate $\alpha = 0.2$. We separated these instances into two sets of~25 instances each, one for training and the other for testing. As a result, irace returned the following parameters:
\begin{description}
  \item[ILS-RN:] $\ell^p_{\min} = 2$ and $\ell^p_{\max} = 4$,
  \item[ILS-CN:] $\ell^p_{\min} = 1$ and $\ell^p_{\max} = 3$,
  \item[GRASP-RN:] $\alpha = 0.38$,
  \item[GRASP-CN:] $\alpha = 0.59$,
  \item[TS-RN:] $t_{\max} = \lceil (|\mathcal{O} + |\mathcal{F}|) \times t \rceil$ with $t=0.9$,
  \item[TS-CN:] $t_{\max} = \lceil (|\mathcal{O} + |\mathcal{F}|) \times t \rceil$ with $t=0.5$,
  \item[SA:] $\bar \ell = 3$, $T_0 = - T_0^p / \ln(T_0^m)$ with $T_0^p = 0.78$ and $T_0^m = 0.79$, $T_f = 10^{-3}$, and $\delta = 0.82$.
\end{description}

With the selected parameters, we ran the metaheuristics on the 50 large instances by varying the learning rate $\alpha \in \{0.1, 0.2, 0.3\}$. Since methods ILS-RN, ILS-CN, GRASP-RN, GRASP-CN, and SA have random components, these methods were applied 5 times to each instance. The same does not apply for methods TS-RN and TS-CN, which were run only once for each instance. Tables~\ref{tab5}, \ref{tab6}, and \ref{tab7} show the results associated with the learning rates $\alpha \in \{0.1, 0.2, 0.3\}$, respectively. In the tables, for each method/instance pair, when applicable, we show the best makespan (among the five runs) and the CPU time the method needed to achieve it. If a method reached its best makespan more than once, the indicated CPU time corresponds in fact to the average CPU time among all the times the method reached the indicated makespan. In the tables, additional lines show some statistics. The line named $\bar C_{\max}$ shows the average of the makespans found by each method (for each instance, the average considers the best over the five runs, when applicable). The line named \#best shows, for each method, the number of instances in which the method found the best solution among the solutions found by the seven methods being considered. The line named gap(\%), shows, for each method, the average gap when, for each instance, the solution found is compared with the best solution found by the seven methods. If the seven methods applied to the same instance found makespans $c_1, \dots, c_7$, the gap of method~$i$ in the considered instance corresponds to $100\% \times (c_i - c_{\min}) / c_{\min}$, where $c_{\min} = \min \{c_1,\dots,c_7\}$. The gap reported in line gap(\%) for each method corresponds to the average over all instances in the table. In the previous sentence, if a method $i$ is run five times per instance, the value $c_i$ made reference to the best over the five runs. So the value reported in line gap(\%) does not refer to what happened in ``the other four runs''. In these other cases, does each method find something close to its best makespan or something much worse? To answer this question, let us now consider that the values $\bar c_1, \dots, \bar c_7$ correspond to the average of the makespans found in the five runs. The line $\overline{\mathrm{gap}}$(\%) shows the same as line gap(\%) but considering $\bar c_1, \dots, \bar c_7$ instead of $c_1, \dots, c_7$. Both lines coincide for the case of TS-RN and TS-CN which have no random components and are therefore only executed once per instance. Thus, when TS-RN or TS-CN is compared with other method, it is the value of $\overline{\mathrm{gap}}$(\%) that must be compared. Otherwise, we would be comparing a method that is run only once against a method that is run five times a keeps the best solution found only.

%Up to this point, the four mentioned lines refer to statistics related to the data shown in the tables. But each method (except for TS-RN and TS-CN) is run five times in each instance; and out of the five runs, only the best makespan is reported in the tables. Nothing was said yet about what happens to the method in the runs where it does not find its best makespan. Does it find something close to its best makespan or something much worse? To answer this question, we calculated the average makespan of the five runs of the same method on the same instance and then calculated the gap from this average makespan to the best makespan found considering all the methods applied to that instance. Finally, for each method, we calculated the average of that gap over all the instances. This value is reported in the line named $\overline{\mathrm{gap}}$(\%). This value is, by definition, greater than or equal to the value reported in the line named gap(\%). If, for a given method, the two values are very similar, it means that, in the runs where it does not find its best solution, the method finds something similar. For the case of TS-RN and TS-CN, which have no random components, the two values are the same. When TS-RN or TS-CN is compared with other method, it is the value of $\overline{\mathrm{gap}}$(\%) that must be compared. Otherwise, we would be comparing a method that is run only once against a method that is run five times a keeps the best solution found only.

\begin{figure}
\centering
\input{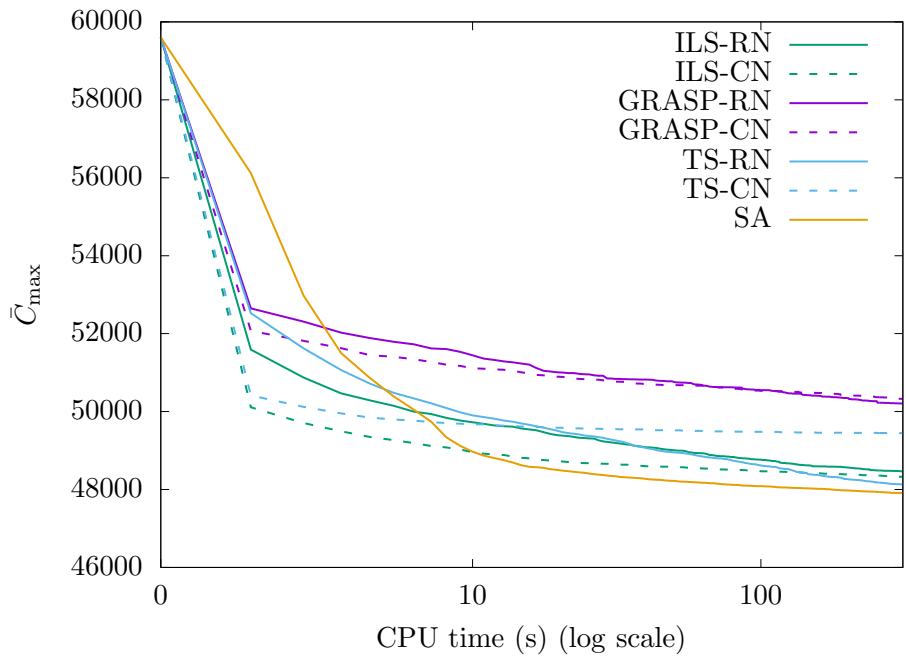}
\caption{This figure shows the value of $\bar C_{\max}$ as a function of the CPU time for each of the seven methods applied to the large-sized instances. The average makespan shown in the graphic is the average of the different values shown in Tables~\ref{tab5}, \ref{tab6}, and \ref{tab7}, that separate the instances by DA type, Y type and learning rates $\alpha \in \{0.1, 0.2, 0.3\}$.}
\label{fig3}
\end{figure}

\begin{figure}
\centering
\input{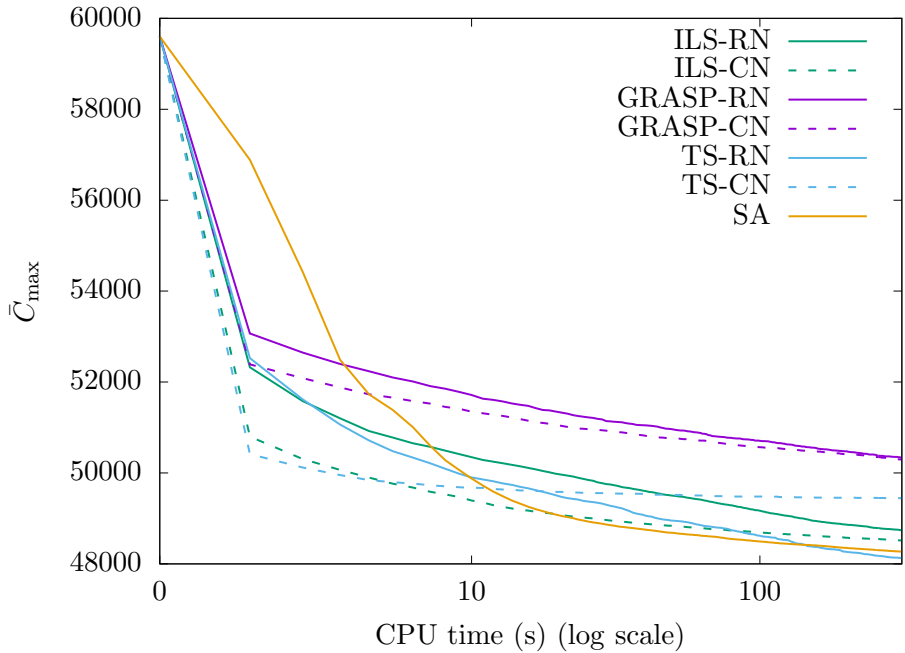}
\caption{This figure is similar to Figure~\ref{fig3}. The difference is that, for methods that solve each instance five times, instead of considering the best of the five makespan, we consider the average of the five makespan. This measure is more fair when in the comparison there are methods that do not posses a random ingredient and, therefore, are run only one time per instance.}
\label{fig4}
\end{figure}

Figure~\ref{fig3} shows, for each method, the average value of $\bar C_{\max}$ as a function of the CPU time. Statistics, consolidating the individual statistics for the DA-type and Y-type instances and the learning rates $\alpha \in \{0.1, 0.2, 0.3\}$, appear in Table~\ref{tab17}. As mentioned above, the considered average is the average of the makespans shown in Tables~\ref{tab5}, \ref{tab6}, and \ref{tab7}, i.e., it is the average of the best makespan over the five runs for each method/instance pair. Using that value in the comparison is harmful for the TS-RN and TS-CN methods. Note that if, we consider five minutes of CPU time, the TN-RN method is the second best, even using this unfair measure in the comparison. Figure~\ref{fig4} shows the same thing except that, for methods that solve each instance five times, instead of considering the best of the five makespan, we consider the average of the five makespan. With this measure, and considering that we have five minutes of CPU available, the TS-RN method appears as the best. The ranking of all the other methods, considering the five-minutos CPU time budget, remains the same.

Figures~\ref{fig3} and~\ref{fig4} show that the ranking of the methods in relation to the average quality of the solutions found is affected by the time limit used in the stopping criterion, since shorter times would yield a different ranking. On the other hand, the figures also show that longer times would have little effect on the comparison. (Note that the time axis in the figures is in logarithmic scale and that the methods almost stabilize at five minutes). Figures~\ref{fig3} and~\ref{fig4} and Table~\ref{tab17} clearly show that GRASP-RN, GRASP-CN and TS-CN were the worst performing methods. There seems to be an explanation for that: lack of diversification. In the case of TS-CN, the difference in efficacy is remarkable when compared to TS-RN. It is worth noting that the cropped neighborhood is radically smaller than the reduced neighborhood and this limits the ability of the tabu search to escape from local solutions. Something similar occurs in the two versions of GRASP, in which local searches are initiated from solutions constructed by a randomized version of the constructive heuristics ECT and EST. Apparently, the randomization of these constructive heuristics was not able to generate sufficiently diverse solutions and repeated runs of the local search ended up repeatedly converging to the same local solutions. On the other side of the spectrum, TS-RN and SA performed outstandingly well, followed closely by the two versions of ILS. The fact that SA (see Figures~\ref{fig3} and~\ref{fig4}) takes a little longer to get good quality solutions is entirely expected, since it is the method that less explores the neighborhood of the current solution in the search for a new solution. Regardless of anything else, it is important to note that all metaheuristics substantially improved the solutions delivered by the constructive heuristics and found by the local search. The constructive heuristics EST and ECT deliver solutions that are, on average, 32.36\% and 33.01\% worse than the best solutions found, respectively. The local search with the reduced neighborhood, which considered as a stand-alone method finds solutions of better quality than the local search with the cropped neighborhood, provides solutions that are on average 16.18\% worse than the best solutions found. On the other hand, all metaheuristics deliver solutions that are, on average, at no more than 6\% of the best solutions found. If we discard the three worst performing methods, this number drops to approximately 2\%. The highlight is for TS-RN and SA, which deliver solutions, on average, at approximately 1\% and 1.5\% of the best solutions found. The two versions of ILS deliver solutions that at approximately 2\% of the best solutions found, in average. SA also stands out for the number of best solutions obtained, 71 out of a total of 150, followed closely by ILS-RN and ILS-CN, with 59 and 58, respectively. TS-RN finds 55 best solutions with the merit of having run only once for each instance, while the others run five times and the best of the five is considered in these statistics.

In order to strengthen the comparison between the methods, we applied the Wilcoxon~\cite{Wilcoxon1945} test for each pair of methods, with a significance level of $\alpha=0.05$, to accept or reject the null hypothesis ``the samples of the two methods come from the same distribution'' or, equivalently, ``the difference between the samples of the two methods follows a symmetrical distribution around zero''. In this test we left out the GRASP-RN and GRASP-CN methods, which clearly performed worse than all the others. Table~\ref{tabwil} shows the results. Suppose that two methods $M_1$ and $M_2$ are compared and that these methods, when applied to a set of $N$ instances, find solutions with makespan $c_{1i}$ and $c_{2i}$, respectively, for $i=1,\dots,N$. Instances $i$ with $c_{2i}=c_{1i}$ are eliminated from the test. For simplicity, let's assume that $N$ already represents the number of instances in which the methods found different makespans. The test assigns each instance $i$ a weight $R_i$ which corresponds to the position of the instance when they are ordered from lowest to highest by the value of $|c_{2i}-c_{1i}|$. In the case of ties, $R_i$ corresponds to the average of the positions of the tied values. With the values of $R_i$ for each instance, we calculate the statistic
\[
W = R^+ - R^-, \text{ where } R^+ = \sum_{\{ i \;|\; c_{2i}>c_{1i}\}} R_i \text{ and } R^- = \sum_{\{ i \;|\; c_{1i}>c_{2i}\}} R_i.
\]
Under the null hypothesis and for large values of $N$, the value of $W$ is divided by its variance $N(N+1)(2N+1)/6$ to obtain $z$. In a two-tailed test, the null hypothesis is rejected if $|z|>z_{\mathrm{critical}}(\alpha/2)$ or, equivalently, if $p < \alpha$. The table shows, for each pair of methods compared, the values of $R^+$, $R^-$, and the $p$-value of $z$. The comparison considering a CPU time limit of 10 seconds shows that ILS-RN is worse than the other four methods compared and that the other four methods are equivalent to each other. This coincides with Figure~\ref{fig4}, which uses the average makespan as a measure of comparison. The comparison considering a CPU time limit of 5 minutes shows that all the methods are different. Furthermore, the two-by-two comparison makes it possible to construct a ranking between the methods in which TS-RN is the best of all, followed by SA, ILS-CN, ILS-RN and TS-CN in that order. This ranking also coincides with that shown on the right-hand side of Figure~\ref{fig4}, which uses the average makespan as a criterion.

\begin{table}
\centering
\begin{tabular}{|c|ccc|ccc|} 
\hline
\multirow{2}{*}{comparison} & \multicolumn{3}{c|}{10 seconds of CPU time limit} & \multicolumn{3}{c|}{5 minutes of CPU time limit} \\ 
\cline{2-7}
& $R^+$ & $R^-$ & $p$-value 
& $R^+$ & $R^-$ & $p$-value \\
\hline
\hline
ILS-RN versus ILS-CN &   804 & 7,974 & 0.000 & 1,402 & 5,385 & 0.000 \\
ILS-RN versus TS-RN  & 2,544 & 6,636 & 0.000 & 1,692 & 6,437 & 0.000 \\
ILS-RN versus TS-CN  & 3,844 & 6,887 & 0.003 & 6,311 & 4,274 & 0.044 \\
ILS-RN versus SA     & 3,579 & 7,006 & 0.001 & 3,221 & 7,075 & 0.000 \\
ILS-CN versus TS-RN  & 5,036 & 4,009 & 0.254 & 1,996 & 5,879 & 0.000 \\
ILS-CN versus TS-CN  & 5,296 & 5,289 & 0.994 & 6,777 & 3,808 & 0.003 \\
ILS-CN versus SA     & 5,624 & 4,961 & 0.513 & 4,102 & 6,051 & 0.047 \\
TS-RN  versus TS-CN  & 5,015 & 5,570 & 0.584 & 9,533 & 1,052 & 0.000 \\
TS-RN  versus SA     & 5,623 & 4,962 & 0.514 & 7,142 & 3,443 & 0.000 \\
TS-CN  versus SA     & 6,452 & 4,724 & 0.102 & 1,939 & 9,236 & 0.000 \\ 

\hline
\end{tabular}
\caption{Details of the Wilcoxon test comparing each pair of methods, to accept or reject the null hypothesis ``the difference between the two methods follows a symmetrical distribution around zero''.}
\label{tabwil}
\end{table}

We close this section by analyzing the behavior of the seven methods in the seven large-sized instances for which proved optimal solutions were reported in~\cite{arbro2023}. Table~\ref{tab8} shows the results. The table shows the instance name and instance learning rate for the seven instances for which the optimal value $C_{\max}^*$ is known. Then, for each method, it shows, as in the previous tables, when applicable, the best makespan found over the five runs in each of these instances. It can be seen that ILS-RN, ILS-CN and SA found the optimal solutions in all the seven instances, followed by TS-RN that found six out of the seven optimal solutions. Independently of finding the optimal solutions, appart from TS-CN, the methods found solutions with an average gap to the optimal solution smaller than 1\%. It must be highlighted that ILS-RN and ILS-CN found the optimal solutions in all the five runs of each instance. An evidence of this is the line ``$\overline{\mathrm{gap}}$ to optimal'' that exhibits a 0.0\% gap.

\begin{table}[ht!]
\begin{center}
\resizebox{\textwidth}{!}{
\begin{tabular}{|c|cc|cc|cc|cc|cc|cc|cc|}
\hline
\multirow{2}{*}{instance} & \multicolumn{2}{c|}{ILS-RN} & \multicolumn{2}{c|}{ILS-CN} & \multicolumn{2}{c|}{GRASP-RN} & \multicolumn{2}{c|}{GRASP-CN} & \multicolumn{2}{c|}{TS-RN} & \multicolumn{2}{c|}{TS-CN} & \multicolumn{2}{c|}{SA} \\
\cline{2-15} 
& $C_{\max}$ & Time & $C_{\max}$ & Time & $C_{\max}$ & Time & $C_{\max}$ & Time & $C_{\max}$ & Time & $C_{\max}$ & Time & $C_{\max}$ & Time \\
\hline
\hline
DAFJS01 & \textbf{23,460} & 12.29 & \textbf{23,460} & 3.39 & 24,172 & 90.48 & 23,961 & 131.52 & 24,424 & 0.05 & 24,661 & 0.01 & 23,992 & 46.48 \\
DAFJS02 & \textbf{26,432} & 60.95 & \textbf{26,432} & 25.84 & 26,989 & 0.01 & 26,869 & 0.28 & \textbf{26,432} & 0.32 & 27,620 & 0.00 & 27,577 & 0.55 \\
DAFJS03 & \textbf{50,343} & 164.37 & 50,462 & 34.55 & 50,424 & 241.44 & 51,028 & 156.99 & \textbf{50,343} & 7.20 & 51,612 & 0.26 & 50,458 & 261.30 \\
DAFJS04 & \textbf{51,872} & 62.24 & 51,992 & 181.28 & 51,927 & 30.24 & 52,358 & 275.80 & 52,291 & 0.11 & 53,348 & 0.00 & 51,927 & 151.36 \\
DAFJS05 & \textbf{34,660} & 167.91 & \textbf{34,660} & 45.10 & 37,081 & 35.04 & 36,883 & 80.14 & 35,278 & 1.96 & 36,009 & 0.26 & 36,141 & 139.08 \\
DAFJS06 & 35,435 & 111.76 & \textbf{35,302} & 176.69 & 37,023 & 122.84 & 36,338 & 98.55 & 35,636 & 15.04 & 37,854 & 0.05 & 35,697 & 256.36 \\
DAFJS07 & 45,805 & 240.68 & 45,306 & 72.33 & 46,624 & 107.20 & 46,791 & 57.95 & 44,281 & 158.67 & 45,060 & 297.79 & \textbf{44,219} & 241.82 \\
DAFJS08 & 53,202 & 52.07 & 53,463 & 41.54 & 53,575 & 174.57 & 54,129 & 129.70 & \textbf{53,002} & 261.88 & 54,018 & 0.69 & 53,292 & 276.00 \\
DAFJS09 & 40,573 & 145.06 & \textbf{40,537} & 166.88 & 41,066 & 119.76 & 40,759 & 184.35 & 41,123 & 4.34 & 42,241 & 0.09 & 41,265 & 4.00 \\
DAFJS10 & 44,523 & 99.32 & \textbf{44,058} & 254.61 & 44,917 & 217.73 & 44,735 & 277.52 & 45,820 & 10.28 & 47,438 & 0.32 & 44,851 & 70.25 \\
DAFJS11 & 55,180 & 121.71 & 55,393 & 151.06 & 58,257 & 219.22 & 58,527 & 228.89 & \textbf{54,351} & 260.24 & 55,337 & 158.73 & 54,708 & 264.50 \\
DAFJS12 & 54,035 & 227.92 & 53,175 & 20.54 & 56,670 & 237.70 & 56,495 & 79.61 & 51,432 & 238.43 & 51,794 & 224.88 & \textbf{51,067} & 234.61 \\
DAFJS13 & 54,654 & 35.35 & 54,371 & 0.05 & 54,392 & 277.63 & 54,583 & 107.84 & \textbf{53,340} & 233.35 & 53,792 & 0.10 & 54,502 & 93.52 \\
DAFJS14 & 61,500 & 133.79 & 60,882 & 169.24 & 61,487 & 185.43 & 61,641 & 265.85 & 59,668 & 196.33 & 61,524 & 0.74 & \textbf{59,573} & 254.00 \\
DAFJS15 & 57,552 & 111.67 & 57,554 & 216.92 & 59,362 & 201.70 & 59,997 & 220.71 & \textbf{55,086} & 253.37 & 55,381 & 115.31 & 55,164 & 197.01 \\
DAFJS16 & 58,828 & 213.11 & 58,327 & 98.92 & 63,365 & 189.53 & 62,494 & 242.24 & 55,927 & 242.72 & 56,758 & 115.48 & \textbf{55,868} & 208.74 \\
DAFJS17 & 65,252 & 138.93 & 65,121 & 1.17 & 64,804 & 269.21 & 65,643 & 74.78 & \textbf{64,262} & 283.76 & 65,525 & 1.33 & 65,602 & 55.29 \\
DAFJS18 & 66,429 & 288.91 & 65,889 & 271.51 & 66,633 & 179.15 & 66,419 & 214.60 & 64,916 & 134.40 & 67,040 & 0.23 & \textbf{64,820} & 219.02 \\
DAFJS19 & 45,416 & 136.34 & 44,744 & 95.28 & 46,757 & 147.62 & 46,301 & 76.22 & 44,633 & 121.09 & 45,602 & 1.30 & \textbf{44,086} & 144.96 \\
DAFJS20 & 59,190 & 298.29 & 58,137 & 241.93 & 59,670 & 43.48 & 59,323 & 79.25 & \textbf{56,432} & 155.97 & 57,141 & 9.73 & 57,245 & 290.79 \\
DAFJS21 & 65,432 & 83.75 & 66,006 & 217.22 & 66,547 & 87.71 & 66,752 & 277.45 & 64,563 & 264.94 & 64,902 & 13.37 & \textbf{64,059} & 124.80 \\
DAFJS22 & 58,025 & 2.64 & 57,426 & 227.59 & 59,477 & 62.55 & 59,086 & 30.34 & 55,246 & 269.80 & \textbf{55,038} & 24.67 & 55,624 & 283.70 \\
DAFJS23 & 42,818 & 226.75 & 42,255 & 149.85 & 43,565 & 66.68 & 43,426 & 150.96 & 41,941 & 57.87 & 42,426 & 1.56 & \textbf{41,286} & 227.89 \\
DAFJS24 & 48,266 & 207.07 & 48,911 & 123.21 & 51,313 & 127.05 & 51,532 & 43.92 & 47,443 & 139.28 & 47,365 & 281.95 & \textbf{46,945} & 217.47 \\
DAFJS25 & 62,813 & 127.68 & 63,039 & 16.33 & 65,321 & 147.01 & 65,866 & 65.19 & 60,824 & 205.95 & 60,308 & 248.14 & \textbf{60,187} & 270.84 \\
DAFJS26 & 62,389 & 131.50 & 61,645 & 158.57 & 65,982 & 299.75 & 65,951 & 196.11 & 58,995 & 185.68 & \textbf{58,749} & 67.28 & 59,450 & 297.29 \\
DAFJS27 & 69,386 & 73.63 & 68,795 & 254.71 & 70,375 & 231.60 & 70,762 & 6.13 & 66,539 & 245.31 & \textbf{66,003} & 148.50 & 66,070 & 204.00 \\
DAFJS28 & 48,091 & 133.90 & 47,738 & 187.43 & 48,669 & 164.03 & 49,309 & 294.79 & 47,428 & 164.05 & 47,551 & 75.18 & \textbf{46,569} & 268.66 \\
DAFJS29 & 57,824 & 58.01 & 57,384 & 259.73 & 59,825 & 173.11 & 59,564 & 234.87 & 55,820 & 10.14 & 56,122 & 37.45 & \textbf{54,686} & 296.94 \\
DAFJS30 & 48,572 & 292.30 & 48,329 & 137.76 & 50,395 & 59.47 & 50,593 & 18.06 & 47,140 & 237.35 & 47,307 & 76.44 & \textbf{46,344} & 232.86 \\
\hline
$\bar{C}_{\max}$ & \multicolumn{2}{c|}{51,598.57} & \multicolumn{2}{c|}{51,359.77} & \multicolumn{2}{c|}{52,888.80} & \multicolumn{2}{c|}{52,937.17} & \multicolumn{2}{c|}{50,487.20} & \multicolumn{2}{c|}{51,184.20} & \multicolumn{2}{c|}{50,442.47} \\
\#best & \multicolumn{2}{c|}{5} & \multicolumn{2}{c|}{6} & \multicolumn{2}{c|}{0} & \multicolumn{2}{c|}{0} & \multicolumn{2}{c|}{8} & \multicolumn{2}{c|}{3} & \multicolumn{2}{c|}{13} \\
gap(\%) & \multicolumn{2}{c|}{2.79} & \multicolumn{2}{c|}{2.33} & \multicolumn{2}{c|}{5.46} & \multicolumn{2}{c|}{5.48} & \multicolumn{2}{c|}{0.91} & \multicolumn{2}{c|}{2.45} & \multicolumn{2}{c|}{0.84} \\
$\overline{\mathrm{gap}}$(\%) & \multicolumn{2}{c|}{3.43} & \multicolumn{2}{c|}{2.94} & \multicolumn{2}{c|}{6.23} & \multicolumn{2}{c|}{6.14} & \multicolumn{2}{c|}{0.91} & \multicolumn{2}{c|}{2.45} & \multicolumn{2}{c|}{1.87} \\
\hline
\multicolumn{15}{c}{}\\
\hline
\multirow{2}{*}{instance} & \multicolumn{2}{c|}{ILS-RN} & \multicolumn{2}{c|}{ILS-CN} & \multicolumn{2}{c|}{GRASP-RN} & \multicolumn{2}{c|}{GRASP-CN} & \multicolumn{2}{c|}{TS-RN} & \multicolumn{2}{c|}{TS-CN} & \multicolumn{2}{c|}{SA} \\
\cline{2-15} 
& $C_{\max}$ & Time & $C_{\max}$ & Time & $C_{\max}$ & Time & $C_{\max}$ & Time & $C_{\max}$ & Time & $C_{\max}$ & Time & $C_{\max}$ & Time \\
\hline
\hline
YFJS01 & \textbf{68,714} & 0.51 & \textbf{68,714} & 0.80 & \textbf{68,714} & 0.12 & 69,189 & 0.02 & \textbf{68,714} & 0.34 & 75,595 & 0.03 & \textbf{68,714} & 0.81 \\
YFJS02 & \textbf{72,465} & 1.40 & \textbf{72,465} & 2.82 & 73,053 & 0.12 & 73,149 & 0.44 & \textbf{72,465} & 0.02 & 73,053 & 0.04 & \textbf{72,465} & 1.79 \\
YFJS03 & \textbf{32,538} & 0.03 & \textbf{32,538} & 0.02 & \textbf{32,538} & 254.14 & \textbf{32,538} & 0.08 & \textbf{32,538} & 0.10 & 33,174 & 0.00 & \textbf{32,538} & 0.47 \\
YFJS04 & \textbf{35,883} & 1.73 & \textbf{35,883} & 0.17 & 36,946 & 0.15 & \textbf{35,883} & 0.10 & \textbf{35,883} & 1.25 & 37,340 & 0.02 & \textbf{35,883} & 0.59 \\
YFJS05 & \textbf{40,186} & 9.99 & \textbf{40,186} & 1.76 & 41,034 & 160.74 & 41,034 & 0.52 & \textbf{40,186} & 4.42 & 41,034 & 0.03 & 41,034 & 0.74 \\
YFJS06 & \textbf{40,441} & 218.75 & \textbf{40,441} & 122.67 & 41,123 & 3.28 & 41,174 & 59.04 & 40,522 & 13.18 & 41,123 & 0.61 & 40,784 & 0.76 \\
YFJS07 & \textbf{40,887} & 131.61 & \textbf{40,887} & 37.21 & 41,625 & 54.84 & 41,354 & 11.63 & \textbf{40,887} & 38.15 & 41,394 & 0.17 & \textbf{40,887} & 16.55 \\
YFJS08 & \textbf{32,573} & 0.69 & \textbf{32,573} & 0.88 & 32,623 & 0.09 & 32,655 & 29.50 & \textbf{32,573} & 0.91 & 34,452 & 0.00 & \textbf{32,573} & 63.60 \\
YFJS09 & \textbf{22,681} & 2.68 & \textbf{22,681} & 0.54 & \textbf{22,681} & 6.30 & \textbf{22,681} & 56.94 & \textbf{22,681} & 0.01 & \textbf{22,681} & 0.24 & \textbf{22,681} & 2.03 \\
YFJS10 & \textbf{37,372} & 11.63 & \textbf{37,372} & 19.65 & 37,685 & 1.92 & 37,685 & 0.62 & 39,621 & 0.03 & 37,482 & 0.00 & \textbf{37,372} & 3.04 \\
YFJS11 & \textbf{47,800} & 12.50 & \textbf{47,800} & 94.06 & 49,493 & 174.55 & 49,208 & 41.96 & 48,195 & 4.95 & 49,721 & 0.04 & \textbf{47,800} & 85.90 \\
YFJS12 & \textbf{46,728} & 26.81 & \textbf{46,728} & 15.79 & 49,411 & 13.84 & 48,774 & 92.81 & \textbf{46,728} & 149.42 & 48,608 & 0.02 & \textbf{46,728} & 63.09 \\
YFJS13 & \textbf{36,911} & 12.30 & \textbf{36,911} & 13.38 & 38,042 & 2.35 & 37,361 & 108.86 & 36,926 & 33.80 & 39,207 & 0.02 & \textbf{36,911} & 76.79 \\
YFJS14 & \textbf{110,185} & 216.19 & 110,558 & 269.22 & 112,062 & 221.06 & 113,009 & 59.06 & 110,625 & 9.33 & 117,667 & 0.03 & 111,282 & 147.69 \\
YFJS15 & 105,935 & 192.42 & 106,221 & 189.75 & 108,296 & 17.04 & 111,374 & 226.73 & \textbf{105,837} & 259.22 & 116,525 & 0.19 & 106,892 & 149.18 \\
YFJS16 & \textbf{105,003} & 56.50 & 105,423 & 119.50 & 111,439 & 107.53 & 111,432 & 284.88 & 105,905 & 197.45 & 115,417 & 0.07 & 105,112 & 56.87 \\
YFJS17 & \textbf{93,469} & 198.07 & 94,558 & 174.10 & 99,912 & 159.01 & 101,069 & 107.68 & 96,737 & 226.68 & 98,329 & 0.99 & 94,950 & 112.29 \\
YFJS18 & \textbf{101,013} & 215.24 & 101,517 & 63.58 & 109,763 & 247.62 & 107,336 & 230.38 & 103,359 & 239.97 & 115,332 & 0.48 & 102,477 & 184.95 \\
YFJS19 & 82,879 & 221.94 & \textbf{81,172} & 115.47 & 94,020 & 109.95 & 92,172 & 67.18 & 86,500 & 242.82 & 88,400 & 2.14 & 81,946 & 93.62 \\
YFJS20 & \textbf{83,689} & 233.90 & 83,731 & 259.71 & 94,228 & 280.03 & 93,689 & 13.39 & 85,867 & 282.07 & 88,516 & 3.63 & 84,707 & 200.84 \\
\hline
$\bar{C}_{\max}$ & \multicolumn{2}{c|}{61,867.60} & \multicolumn{2}{c|}{61,917.95} & \multicolumn{2}{c|}{64,734.40} & \multicolumn{2}{c|}{64,638.30} & \multicolumn{2}{c|}{62,637.45} & \multicolumn{2}{c|}{65,752.50} & \multicolumn{2}{c|}{62,186.80} \\
\#best & \multicolumn{2}{c|}{18} & \multicolumn{2}{c|}{14} & \multicolumn{2}{c|}{3} & \multicolumn{2}{c|}{3} & \multicolumn{2}{c|}{10} & \multicolumn{2}{c|}{1} & \multicolumn{2}{c|}{11} \\
gap(\%) & \multicolumn{2}{c|}{0.11} & \multicolumn{2}{c|}{0.14} & \multicolumn{2}{c|}{3.84} & \multicolumn{2}{c|}{3.51} & \multicolumn{2}{c|}{1.17} & \multicolumn{2}{c|}{5.15} & \multicolumn{2}{c|}{0.51} \\
$\overline{\mathrm{gap}}$(\%) & \multicolumn{2}{c|}{0.40} & \multicolumn{2}{c|}{0.23} & \multicolumn{2}{c|}{4.59} & \multicolumn{2}{c|}{4.02} & \multicolumn{2}{c|}{1.17} & \multicolumn{2}{c|}{5.15} & \multicolumn{2}{c|}{0.97} \\
\hline
\end{tabular}}
\end{center}
\caption{Results of applying the metaheuristics to the large-sized instances with learning rate $\alpha = 0.1$.}
\label{tab5}
\end{table}

\begin{table}[ht!]
\begin{center}
\resizebox{\textwidth}{!}{
\begin{tabular}{|c|cc|cc|cc|cc|cc|cc|cc|}
\hline
\multirow{2}{*}{instance} & \multicolumn{2}{c|}{ILS-RN} & \multicolumn{2}{c|}{ILS-CN} & \multicolumn{2}{c|}{GRASP-RN} & \multicolumn{2}{c|}{GRASP-CN} & \multicolumn{2}{c|}{TS-RN} & \multicolumn{2}{c|}{TS-CN} & \multicolumn{2}{c|}{SA} \\
\cline{2-15} 
& $C_{\max}$ & Time & $C_{\max}$ & Time & $C_{\max}$ & Time & $C_{\max}$ & Time & $C_{\max}$ & Time & $C_{\max}$ & Time & $C_{\max}$ & Time \\
\hline
\hline
DAFJS01 & \textbf{21,683} & 17.45 & \textbf{21,683} & 7.46 & 21,861 & 8.41 & 21,861 & 136.04 & 22,286 & 0.03 & 22,703 & 0.01 & 21,936 & 0.38 \\
DAFJS02 & \textbf{24,235} & 48.48 & \textbf{24,235} & 55.62 & 24,618 & 0.05 & 24,571 & 4.00 & 24,875 & 0.14 & 24,477 & 0.00 & 25,141 & 0.51 \\
DAFJS03 & \textbf{44,266} & 213.30 & 44,453 & 208.60 & 44,481 & 211.71 & 45,079 & 169.49 & 44,419 & 12.80 & 45,275 & 0.37 & \textbf{44,266} & 147.45 \\
DAFJS04 & \textbf{44,713} & 54.08 & 44,938 & 53.17 & \textbf{44,713} & 298.64 & 45,644 & 29.88 & 44,821 & 1.41 & 46,807 & 0.28 & 44,816 & 239.37 \\
DAFJS05 & 30,904 & 28.28 & \textbf{30,773} & 144.29 & 31,848 & 54.88 & 32,065 & 22.88 & 32,028 & 0.92 & 31,772 & 0.66 & 31,295 & 61.32 \\
DAFJS06 & 30,602 & 65.94 & \textbf{30,553} & 16.14 & 31,994 & 90.56 & 31,768 & 145.27 & 31,521 & 5.17 & 31,039 & 0.52 & 31,216 & 44.70 \\
DAFJS07 & 40,343 & 46.23 & 39,788 & 223.54 & 40,321 & 121.98 & 40,891 & 84.05 & \textbf{38,731} & 37.49 & 39,488 & 223.99 & 38,833 & 112.43 \\
DAFJS08 & 45,636 & 273.46 & 45,784 & 218.76 & 45,955 & 9.16 & 46,589 & 7.49 & \textbf{45,162} & 133.36 & 46,129 & 119.13 & 45,302 & 128.92 \\
DAFJS09 & \textbf{35,098} & 25.22 & 35,244 & 74.17 & 35,759 & 18.62 & 35,507 & 224.49 & 35,161 & 5.16 & 36,004 & 0.20 & 35,418 & 193.84 \\
DAFJS10 & 37,326 & 238.28 & 37,172 & 198.27 & 37,736 & 257.38 & 37,856 & 76.13 & \textbf{36,791} & 293.02 & 38,453 & 0.75 & 37,183 & 11.88 \\
DAFJS11 & 47,417 & 203.56 & 47,221 & 248.85 & 49,271 & 182.94 & 49,049 & 30.18 & 45,861 & 300.00 & 46,470 & 169.13 & \textbf{45,822} & 242.59 \\
DAFJS12 & 45,613 & 256.81 & 45,394 & 45.36 & 48,195 & 14.33 & 47,896 & 195.90 & 44,167 & 232.32 & 43,932 & 259.97 & \textbf{43,227} & 261.22 \\
DAFJS13 & 45,230 & 112.83 & 45,187 & 255.42 & 45,154 & 195.56 & 45,422 & 200.48 & \textbf{44,180} & 128.30 & 46,329 & 0.17 & 45,309 & 154.25 \\
DAFJS14 & 50,958 & 250.28 & 50,897 & 103.42 & 51,082 & 235.59 & 50,954 & 267.23 & \textbf{49,611} & 171.05 & 51,800 & 0.36 & 50,555 & 271.75 \\
DAFJS15 & 49,326 & 277.73 & 49,120 & 158.41 & 50,696 & 137.86 & 50,910 & 95.17 & 47,211 & 221.78 & 47,414 & 152.93 & \textbf{47,151} & 282.47 \\
DAFJS16 & 50,300 & 40.01 & 49,896 & 227.89 & 53,242 & 121.42 & 54,087 & 243.46 & 48,972 & 297.39 & 48,312 & 169.85 & \textbf{47,661} & 169.18 \\
DAFJS17 & 53,388 & 122.05 & 52,783 & 157.04 & 52,905 & 171.59 & 53,513 & 237.27 & \textbf{51,715} & 55.57 & 52,459 & 2.64 & 52,081 & 94.05 \\
DAFJS18 & 54,313 & 275.00 & 54,237 & 103.14 & 54,569 & 292.99 & 55,154 & 147.54 & \textbf{52,624} & 164.38 & 54,008 & 1.39 & 53,942 & 183.68 \\
DAFJS19 & 39,480 & 94.81 & 39,312 & 77.07 & 41,046 & 124.48 & 40,931 & 162.08 & 38,740 & 130.01 & 38,941 & 6.77 & \textbf{38,012} & 150.83 \\
DAFJS20 & 49,099 & 259.17 & 49,226 & 160.00 & 49,290 & 184.19 & 50,041 & 197.61 & 48,088 & 230.79 & 49,274 & 0.55 & \textbf{47,680} & 282.84 \\
DAFJS21 & 54,311 & 5.09 & 54,163 & 118.65 & 54,851 & 216.71 & 55,002 & 20.65 & \textbf{52,681} & 264.76 & 53,996 & 9.95 & 53,196 & 291.79 \\
DAFJS22 & 46,679 & 213.56 & 46,468 & 161.65 & 48,169 & 285.94 & 47,938 & 42.82 & 46,088 & 288.61 & \textbf{44,851} & 59.63 & 45,011 & 251.83 \\
DAFJS23 & 37,459 & 248.79 & 36,733 & 262.51 & 38,106 & 230.28 & 38,487 & 243.86 & \textbf{36,137} & 166.72 & 36,190 & 13.53 & 36,373 & 49.89 \\
DAFJS24 & 42,129 & 124.24 & 41,596 & 275.85 & 43,617 & 113.82 & 43,693 & 295.60 & \textbf{40,079} & 216.52 & 41,307 & 26.26 & 40,453 & 297.11 \\
DAFJS25 & 53,096 & 203.63 & 52,727 & 274.34 & 54,913 & 107.96 & 55,011 & 220.93 & 52,604 & 264.38 & 50,965 & 248.46 & \textbf{50,887} & 116.74 \\
DAFJS26 & 52,492 & 211.10 & 52,147 & 218.93 & 54,563 & 21.44 & 54,439 & 241.66 & \textbf{50,006} & 140.87 & 50,465 & 74.69 & 50,062 & 110.48 \\
DAFJS27 & 57,249 & 244.03 & 56,428 & 299.00 & 58,674 & 259.45 & 59,078 & 30.84 & 55,599 & 293.44 & 55,398 & 29.41 & \textbf{55,291} & 270.23 \\
DAFJS28 & 41,928 & 91.85 & 42,013 & 119.00 & 42,847 & 90.69 & 43,197 & 208.33 & 40,686 & 26.31 & 41,780 & 30.78 & \textbf{40,431} & 212.68 \\
DAFJS29 & 49,863 & 282.73 & 49,538 & 162.49 & 51,391 & 213.00 & 51,881 & 188.89 & 48,953 & 155.28 & 48,942 & 2.58 & \textbf{47,091} & 293.59 \\
DAFJS30 & 42,001 & 145.84 & 42,138 & 130.38 & 43,276 & 119.87 & 43,617 & 4.51 & 40,839 & 246.80 & 40,746 & 63.14 & \textbf{40,325} & 137.56 \\
\hline
$\bar{C}_{\max}$ & \multicolumn{2}{c|}{43,904.57} & \multicolumn{2}{c|}{43,728.23} & \multicolumn{2}{c|}{44,838.10} & \multicolumn{2}{c|}{45,071.03} & \multicolumn{2}{c|}{43,021.20} & \multicolumn{2}{c|}{43,524.20} & \multicolumn{2}{c|}{42,865.53} \\
\#best & \multicolumn{2}{c|}{5} & \multicolumn{2}{c|}{4} & \multicolumn{2}{c|}{1} & \multicolumn{2}{c|}{0} & \multicolumn{2}{c|}{11} & \multicolumn{2}{c|}{1} & \multicolumn{2}{c|}{12} \\
gap(\%) & \multicolumn{2}{c|}{2.91} & \multicolumn{2}{c|}{2.52} & \multicolumn{2}{c|}{5.12} & \multicolumn{2}{c|}{5.62} & \multicolumn{2}{c|}{1.13} & \multicolumn{2}{c|}{2.27} & \multicolumn{2}{c|}{0.74} \\
$\overline{\mathrm{gap}}$(\%) & \multicolumn{2}{c|}{3.49} & \multicolumn{2}{c|}{3.15} & \multicolumn{2}{c|}{5.99} & \multicolumn{2}{c|}{6.30} & \multicolumn{2}{c|}{1.13} & \multicolumn{2}{c|}{2.27} & \multicolumn{2}{c|}{1.69} \\
\hline
\multicolumn{15}{c}{}\\
\hline
\multirow{2}{*}{instance} & \multicolumn{2}{c|}{ILS-RN} & \multicolumn{2}{c|}{ILS-CN} & \multicolumn{2}{c|}{GRASP-RN} & \multicolumn{2}{c|}{GRASP-CN} & \multicolumn{2}{c|}{TS-RN} & \multicolumn{2}{c|}{TS-CN} & \multicolumn{2}{c|}{SA} \\
\cline{2-15} 
& $C_{\max}$ & Time & $C_{\max}$ & Time & $C_{\max}$ & Time & $C_{\max}$ & Time & $C_{\max}$ & Time & $C_{\max}$ & Time & $C_{\max}$ & Time \\
\hline
\hline
YFJS01 & \textbf{61,606} & 1.10 & \textbf{61,606} & 1.12 & 64,456 & 0.19 & 62,327 & 9.05 & \textbf{61,606} & 3.43 & 62,366 & 0.05 & 62,124 & 2.24 \\
YFJS02 & \textbf{64,172} & 0.59 & \textbf{64,172} & 1.35 & 65,967 & 0.99 & 66,394 & 0.87 & \textbf{64,172} & 1.65 & 64,805 & 0.16 & \textbf{64,172} & 28.72 \\
YFJS03 & \textbf{30,073} & 0.36 & \textbf{30,073} & 0.26 & \textbf{30,073} & 0.08 & \textbf{30,073} & 0.02 & \textbf{30,073} & 0.01 & \textbf{30,073} & 0.00 & \textbf{30,073} & 0.59 \\
YFJS04 & \textbf{32,670} & 0.51 & \textbf{32,670} & 0.37 & 33,302 & 226.25 & 32,802 & 166.10 & \textbf{32,670} & 1.53 & 33,302 & 0.01 & \textbf{32,670} & 0.65 \\
YFJS05 & \textbf{37,182} & 7.66 & \textbf{37,182} & 6.16 & 37,670 & 0.33 & 37,255 & 0.79 & \textbf{37,182} & 4.32 & 38,682 & 0.01 & 37,682 & 0.68 \\
YFJS06 & \textbf{36,749} & 68.53 & \textbf{36,749} & 11.18 & 37,293 & 94.00 & 36,812 & 17.99 & \textbf{36,749} & 137.03 & 37,943 & 0.17 & 37,107 & 1.01 \\
YFJS07 & \textbf{36,611} & 102.75 & \textbf{36,611} & 57.41 & 37,159 & 3.75 & 37,159 & 8.11 & \textbf{36,611} & 77.85 & 37,159 & 0.34 & \textbf{36,611} & 1.07 \\
YFJS08 & \textbf{30,184} & 0.23 & \textbf{30,184} & 0.51 & 30,267 & 0.46 & 30,384 & 227.80 & 30,267 & 1.93 & 32,595 & 0.00 & \textbf{30,184} & 1.10 \\
YFJS09 & \textbf{21,363} & 12.05 & \textbf{21,363} & 1.63 & \textbf{21,363} & 130.79 & \textbf{21,363} & 141.45 & \textbf{21,363} & 10.87 & \textbf{21,363} & 0.89 & 21,450 & 2.64 \\
YFJS10 & \textbf{34,364} & 10.57 & \textbf{34,364} & 8.30 & 35,389 & 0.03 & 34,689 & 11.28 & 35,251 & 0.14 & 35,116 & 0.00 & \textbf{34,364} & 38.75 \\
YFJS11 & \textbf{43,454} & 158.07 & 43,544 & 164.27 & 45,613 & 12.83 & 45,649 & 33.90 & 44,094 & 0.53 & 46,983 & 0.00 & \textbf{43,454} & 81.93 \\
YFJS12 & \textbf{42,232} & 177.44 & 42,397 & 131.72 & 45,207 & 91.71 & 44,084 & 69.04 & 42,597 & 221.30 & 46,047 & 0.01 & \textbf{42,232} & 4.28 \\
YFJS13 & \textbf{33,875} & 34.01 & \textbf{33,875} & 18.97 & 34,401 & 4.99 & 34,240 & 50.77 & \textbf{33,875} & 34.33 & 34,335 & 0.20 & \textbf{33,875} & 3.48 \\
YFJS14 & \textbf{92,522} & 58.55 & 93,435 & 164.41 & 96,665 & 218.81 & 96,415 & 81.18 & 93,409 & 43.69 & 97,155 & 0.11 & 93,812 & 117.48 \\
YFJS15 & 91,555 & 188.91 & 91,006 & 203.20 & 96,697 & 210.91 & 97,126 & 71.71 & \textbf{90,801} & 216.11 & 95,840 & 0.40 & 91,216 & 253.51 \\
YFJS16 & 90,957 & 133.94 & 90,591 & 138.68 & 95,294 & 289.35 & 96,871 & 283.70 & 91,357 & 103.88 & 100,478 & 0.11 & \textbf{90,266} & 150.82 \\
YFJS17 & 78,747 & 203.06 & \textbf{78,585} & 294.38 & 85,736 & 0.00 & 85,736 & 0.00 & 78,982 & 276.75 & 83,631 & 0.18 & 79,177 & 283.88 \\
YFJS18 & 85,415 & 255.16 & \textbf{85,191} & 56.76 & 92,382 & 86.94 & 94,410 & 189.36 & 92,710 & 70.33 & 97,864 & 0.04 & 86,119 & 263.05 \\
YFJS19 & 70,527 & 299.05 & \textbf{68,718} & 122.25 & 78,364 & 245.26 & 79,427 & 14.98 & 74,284 & 184.84 & 71,701 & 8.45 & 69,662 & 57.52 \\
YFJS20 & 72,991 & 221.56 & \textbf{70,286} & 170.65 & 80,074 & 29.02 & 81,808 & 36.17 & 71,131 & 268.24 & 74,652 & 0.97 & 73,215 & 138.10 \\
\hline
$\bar{C}_{\max}$ & \multicolumn{2}{c|}{54,362.45} & \multicolumn{2}{c|}{54,130.10} & \multicolumn{2}{c|}{57,168.60} & \multicolumn{2}{c|}{57,251.20} & \multicolumn{2}{c|}{54,959.20} & \multicolumn{2}{c|}{57,104.50} & \multicolumn{2}{c|}{54,473.25} \\
\#best & \multicolumn{2}{c|}{14} & \multicolumn{2}{c|}{15} & \multicolumn{2}{c|}{2} & \multicolumn{2}{c|}{2} & \multicolumn{2}{c|}{10} & \multicolumn{2}{c|}{2} & \multicolumn{2}{c|}{10} \\
gap(\%) & \multicolumn{2}{c|}{0.43} & \multicolumn{2}{c|}{0.11} & \multicolumn{2}{c|}{4.63} & \multicolumn{2}{c|}{4.47} & \multicolumn{2}{c|}{1.30} & \multicolumn{2}{c|}{4.77} & \multicolumn{2}{c|}{0.64} \\
$\overline{\mathrm{gap}}$(\%) & \multicolumn{2}{c|}{0.80} & \multicolumn{2}{c|}{0.26} & \multicolumn{2}{c|}{5.67} & \multicolumn{2}{c|}{4.97} & \multicolumn{2}{c|}{1.30} & \multicolumn{2}{c|}{4.77} & \multicolumn{2}{c|}{1.36} \\
\hline
\end{tabular}}
\end{center}
\caption{Results of applying the metaheuristics to the large-sized instances with learning rate $\alpha = 0.2$.}
\label{tab6}
\end{table}

\begin{table}[ht!]
\begin{center}
\resizebox{\textwidth}{!}{
\begin{tabular}{|c|cc|cc|cc|cc|cc|cc|cc|}
\hline
\multirow{2}{*}{instance} & \multicolumn{2}{c|}{ILS-RN} & \multicolumn{2}{c|}{ILS-CN} & \multicolumn{2}{c|}{GRASP-RN} & \multicolumn{2}{c|}{GRASP-CN} & \multicolumn{2}{c|}{TS-RN} & \multicolumn{2}{c|}{TS-CN} & \multicolumn{2}{c|}{SA} \\
\cline{2-15} 
& $C_{\max}$ & Time & $C_{\max}$ & Time & $C_{\max}$ & Time & $C_{\max}$ & Time & $C_{\max}$ & Time & $C_{\max}$ & Time & $C_{\max}$ & Time \\
\hline
\hline
DAFJS01 & \textbf{19,716} & 42.15 & \textbf{19,716} & 93.26 & 19,969 & 181.01 & 19,994 & 1.14 & 20,137 & 0.35 & 21,082 & 0.01 & 20,092 & 0.59 \\
DAFJS02 & \textbf{22,147} & 44.81 & \textbf{22,147} & 33.01 & 22,404 & 1.20 & 22,236 & 0.04 & 22,191 & 0.44 & 22,429 & 0.00 & 22,695 & 25.70 \\
DAFJS03 & 38,916 & 229.98 & 38,994 & 282.80 & 39,280 & 94.82 & 40,454 & 111.31 & 38,928 & 105.79 & 40,724 & 7.17 & \textbf{38,853} & 42.33 \\
DAFJS04 & \textbf{39,117} & 275.95 & 39,492 & 39.94 & 39,358 & 60.46 & 39,943 & 283.63 & 39,347 & 1.53 & 40,425 & 0.53 & 39,285 & 246.66 \\
DAFJS05 & \textbf{27,192} & 175.00 & 27,221 & 266.61 & 27,438 & 55.27 & 28,152 & 92.66 & 27,401 & 0.83 & 28,767 & 0.14 & 27,453 & 0.86 \\
DAFJS06 & 26,565 & 224.96 & 26,552 & 138.78 & 27,706 & 107.01 & 27,111 & 240.56 & \textbf{26,535} & 10.26 & 28,134 & 0.06 & 26,763 & 23.70 \\
DAFJS07 & 35,287 & 15.78 & 35,262 & 252.14 & 36,318 & 168.90 & 36,309 & 197.95 & 35,735 & 8.84 & 35,074 & 74.05 & \textbf{34,493} & 271.96 \\
DAFJS08 & 39,120 & 225.81 & 39,413 & 106.90 & 39,394 & 121.44 & 40,345 & 135.78 & 38,629 & 192.08 & 39,620 & 201.57 & \textbf{38,605} & 53.16 \\
DAFJS09 & 30,627 & 93.31 & 30,945 & 4.12 & 31,421 & 14.02 & 31,004 & 255.55 & \textbf{30,353} & 79.14 & 32,242 & 0.11 & 30,833 & 46.92 \\
DAFJS10 & \textbf{30,880} & 131.80 & 31,191 & 202.19 & 31,544 & 192.89 & 31,802 & 279.43 & 30,983 & 182.18 & 31,853 & 1.87 & 31,252 & 197.99 \\
DAFJS11 & 40,212 & 68.69 & 40,549 & 187.16 & 41,870 & 242.15 & 42,027 & 84.21 & 39,227 & 253.69 & 39,751 & 124.91 & \textbf{38,537} & 252.92 \\
DAFJS12 & 38,128 & 83.41 & 38,671 & 61.78 & 40,625 & 271.43 & 40,937 & 93.63 & 37,351 & 174.33 & 37,394 & 263.18 & \textbf{37,030} & 286.43 \\
DAFJS13 & 37,672 & 48.20 & 37,699 & 163.35 & 37,806 & 244.56 & 38,010 & 273.30 & 37,324 & 59.41 & 38,811 & 0.32 & \textbf{37,197} & 242.77 \\
DAFJS14 & 42,469 & 275.82 & 42,260 & 88.40 & 42,350 & 198.57 & 42,507 & 278.76 & \textbf{41,203} & 164.66 & 41,969 & 0.43 & 41,854 & 234.65 \\
DAFJS15 & 42,125 & 156.71 & 42,172 & 190.76 & 43,215 & 186.52 & 43,776 & 40.58 & \textbf{40,092} & 252.92 & 40,345 & 40.70 & 40,465 & 180.67 \\
DAFJS16 & 43,378 & 39.90 & 42,572 & 191.37 & 45,354 & 201.36 & 45,439 & 214.16 & 41,337 & 298.16 & 41,201 & 18.64 & \textbf{40,487} & 102.62 \\
DAFJS17 & 43,350 & 289.76 & 43,143 & 193.23 & 43,698 & 248.64 & 43,121 & 173.83 & \textbf{42,367} & 156.85 & 42,622 & 6.69 & 43,083 & 219.02 \\
DAFJS18 & 44,635 & 0.79 & 45,043 & 140.76 & 45,611 & 263.32 & 45,211 & 76.37 & \textbf{43,899} & 105.90 & 43,981 & 0.84 & 44,417 & 250.00 \\
DAFJS19 & 34,218 & 89.09 & 34,104 & 52.36 & 35,401 & 198.41 & 35,425 & 94.86 & 33,501 & 43.82 & 34,187 & 0.66 & \textbf{33,250} & 238.49 \\
DAFJS20 & 41,035 & 31.13 & 41,242 & 196.95 & 41,550 & 118.66 & 42,268 & 165.12 & 40,059 & 214.87 & 40,245 & 8.86 & \textbf{39,812} & 262.00 \\
DAFJS21 & 45,000 & 50.80 & 44,519 & 215.72 & 45,409 & 141.46 & 44,883 & 249.15 & 44,190 & 297.58 & 43,840 & 25.11 & \textbf{43,684} & 173.57 \\
DAFJS22 & 37,634 & 290.20 & 37,610 & 55.93 & 38,705 & 135.46 & 39,392 & 124.86 & 37,083 & 120.04 & \textbf{36,276} & 131.71 & 36,866 & 124.16 \\
DAFJS23 & 33,077 & 95.58 & 32,960 & 10.78 & 33,756 & 107.08 & 33,543 & 44.93 & \textbf{31,860} & 39.67 & 32,412 & 15.90 & 32,017 & 188.80 \\
DAFJS24 & 36,252 & 239.59 & 36,121 & 68.70 & 37,838 & 92.47 & 37,964 & 29.84 & 35,122 & 299.58 & 35,273 & 54.63 & \textbf{34,422} & 268.69 \\
DAFJS25 & 44,698 & 215.09 & 45,061 & 277.98 & 46,361 & 44.37 & 47,007 & 1.62 & 43,654 & 198.69 & 43,176 & 115.77 & \textbf{43,050} & 244.18 \\
DAFJS26 & 44,105 & 263.42 & 43,877 & 72.88 & 47,208 & 44.13 & 46,799 & 79.24 & 43,025 & 189.92 & 42,082 & 249.03 & \textbf{42,080} & 211.51 \\
DAFJS27 & 47,643 & 81.38 & 47,648 & 269.30 & 48,748 & 188.12 & 49,400 & 275.14 & \textbf{46,192} & 284.20 & 46,887 & 17.37 & 46,510 & 204.55 \\
DAFJS28 & 36,496 & 246.88 & 36,494 & 299.25 & 37,703 & 157.63 & 37,808 & 192.62 & 35,694 & 31.02 & 36,285 & 31.11 & \textbf{35,637} & 191.18 \\
DAFJS29 & 42,923 & 148.58 & 43,334 & 244.96 & 45,174 & 110.20 & 44,735 & 40.40 & 42,218 & 267.02 & 42,648 & 22.61 & \textbf{41,410} & 281.11 \\
DAFJS30 & 37,051 & 148.59 & 37,006 & 190.53 & 38,067 & 194.66 & 37,421 & 258.42 & 35,609 & 250.85 & 35,593 & 14.62 & \textbf{35,242} & 283.70 \\
\hline
$\bar{C}_{\max}$ & \multicolumn{2}{c|}{37,388.93} & \multicolumn{2}{c|}{37,433.93} & \multicolumn{2}{c|}{38,376.03} & \multicolumn{2}{c|}{38,500.77} & \multicolumn{2}{c|}{36,708.20} & \multicolumn{2}{c|}{37,177.60} & \multicolumn{2}{c|}{36,579.13} \\
\#best & \multicolumn{2}{c|}{5} & \multicolumn{2}{c|}{2} & \multicolumn{2}{c|}{0} & \multicolumn{2}{c|}{0} & \multicolumn{2}{c|}{8} & \multicolumn{2}{c|}{1} & \multicolumn{2}{c|}{16} \\
gap(\%) & \multicolumn{2}{c|}{2.58} & \multicolumn{2}{c|}{2.71} & \multicolumn{2}{c|}{5.26} & \multicolumn{2}{c|}{5.57} & \multicolumn{2}{c|}{0.89} & \multicolumn{2}{c|}{2.38} & \multicolumn{2}{c|}{0.59} \\
$\overline{\mathrm{gap}}$(\%) & \multicolumn{2}{c|}{3.31} & \multicolumn{2}{c|}{3.17} & \multicolumn{2}{c|}{6.04} & \multicolumn{2}{c|}{6.23} & \multicolumn{2}{c|}{0.89} & \multicolumn{2}{c|}{2.38} & \multicolumn{2}{c|}{1.53} \\
\hline
\multicolumn{15}{c}{}\\
\hline
\multirow{2}{*}{instance} & \multicolumn{2}{c|}{ILS-RN} & \multicolumn{2}{c|}{ILS-CN} & \multicolumn{2}{c|}{GRASP-RN} & \multicolumn{2}{c|}{GRASP-CN} & \multicolumn{2}{c|}{TS-RN} & \multicolumn{2}{c|}{TS-CN} & \multicolumn{2}{c|}{SA} \\
\cline{2-15} 
& $C_{\max}$ & Time & $C_{\max}$ & Time & $C_{\max}$ & Time & $C_{\max}$ & Time & $C_{\max}$ & Time & $C_{\max}$ & Time & $C_{\max}$ & Time \\
\hline
\hline
YFJS01 & \textbf{55,128} & 1.60 & \textbf{55,128} & 1.29 & 55,458 & 0.18 & 56,790 & 2.84 & 55,693 & 0.92 & 57,900 & 0.02 & \textbf{55,128} & 0.73 \\
YFJS02 & \textbf{57,187} & 1.20 & \textbf{57,187} & 2.25 & 59,000 & 0.32 & 58,470 & 0.21 & \textbf{57,187} & 1.60 & 57,634 & 0.18 & \textbf{57,187} & 10.48 \\
YFJS03 & \textbf{27,686} & 0.03 & \textbf{27,686} & 0.03 & \textbf{27,686} & 0.00 & \textbf{27,686} & 0.01 & \textbf{27,686} & 0.02 & 28,772 & 0.00 & \textbf{27,686} & 0.47 \\
YFJS04 & \textbf{29,692} & 0.29 & \textbf{29,692} & 1.66 & \textbf{29,692} & 23.14 & 30,683 & 160.12 & \textbf{29,692} & 0.20 & \textbf{29,692} & 0.05 & \textbf{29,692} & 0.49 \\
YFJS05 & \textbf{33,736} & 2.48 & \textbf{33,736} & 3.45 & 34,430 & 2.41 & 33,856 & 182.52 & \textbf{33,736} & 12.59 & 35,577 & 0.01 & \textbf{33,736} & 0.90 \\
YFJS06 & \textbf{33,276} & 35.79 & \textbf{33,276} & 111.43 & 33,557 & 23.84 & 33,557 & 6.56 & \textbf{33,276} & 32.04 & 33,892 & 0.09 & 33,892 & 0.78 \\
YFJS07 & \textbf{33,011} & 91.58 & \textbf{33,011} & 82.50 & 33,352 & 30.98 & 33,059 & 126.97 & 33,092 & 4.76 & 33,317 & 0.01 & 33,092 & 0.55 \\
YFJS08 & \textbf{28,192} & 0.35 & \textbf{28,192} & 0.62 & 28,597 & 0.07 & 28,355 & 0.14 & \textbf{28,192} & 0.09 & 28,638 & 0.00 & \textbf{28,192} & 0.82 \\
YFJS09 & \textbf{20,077} & 3.69 & \textbf{20,077} & 2.38 & \textbf{20,077} & 33.66 & \textbf{20,077} & 22.83 & \textbf{20,077} & 9.58 & \textbf{20,077} & 0.23 & 20,121 & 165.05 \\
YFJS10 & \textbf{31,559} & 33.20 & \textbf{31,559} & 119.22 & 32,582 & 0.72 & 32,159 & 0.57 & 31,998 & 0.77 & 33,111 & 0.00 & \textbf{31,559} & 74.40 \\
YFJS11 & 40,204 & 102.14 & \textbf{40,103} & 105.15 & 41,079 & 136.76 & 40,619 & 20.45 & 40,350 & 6.53 & 41,547 & 0.00 & 40,319 & 2.58 \\
YFJS12 & \textbf{38,003} & 173.04 & \textbf{38,003} & 154.12 & 39,993 & 21.82 & 38,618 & 35.05 & 38,392 & 268.63 & 41,620 & 0.02 & 38,392 & 103.93 \\
YFJS13 & \textbf{30,711} & 88.73 & \textbf{30,711} & 89.29 & 31,823 & 1.59 & 31,518 & 55.89 & 30,714 & 180.10 & 33,449 & 0.01 & 30,874 & 1.20 \\
YFJS14 & 79,051 & 186.80 & 79,386 & 200.32 & 81,558 & 137.18 & 84,481 & 176.28 & \textbf{78,762} & 142.02 & 83,816 & 0.18 & 79,486 & 38.18 \\
YFJS15 & 79,374 & 139.96 & 78,660 & 198.52 & 84,181 & 133.33 & 84,277 & 230.45 & 79,092 & 191.45 & 87,555 & 0.10 & \textbf{78,545} & 69.16 \\
YFJS16 & 78,426 & 106.57 & 78,333 & 198.20 & 83,753 & 2.48 & 85,453 & 284.11 & 78,491 & 112.32 & 83,317 & 0.43 & \textbf{78,004} & 83.92 \\
YFJS17 & 66,188 & 295.79 & \textbf{66,029} & 82.47 & 73,682 & 0.00 & 73,682 & 0.00 & 66,802 & 275.82 & 71,161 & 0.04 & 66,388 & 158.14 \\
YFJS18 & 73,732 & 298.52 & \textbf{71,791} & 115.35 & 80,279 & 246.11 & 82,201 & 17.58 & 76,031 & 298.20 & 77,670 & 1.43 & 73,078 & 243.24 \\
YFJS19 & 61,283 & 149.66 & \textbf{58,933} & 123.03 & 70,115 & 118.17 & 68,545 & 185.31 & 60,919 & 297.73 & 61,788 & 4.36 & 59,633 & 299.10 \\
YFJS20 & 60,960 & 148.73 & \textbf{59,846} & 186.72 & 68,963 & 197.07 & 71,916 & 244.37 & 60,137 & 237.53 & 62,668 & 1.26 & 61,277 & 170.47 \\
\hline
$\bar{C}_{\max}$ & \multicolumn{2}{c|}{47,873.80} & \multicolumn{2}{c|}{47,566.95} & \multicolumn{2}{c|}{50,492.85} & \multicolumn{2}{c|}{50,800.10} & \multicolumn{2}{c|}{48,015.95} & \multicolumn{2}{c|}{50,160.05} & \multicolumn{2}{c|}{47,814.05} \\
\#best & \multicolumn{2}{c|}{12} & \multicolumn{2}{c|}{17} & \multicolumn{2}{c|}{3} & \multicolumn{2}{c|}{2} & \multicolumn{2}{c|}{8} & \multicolumn{2}{c|}{2} & \multicolumn{2}{c|}{9} \\
gap(\%) & \multicolumn{2}{c|}{0.55} & \multicolumn{2}{c|}{0.07} & \multicolumn{2}{c|}{4.97} & \multicolumn{2}{c|}{5.23} & \multicolumn{2}{c|}{0.83} & \multicolumn{2}{c|}{4.84} & \multicolumn{2}{c|}{0.56} \\
$\overline{\mathrm{gap}}$(\%) & \multicolumn{2}{c|}{0.94} & \multicolumn{2}{c|}{0.25} & \multicolumn{2}{c|}{5.87} & \multicolumn{2}{c|}{5.98} & \multicolumn{2}{c|}{0.83} & \multicolumn{2}{c|}{4.84} & \multicolumn{2}{c|}{1.39} \\
\hline
\end{tabular}}
\end{center}
\caption{Results of applying the metaheuristics to the large-sized instances with learning rate $\alpha = 0.3$.}
\label{tab7}
\end{table}

\begin{table}[ht!]
\begin{center}
\resizebox{\textwidth}{!}{
\begin{tabular}{|ccc|ccccccc|}
\hline
instance & $\alpha$ & $C_{\max}^\star$ & ILS-RN & ILS-CN & GRASP-RN & GRASP-CN & TS-RN & TS-CN & SA \\
\hline
\hline
YFJS03 & 0.1 & 32,538 & 32,538 & 32,538 & 32,538 & 32,538 & 32,538 & 33,174 & 32,538 \\
YFJS04 & 0.1 & 35,883 & 35,883 & 35,883 & 36,946 & 35,883 & 35,883 & 37,340 & 35,883 \\
YFJS08 & 0.1 & 32,573 & 32,573 & 32,573 & 32,623 & 32,655 & 32,573 & 34,452 & 32,573 \\
YFJS03 & 0.2 & 30,073 & 30,073 & 30,073 & 30,073 & 30,073 & 30,073 & 30,073 & 30,073 \\
YFJS08 & 0.2 & 30,184 & 30,184 & 30,184 & 30,267 & 30,384 & 30,267 & 32,595 & 30,184 \\
YFJS03 & 0.3 & 27,686 & 27,686 & 27,686 & 27,686 & 27,686 & 27,686 & 28,772 & 27,686 \\
YFJS08 & 0.3 & 28,192 & 28,192 & 28,192 & 28,597 & 28,355 & 28,192 & 28,638 & 28,192 \\
\hline
%\multicolumn{3}{|c|}{$\bar{C}_{\max}$} & 31,018.43 & 31,018.43 & 31,247.14 & 31,082.00 & 31,030.29 & 32,149.14 & 31,018.43 \\
\multicolumn{3}{|c|}{\#optimal} & 7 & 7 & 3 & 4 & 6 & 1 & 7 \\
\multicolumn{3}{|c|}{gap to optimal (\%)} & 0.00 & 0.00 & 0.69 & 0.21 & 0.04 & 3.61 & 0.00 \\
\multicolumn{3}{|c|}{$\overline{\mathrm{gap}}$ to optimal (\%)} & 0.00 & 0.00 & 0.73 & 0.31 & 0.04 & 3.61 & 0.24 \\
\hline
\end{tabular}}
\end{center}
\caption{Performance of the metaheuristics in the seven large-sized instances for which proved optimal solutions are known.}
\label{tab8}
\end{table}

\subsection{Experiments with the metaheuristics in the small-sized instances}

In this section we analyze the behavior of the metaheuristics on the small-sized instances. It is important to note that for all metaheuristics we considered their parameters exactly as defined in the previous section, i.e., with the calibration done for the large-sized instances. Tables~\ref{tab9}, \ref{tab10}, and~\ref{tab11} show the results for the instances with learning effect rate $\alpha \in \{ 0.1, 0.2, 0.3 \}$, respectively. Unlike in the large-sized instances, the ILS-RN, ILS-CN, GRASP-RN, GRASP-CN, GRASP-CN, TS-RN and SA methods stand out in terms of performance, while only TS-CN underperforms. As in the case of the large-sized instances, a summary consolidating the data of Tables ~\ref{tab9}, \ref{tab10}, and~\ref{tab11} can be seen in Table~\ref{tab17}. What is most interesting in the case of the small-sized instances is that proven optimal solutions for 173, out of the 180 instances, were reported in~\cite{arbro2023}. Table~\ref{tab12} shows the comparison with the optimal solutions in those 173 instances. In the table, for each method, we show the number of instances in which the method found an optimal solution (\#optimal). We also show the average gap in relation to the optimal solution when, for each method and instance, we consider the best solution of those found in the five runs of that method/instance pair (line ``gap to optimal''). Finally, in another line of the table we show the gap in relation to the optimal solution when, for a method/instance pair, we consider the average makespan of the five runs of that pair (line ``$\overline{\mathrm{gap}}$ to optimal''). The numbers in the table show that ILS-RN and ILS-CN found optimal solutions in \textit{all} 173 instances and, moreover, in all five runs of each method in each instance. Leaving aside TS-CN, the other methods also perform well, finding a large number of optimal solutions and with average gaps relative to the optimal solution of less than 1\%.

\begin{table}[ht!]
\begin{center}
\resizebox{!}{0.45\textheight}{
\begin{tabular}{|c|cc|cc|cc|cc|cc|cc|cc|}
\hline
\multirow{2}{*}{instance} & \multicolumn{2}{c|}{ILS-RN} & \multicolumn{2}{c|}{ILS-CN} & \multicolumn{2}{c|}{GRASP-RN} & \multicolumn{2}{c|}{GRASP-CN} & \multicolumn{2}{c|}{TS-RN} & \multicolumn{2}{c|}{TS-CN} & \multicolumn{2}{c|}{SA} \\
\cline{2-15} 
& $C_{\max}$ & Time & $C_{\max}$ & Time & $C_{\max}$ & Time & $C_{\max}$ & Time & $C_{\max}$ & Time & $C_{\max}$ & Time & $C_{\max}$ & Time \\
\hline
\hline
miniDAFJS01 & \textbf{22,875} & 0.01 & \textbf{22,875} & 0.05 & \textbf{22,875} & 8.09 & \textbf{22,875} & 16.21 & \textbf{22,875} & 0.00 & 23,264 & 0.00 & \textbf{22,875} & 0.40 \\
miniDAFJS02 & \textbf{22,708} & 0.00 & \textbf{22,708} & 0.00 & \textbf{22,708} & 0.03 & \textbf{22,708} & 0.04 & \textbf{22,708} & 0.00 & 23,242 & 0.00 & \textbf{22,708} & 0.14 \\
miniDAFJS03 & \textbf{18,363} & 0.00 & \textbf{18,363} & 0.00 & \textbf{18,363} & 0.00 & \textbf{18,363} & 0.00 & \textbf{18,363} & 0.00 & \textbf{18,363} & 0.00 & \textbf{18,363} & 0.00 \\
miniDAFJS04 & \textbf{20,498} & 0.00 & \textbf{20,498} & 0.00 & \textbf{20,498} & 0.00 & \textbf{20,498} & 0.00 & \textbf{20,498} & 0.00 & 21,172 & 0.00 & \textbf{20,498} & 0.10 \\
miniDAFJS05 & \textbf{20,593} & 0.01 & \textbf{20,593} & 0.00 & \textbf{20,593} & 2.01 & \textbf{20,593} & 6.26 & \textbf{20,593} & 0.01 & 21,324 & 0.00 & \textbf{20,593} & 0.32 \\
miniDAFJS06 & \textbf{22,867} & 0.01 & \textbf{22,867} & 0.00 & \textbf{22,867} & 0.69 & \textbf{22,867} & 1.04 & 23,370 & 0.00 & 23,370 & 0.00 & \textbf{22,867} & 0.33 \\
miniDAFJS07 & \textbf{25,715} & 0.01 & \textbf{25,715} & 0.01 & \textbf{25,715} & 0.15 & \textbf{25,715} & 0.31 & \textbf{25,715} & 0.00 & 26,978 & 0.00 & \textbf{25,715} & 0.18 \\
miniDAFJS08 & \textbf{19,878} & 0.00 & \textbf{19,878} & 0.00 & \textbf{19,878} & 0.00 & \textbf{19,878} & 0.00 & \textbf{19,878} & 0.00 & \textbf{19,878} & 0.00 & \textbf{19,878} & 0.00 \\
miniDAFJS09 & \textbf{24,267} & 0.03 & \textbf{24,267} & 0.16 & \textbf{24,267} & 1.30 & \textbf{24,267} & 8.19 & \textbf{24,267} & 0.00 & 24,304 & 0.00 & \textbf{24,267} & 0.41 \\
miniDAFJS10 & \textbf{20,336} & 0.00 & \textbf{20,336} & 0.00 & \textbf{20,336} & 0.00 & \textbf{20,336} & 0.00 & 20,359 & 0.00 & 20,873 & 0.00 & \textbf{20,336} & 0.24 \\
miniDAFJS11 & \textbf{29,968} & 0.01 & \textbf{29,968} & 0.02 & \textbf{29,968} & 1.11 & \textbf{29,968} & 0.00 & 33,263 & 0.00 & 33,689 & 0.00 & \textbf{29,968} & 0.42 \\
miniDAFJS12 & \textbf{18,670} & 0.01 & \textbf{18,670} & 0.01 & \textbf{18,670} & 0.00 & \textbf{18,670} & 2.87 & 20,222 & 0.00 & 20,342 & 0.00 & \textbf{18,670} & 0.26 \\
miniDAFJS13 & \textbf{16,313} & 0.00 & \textbf{16,313} & 0.00 & \textbf{16,313} & 0.00 & \textbf{16,313} & 0.00 & \textbf{16,313} & 0.00 & 17,091 & 0.00 & \textbf{16,313} & 0.01 \\
miniDAFJS14 & \textbf{23,140} & 0.00 & \textbf{23,140} & 0.00 & \textbf{23,140} & 0.00 & \textbf{23,140} & 3.36 & \textbf{23,140} & 0.00 & \textbf{23,140} & 0.00 & \textbf{23,140} & 0.30 \\
miniDAFJS15 & \textbf{21,715} & 0.00 & \textbf{21,715} & 0.00 & \textbf{21,715} & 0.00 & \textbf{21,715} & 0.24 & \textbf{21,715} & 0.00 & 21,882 & 0.00 & \textbf{21,715} & 0.16 \\
miniDAFJS16 & \textbf{25,426} & 0.00 & \textbf{25,426} & 0.00 & \textbf{25,426} & 0.00 & \textbf{25,426} & 0.00 & \textbf{25,426} & 0.00 & 25,691 & 0.00 & \textbf{25,426} & 0.26 \\
miniDAFJS17 & \textbf{20,155} & 0.00 & \textbf{20,155} & 0.00 & \textbf{20,155} & 0.00 & \textbf{20,155} & 0.00 & \textbf{20,155} & 0.00 & \textbf{20,155} & 0.00 & \textbf{20,155} & 0.19 \\
miniDAFJS18 & \textbf{18,135} & 0.00 & \textbf{18,135} & 0.00 & \textbf{18,135} & 0.00 & \textbf{18,135} & 0.00 & \textbf{18,135} & 0.00 & 18,445 & 0.00 & \textbf{18,135} & 0.25 \\
miniDAFJS19 & \textbf{20,945} & 0.00 & \textbf{20,945} & 0.00 & \textbf{20,945} & 0.00 & \textbf{20,945} & 0.00 & \textbf{20,945} & 0.00 & 21,293 & 0.00 & \textbf{20,945} & 0.18 \\
miniDAFJS20 & \textbf{21,838} & 0.02 & \textbf{21,838} & 0.01 & \textbf{21,838} & 12.83 & \textbf{21,838} & 0.60 & \textbf{21,838} & 0.00 & 22,675 & 0.00 & \textbf{21,838} & 0.30 \\
miniDAFJS21 & \textbf{23,344} & 0.60 & \textbf{23,344} & 0.14 & \textbf{23,344} & 50.28 & \textbf{23,344} & 74.78 & 23,456 & 0.02 & 24,554 & 0.00 & 23,563 & 35.44 \\
miniDAFJS22 & \textbf{25,923} & 0.00 & \textbf{25,923} & 0.00 & \textbf{25,923} & 0.00 & \textbf{25,923} & 0.03 & \textbf{25,923} & 0.00 & \textbf{25,923} & 0.00 & \textbf{25,923} & 0.41 \\
miniDAFJS23 & \textbf{24,038} & 3.23 & \textbf{24,038} & 0.85 & 24,998 & 0.03 & 24,998 & 0.01 & 24,585 & 0.11 & 24,438 & 0.00 & 24,438 & 0.49 \\
miniDAFJS24 & \textbf{24,579} & 0.04 & \textbf{24,579} & 0.01 & \textbf{24,579} & 0.01 & \textbf{24,579} & 0.02 & 25,521 & 0.01 & 25,694 & 0.00 & \textbf{24,579} & 0.43 \\
miniDAFJS25 & \textbf{21,143} & 0.01 & \textbf{21,143} & 0.00 & \textbf{21,143} & 23.97 & \textbf{21,143} & 0.00 & \textbf{21,143} & 0.00 & 22,682 & 0.00 & \textbf{21,143} & 0.33 \\
miniDAFJS26 & \textbf{21,120} & 0.16 & \textbf{21,120} & 0.11 & \textbf{21,120} & 42.71 & \textbf{21,120} & 43.64 & \textbf{21,120} & 0.00 & 21,468 & 0.00 & \textbf{21,120} & 0.57 \\
miniDAFJS27 & \textbf{22,050} & 0.14 & \textbf{22,050} & 0.20 & 22,106 & 0.12 & \textbf{22,050} & 0.24 & 22,106 & 0.03 & 22,106 & 0.00 & 22,106 & 0.41 \\
miniDAFJS28 & \textbf{22,708} & 0.03 & \textbf{22,708} & 0.02 & \textbf{22,708} & 0.00 & \textbf{22,708} & 0.00 & \textbf{22,708} & 0.00 & 23,120 & 0.00 & \textbf{22,708} & 0.33 \\
miniDAFJS29 & \textbf{20,278} & 0.01 & \textbf{20,278} & 0.01 & \textbf{20,278} & 0.01 & \textbf{20,278} & 0.00 & \textbf{20,278} & 0.00 & \textbf{20,278} & 0.00 & \textbf{20,278} & 0.33 \\
miniDAFJS30 & \textbf{23,558} & 0.05 & \textbf{23,558} & 0.03 & \textbf{23,558} & 33.60 & \textbf{23,558} & 48.64 & \textbf{23,558} & 0.09 & \textbf{23,558} & 0.01 & \textbf{23,558} & 0.44 \\
\hline 
$\bar{C}_{\max}$ & \multicolumn{2}{c|}{22,104.87} & \multicolumn{2}{c|}{22,104.87} & \multicolumn{2}{c|}{22,138.73} & \multicolumn{2}{c|}{22,136.87} & \multicolumn{2}{c|}{22,339.20} & \multicolumn{2}{c|}{22,699.73} & \multicolumn{2}{c|}{22,127.37} \\
\#best & \multicolumn{2}{c|}{30} & \multicolumn{2}{c|}{30} & \multicolumn{2}{c|}{28} & \multicolumn{2}{c|}{29} & \multicolumn{2}{c|}{22} & \multicolumn{2}{c|}{7} & \multicolumn{2}{c|}{27} \\
gap(\%) & \multicolumn{2}{c|}{0.00} & \multicolumn{2}{c|}{0.00} & \multicolumn{2}{c|}{0.14} & \multicolumn{2}{c|}{0.13} & \multicolumn{2}{c|}{0.95} & \multicolumn{2}{c|}{2.61} & \multicolumn{2}{c|}{0.10} \\
$\overline{\mathrm{gap}}$(\%) & \multicolumn{2}{c|}{0.00} & \multicolumn{2}{c|}{0.00} & \multicolumn{2}{c|}{0.18} & \multicolumn{2}{c|}{0.17} & \multicolumn{2}{c|}{0.95} & \multicolumn{2}{c|}{2.61} & \multicolumn{2}{c|}{0.35} \\
\hline
\multicolumn{15}{c}{}\\
\hline
\multirow{2}{*}{instance} & \multicolumn{2}{c|}{ILS-RN} & \multicolumn{2}{c|}{ILS-CN} & \multicolumn{2}{c|}{GRASP-RN} & \multicolumn{2}{c|}{GRASP-CN} & \multicolumn{2}{c|}{TS-RN} & \multicolumn{2}{c|}{TS-CN} & \multicolumn{2}{c|}{SA} \\
\cline{2-15} 
& $C_{\max}$ & Time & $C_{\max}$ & Time & $C_{\max}$ & Time & $C_{\max}$ & Time & $C_{\max}$ & Time & $C_{\max}$ & Time & $C_{\max}$ & Time \\
\hline
\hline
miniYFJS01 & \textbf{35,046} & 0.00 & \textbf{35,046} & 0.00 & \textbf{35,046} & 0.00 & \textbf{35,046} & 0.00 & \textbf{35,046} & 0.00 & 35,243 & 0.00 & \textbf{35,046} & 0.20 \\
miniYFJS02 & \textbf{24,359} & 0.01 & \textbf{24,359} & 0.01 & \textbf{24,359} & 0.00 & \textbf{24,359} & 0.00 & \textbf{24,359} & 0.00 & 25,654 & 0.00 & \textbf{24,359} & 0.29 \\
miniYFJS03 & \textbf{47,391} & 0.00 & \textbf{47,391} & 0.01 & \textbf{47,391} & 0.00 & \textbf{47,391} & 0.14 & 52,098 & 0.00 & 53,381 & 0.00 & \textbf{47,391} & 0.35 \\
miniYFJS04 & \textbf{25,394} & 0.00 & \textbf{25,394} & 0.00 & \textbf{25,394} & 0.00 & \textbf{25,394} & 0.00 & \textbf{25,394} & 0.00 & \textbf{25,394} & 0.00 & \textbf{25,394} & 0.00 \\
miniYFJS05 & \textbf{23,985} & 0.00 & \textbf{23,985} & 0.00 & \textbf{23,985} & 0.00 & \textbf{23,985} & 0.00 & \textbf{23,985} & 0.00 & 26,371 & 0.00 & \textbf{23,985} & 0.24 \\
miniYFJS06 & \textbf{29,469} & 0.01 & \textbf{29,469} & 0.01 & \textbf{29,469} & 0.35 & \textbf{29,469} & 0.25 & 30,508 & 0.00 & 30,952 & 0.00 & \textbf{29,469} & 0.30 \\
miniYFJS07 & \textbf{45,705} & 0.01 & \textbf{45,705} & 0.05 & \textbf{45,705} & 0.01 & \textbf{45,705} & 14.81 & \textbf{45,705} & 0.01 & 47,524 & 0.00 & \textbf{45,705} & 0.43 \\
miniYFJS08 & \textbf{33,829} & 0.01 & \textbf{33,829} & 0.03 & \textbf{33,829} & 0.00 & \textbf{33,829} & 0.16 & 33,883 & 0.00 & 34,338 & 0.00 & \textbf{33,829} & 0.38 \\
miniYFJS09 & \textbf{37,049} & 0.00 & \textbf{37,049} & 0.00 & \textbf{37,049} & 0.00 & \textbf{37,049} & 0.00 & \textbf{37,049} & 0.00 & 37,354 & 0.00 & \textbf{37,049} & 0.30 \\
miniYFJS10 & \textbf{27,310} & 0.01 & \textbf{27,310} & 0.00 & \textbf{27,310} & 0.01 & \textbf{27,310} & 0.00 & \textbf{27,310} & 0.00 & \textbf{27,310} & 0.00 & \textbf{27,310} & 0.33 \\
miniYFJS11 & \textbf{41,300} & 0.00 & \textbf{41,300} & 0.01 & \textbf{41,300} & 0.00 & \textbf{41,300} & 0.05 & \textbf{41,300} & 0.00 & 41,423 & 0.00 & \textbf{41,300} & 0.39 \\
miniYFJS12 & \textbf{30,145} & 0.11 & \textbf{30,145} & 0.06 & \textbf{30,145} & 55.26 & \textbf{30,145} & 64.48 & \textbf{30,145} & 0.03 & 33,329 & 0.00 & \textbf{30,145} & 0.54 \\
miniYFJS13 & \textbf{30,962} & 0.01 & \textbf{30,962} & 0.00 & \textbf{30,962} & 0.00 & \textbf{30,962} & 0.01 & \textbf{30,962} & 0.01 & \textbf{30,962} & 0.00 & \textbf{30,962} & 0.35 \\
miniYFJS14 & \textbf{31,398} & 0.01 & \textbf{31,398} & 0.01 & \textbf{31,398} & 0.00 & \textbf{31,398} & 0.00 & \textbf{31,398} & 0.00 & \textbf{31,398} & 0.00 & \textbf{31,398} & 0.38 \\
miniYFJS15 & \textbf{45,442} & 0.00 & \textbf{45,442} & 0.00 & \textbf{45,442} & 0.00 & \textbf{45,442} & 0.00 & \textbf{45,442} & 0.00 & \textbf{45,442} & 0.00 & \textbf{45,442} & 0.39 \\
miniYFJS16 & \textbf{33,791} & 0.00 & \textbf{33,791} & 0.02 & \textbf{33,791} & 0.00 & \textbf{33,791} & 19.85 & \textbf{33,791} & 0.00 & \textbf{33,791} & 0.00 & \textbf{33,791} & 0.41 \\
miniYFJS17 & \textbf{42,838} & 0.02 & \textbf{42,838} & 0.02 & \textbf{42,838} & 0.00 & \textbf{42,838} & 13.56 & \textbf{42,838} & 0.00 & 44,181 & 0.00 & \textbf{42,838} & 0.46 \\
miniYFJS18 & \textbf{28,247} & 0.01 & \textbf{28,247} & 0.01 & \textbf{28,247} & 0.01 & \textbf{28,247} & 0.00 & \textbf{28,247} & 0.00 & \textbf{28,247} & 0.00 & \textbf{28,247} & 0.41 \\
miniYFJS19 & \textbf{33,601} & 0.04 & \textbf{33,601} & 0.05 & \textbf{33,601} & 0.00 & \textbf{33,601} & 0.01 & \textbf{33,601} & 0.01 & 36,706 & 0.00 & \textbf{33,601} & 0.45 \\
miniYFJS20 & \textbf{30,837} & 0.00 & \textbf{30,837} & 0.00 & \textbf{30,837} & 0.00 & \textbf{30,837} & 0.00 & \textbf{30,837} & 0.02 & 33,385 & 0.00 & \textbf{30,837} & 0.41 \\
miniYFJS21 & \textbf{37,096} & 0.74 & \textbf{37,096} & 0.34 & 37,308 & 0.03 & 37,275 & 0.23 & \textbf{37,096} & 0.41 & 37,275 & 0.02 & 37,308 & 70.41 \\
miniYFJS22 & \textbf{34,282} & 0.03 & \textbf{34,282} & 0.03 & \textbf{34,282} & 0.01 & \textbf{34,282} & 0.08 & \textbf{34,282} & 0.00 & 34,917 & 0.00 & \textbf{34,282} & 0.48 \\
miniYFJS23 & \textbf{42,079} & 0.04 & \textbf{42,079} & 0.02 & \textbf{42,079} & 0.07 & \textbf{42,079} & 38.95 & \textbf{42,079} & 0.01 & \textbf{42,079} & 0.00 & \textbf{42,079} & 0.51 \\
miniYFJS24 & \textbf{30,905} & 0.11 & \textbf{30,905} & 0.08 & 30,965 & 69.07 & \textbf{30,905} & 53.76 & \textbf{30,905} & 0.14 & 33,144 & 0.00 & \textbf{30,905} & 0.56 \\
miniYFJS25 & \textbf{36,170} & 0.10 & \textbf{36,170} & 0.13 & \textbf{36,170} & 76.58 & \textbf{36,170} & 1.25 & \textbf{36,170} & 0.11 & 37,251 & 0.00 & \textbf{36,170} & 0.94 \\
miniYFJS26 & \textbf{51,466} & 0.01 & \textbf{51,466} & 0.02 & \textbf{51,466} & 0.83 & \textbf{51,466} & 1.23 & 54,534 & 0.03 & 51,793 & 0.00 & \textbf{51,466} & 0.49 \\
miniYFJS27 & \textbf{36,719} & 0.01 & \textbf{36,719} & 0.01 & \textbf{36,719} & 0.00 & \textbf{36,719} & 0.02 & \textbf{36,719} & 0.01 & 36,823 & 0.00 & \textbf{36,719} & 0.47 \\
miniYFJS28 & \textbf{34,509} & 0.04 & \textbf{34,509} & 0.03 & \textbf{34,509} & 87.68 & \textbf{34,509} & 52.00 & \textbf{34,509} & 0.04 & \textbf{34,509} & 0.02 & \textbf{34,509} & 31.22 \\
miniYFJS29 & \textbf{39,798} & 0.06 & \textbf{39,798} & 0.14 & \textbf{39,798} & 0.07 & \textbf{39,798} & 0.05 & 39,949 & 0.01 & 40,891 & 0.00 & \textbf{39,798} & 0.63 \\
miniYFJS30 & \textbf{33,974} & 0.03 & \textbf{33,974} & 0.06 & \textbf{33,974} & 0.05 & \textbf{33,974} & 0.08 & \textbf{33,974} & 0.12 & \textbf{33,974} & 0.00 & \textbf{33,974} & 0.52 \\
\hline
$\bar{C}_{\max}$ & \multicolumn{2}{c|}{35,169.87} & \multicolumn{2}{c|}{35,169.87} & \multicolumn{2}{c|}{35,178.93} & \multicolumn{2}{c|}{35,175.83} & \multicolumn{2}{c|}{35,470.50} & \multicolumn{2}{c|}{36,168.03} & \multicolumn{2}{c|}{35,176.93} \\
\#best & \multicolumn{2}{c|}{30} & \multicolumn{2}{c|}{30} & \multicolumn{2}{c|}{28} & \multicolumn{2}{c|}{29} & \multicolumn{2}{c|}{25} & \multicolumn{2}{c|}{10} & \multicolumn{2}{c|}{29} \\
gap(\%) & \multicolumn{2}{c|}{0.00} & \multicolumn{2}{c|}{0.00} & \multicolumn{2}{c|}{0.03} & \multicolumn{2}{c|}{0.02} & \multicolumn{2}{c|}{0.67} & \multicolumn{2}{c|}{2.92} & \multicolumn{2}{c|}{0.02} \\
$\overline{\mathrm{gap}}$(\%) & \multicolumn{2}{c|}{0.00} & \multicolumn{2}{c|}{0.00} & \multicolumn{2}{c|}{0.39} & \multicolumn{2}{c|}{0.06} & \multicolumn{2}{c|}{0.67} & \multicolumn{2}{c|}{2.92} & \multicolumn{2}{c|}{0.18} \\
\hline
\end{tabular}}
\end{center}
\caption{Results of applying the metaheuristics to the small-sized instances with learning rate $\alpha = 0.1$.}
\label{tab9}
\end{table}

\begin{table}[ht!]
\begin{center}
\resizebox{!}{0.45\textheight}{
\begin{tabular}{|c|cc|cc|cc|cc|cc|cc|cc|}
\hline
\multirow{2}{*}{instance} & \multicolumn{2}{c|}{ILS-RN} & \multicolumn{2}{c|}{ILS-CN} & \multicolumn{2}{c|}{GRASP-RN} & \multicolumn{2}{c|}{GRASP-CN} & \multicolumn{2}{c|}{TS-RN} & \multicolumn{2}{c|}{TS-CN} & \multicolumn{2}{c|}{SA} \\
\cline{2-15} 
& $C_{\max}$ & Time & $C_{\max}$ & Time & $C_{\max}$ & Time & $C_{\max}$ & Time & $C_{\max}$ & Time & $C_{\max}$ & Time & $C_{\max}$ & Time \\
\hline
\hline
miniDAFJS01 & \textbf{21,327} & 0.02 & \textbf{21,327} & 0.03 & \textbf{21,327} & 0.03 & \textbf{21,327} & 28.07 & 21,546 & 0.00 & 21,546 & 0.00 & \textbf{21,327} & 0.33 \\
miniDAFJS02 & \textbf{20,635} & 0.01 & \textbf{20,635} & 0.01 & \textbf{20,635} & 0.05 & \textbf{20,635} & 0.04 & \textbf{20,635} & 0.00 & 22,161 & 0.00 & \textbf{20,635} & 0.19 \\
miniDAFJS03 & \textbf{17,972} & 0.00 & \textbf{17,972} & 0.00 & \textbf{17,972} & 0.00 & \textbf{17,972} & 0.00 & \textbf{17,972} & 0.00 & \textbf{17,972} & 0.00 & \textbf{17,972} & 0.00 \\
miniDAFJS04 & \textbf{19,602} & 0.00 & \textbf{19,602} & 0.00 & \textbf{19,602} & 0.00 & \textbf{19,602} & 0.00 & \textbf{19,602} & 0.00 & 20,691 & 0.00 & \textbf{19,602} & 0.11 \\
miniDAFJS05 & \textbf{18,803} & 0.00 & \textbf{18,803} & 0.01 & \textbf{18,803} & 3.57 & \textbf{18,803} & 4.21 & \textbf{18,803} & 0.01 & \textbf{18,803} & 0.00 & \textbf{18,803} & 0.28 \\
miniDAFJS06 & \textbf{20,568} & 0.04 & \textbf{20,568} & 0.04 & \textbf{20,568} & 0.01 & \textbf{20,568} & 0.01 & 20,783 & 0.00 & 20,783 & 0.00 & \textbf{20,568} & 0.31 \\
miniDAFJS07 & \textbf{24,715} & 0.00 & \textbf{24,715} & 0.01 & \textbf{24,715} & 0.10 & \textbf{24,715} & 0.63 & \textbf{24,715} & 0.00 & \textbf{24,715} & 0.00 & \textbf{24,715} & 0.22 \\
miniDAFJS08 & \textbf{18,857} & 0.00 & \textbf{18,857} & 0.00 & \textbf{18,857} & 0.00 & \textbf{18,857} & 0.00 & \textbf{18,857} & 0.00 & \textbf{18,857} & 0.00 & \textbf{18,857} & 0.00 \\
miniDAFJS09 & \textbf{22,660} & 0.02 & \textbf{22,660} & 0.05 & \textbf{22,660} & 4.57 & \textbf{22,660} & 2.76 & 22,730 & 0.01 & 22,730 & 0.00 & \textbf{22,660} & 0.31 \\
miniDAFJS10 & \textbf{18,823} & 0.00 & \textbf{18,823} & 0.00 & \textbf{18,823} & 0.00 & \textbf{18,823} & 0.00 & \textbf{18,823} & 0.00 & 20,180 & 0.00 & \textbf{18,823} & 0.18 \\
miniDAFJS11 & \textbf{27,455} & 0.02 & \textbf{27,455} & 0.01 & \textbf{27,455} & 2.05 & \textbf{27,455} & 0.00 & 30,550 & 0.00 & 30,941 & 0.00 & \textbf{27,455} & 0.35 \\
miniDAFJS12 & \textbf{17,874} & 0.01 & \textbf{17,874} & 0.01 & \textbf{17,874} & 0.92 & \textbf{17,874} & 1.47 & 19,624 & 0.00 & 18,974 & 0.00 & \textbf{17,874} & 0.24 \\
miniDAFJS13 & \textbf{15,143} & 0.00 & \textbf{15,143} & 0.00 & \textbf{15,143} & 0.00 & \textbf{15,143} & 0.00 & \textbf{15,143} & 0.00 & 16,593 & 0.00 & \textbf{15,143} & 0.01 \\
miniDAFJS14 & \textbf{21,817} & 0.00 & \textbf{21,817} & 0.00 & \textbf{21,817} & 0.00 & \textbf{21,817} & 3.31 & \textbf{21,817} & 0.00 & 21,990 & 0.00 & \textbf{21,817} & 0.27 \\
miniDAFJS15 & \textbf{20,236} & 0.00 & \textbf{20,236} & 0.00 & \textbf{20,236} & 0.00 & \textbf{20,236} & 0.00 & \textbf{20,236} & 0.00 & \textbf{20,236} & 0.00 & \textbf{20,236} & 0.14 \\
miniDAFJS16 & \textbf{24,114} & 0.00 & \textbf{24,114} & 0.00 & \textbf{24,114} & 0.00 & \textbf{24,114} & 0.01 & \textbf{24,114} & 0.00 & 24,593 & 0.00 & \textbf{24,114} & 0.27 \\
miniDAFJS17 & \textbf{19,145} & 0.00 & \textbf{19,145} & 0.00 & \textbf{19,145} & 0.00 & \textbf{19,145} & 0.00 & \textbf{19,145} & 0.00 & \textbf{19,145} & 0.00 & \textbf{19,145} & 0.15 \\
miniDAFJS18$^+$ & \textbf{17,270} & 0.00 & \textbf{17,270} & 0.00 & \textbf{17,270} & 0.00 & \textbf{17,270} & 0.00 & \textbf{17,270} & 0.00 & 17,966 & 0.00 & \textbf{17,270} & 0.25 \\
miniDAFJS19 & \textbf{19,642} & 0.00 & \textbf{19,642} & 0.00 & \textbf{19,642} & 0.00 & \textbf{19,642} & 0.00 & \textbf{19,642} & 0.00 & 20,107 & 0.00 & \textbf{19,642} & 0.16 \\
miniDAFJS20 & \textbf{20,086} & 0.04 & \textbf{20,086} & 0.01 & \textbf{20,086} & 24.76 & \textbf{20,086} & 1.87 & 21,286 & 0.00 & 20,646 & 0.00 & \textbf{20,086} & 0.28 \\
miniDAFJS21 & \textbf{21,352} & 0.21 & \textbf{21,352} & 0.92 & 21,433 & 43.93 & 21,536 & 93.99 & \textbf{21,352} & 0.04 & 21,783 & 0.01 & 21,433 & 1.51 \\
miniDAFJS22 & \textbf{23,852} & 0.01 & \textbf{23,852} & 0.02 & \textbf{23,852} & 38.91 & \textbf{23,852} & 24.30 & 24,273 & 0.00 & 24,273 & 0.00 & \textbf{23,852} & 0.41 \\
miniDAFJS23$^+$ & \textbf{22,390} & 2.60 & \textbf{22,390} & 0.82 & 23,228 & 0.02 & 23,300 & 0.02 & \textbf{22,390} & 0.31 & 22,491 & 0.02 & 23,120 & 0.49 \\
miniDAFJS24 & \textbf{22,521} & 0.07 & \textbf{22,521} & 0.07 & \textbf{22,521} & 19.38 & \textbf{22,521} & 44.36 & \textbf{22,521} & 0.03 & 22,551 & 0.00 & \textbf{22,521} & 0.40 \\
miniDAFJS25 & \textbf{19,809} & 0.01 & \textbf{19,809} & 0.03 & \textbf{19,809} & 38.66 & \textbf{19,809} & 0.00 & 19,913 & 0.01 & 19,913 & 0.00 & \textbf{19,809} & 0.35 \\
miniDAFJS26 & \textbf{19,724} & 0.20 & \textbf{19,724} & 0.23 & \textbf{19,724} & 24.84 & \textbf{19,724} & 48.07 & \textbf{19,724} & 0.00 & 19,994 & 0.00 & 19,773 & 0.33 \\
miniDAFJS27 & \textbf{20,245} & 0.13 & \textbf{20,245} & 0.03 & \textbf{20,245} & 0.07 & \textbf{20,245} & 35.36 & \textbf{20,245} & 0.02 & 20,646 & 0.00 & \textbf{20,245} & 0.47 \\
miniDAFJS28 & \textbf{20,635} & 0.02 & \textbf{20,635} & 0.04 & \textbf{20,635} & 0.00 & \textbf{20,635} & 0.00 & \textbf{20,635} & 0.00 & 21,693 & 0.00 & \textbf{20,635} & 0.38 \\
miniDAFJS29 & \textbf{19,201} & 0.01 & \textbf{19,201} & 0.00 & \textbf{19,201} & 1.39 & \textbf{19,201} & 0.00 & \textbf{19,201} & 0.00 & 19,448 & 0.00 & \textbf{19,201} & 0.31 \\
miniDAFJS30 & \textbf{21,552} & 0.06 & \textbf{21,552} & 0.21 & 21,698 & 46.08 & \textbf{21,552} & 63.33 & \textbf{21,552} & 0.15 & 21,645 & 0.00 & \textbf{21,552} & 0.45 \\ \hline
$\bar{C}_{\max}$ & \multicolumn{2}{c|}{20,600.83} & \multicolumn{2}{c|}{20,600.83} & \multicolumn{2}{c|}{20,636.33} & \multicolumn{2}{c|}{20,637.30} & \multicolumn{2}{c|}{20,836.63} & \multicolumn{2}{c|}{21,135.53} & \multicolumn{2}{c|}{20,629.50} \\
\#best & \multicolumn{2}{c|}{30} & \multicolumn{2}{c|}{30} & \multicolumn{2}{c|}{27} & \multicolumn{2}{c|}{28} & \multicolumn{2}{c|}{22} & \multicolumn{2}{c|}{6} & \multicolumn{2}{c|}{27} \\
gap(\%) & \multicolumn{2}{c|}{0.00} & \multicolumn{2}{c|}{0.00} & \multicolumn{2}{c|}{0.16} & \multicolumn{2}{c|}{0.16} & \multicolumn{2}{c|}{1.06} & \multicolumn{2}{c|}{2.60} & \multicolumn{2}{c|}{0.13} \\
$\overline{\mathrm{gap}}$(\%) & \multicolumn{2}{c|}{0.00} & \multicolumn{2}{c|}{0.00} & \multicolumn{2}{c|}{0.24} & \multicolumn{2}{c|}{0.28} & \multicolumn{2}{c|}{1.06} & \multicolumn{2}{c|}{2.60} & \multicolumn{2}{c|}{0.42} \\ 
\hline
\multicolumn{15}{c}{}\\
\hline
\multirow{2}{*}{instance} & \multicolumn{2}{c|}{ILS--RN} & \multicolumn{2}{c|}{ILS--CN} & \multicolumn{2}{c|}{GRASP--RN} & \multicolumn{2}{c|}{GRASP--CN} & \multicolumn{2}{c|}{TS--RN} & \multicolumn{2}{c|}{TS--CN} & \multicolumn{2}{c|}{SA} \\
\cline{2-15} 
& $C_{\max}$ & time & $C_{\max}$ & time & $C_{\max}$ & time & $C_{\max}$ & time & $C_{\max}$ & time & $C_{\max}$ & time & $C_{\max}$ & time \\
\hline
\hline
miniYFJS01 & \textbf{33,132} & 0.00 & \textbf{33,132} & 0.00 & \textbf{33,132} & 0.00 & \textbf{33,132} & 4.67 & \textbf{33,132} & 0.00 & 34,443 & 0.00 & \textbf{33,132} & 0.28 \\
miniYFJS02 & \textbf{23,100} & 0.00 & \textbf{23,100} & 0.01 & \textbf{23,100} & 0.00 & \textbf{23,100} & 0.00 & \textbf{23,100} & 0.00 & 25,032 & 0.00 & \textbf{23,100} & 0.30 \\
miniYFJS03 & \textbf{42,896} & 0.00 & \textbf{42,896} & 0.01 & \textbf{42,896} & 0.00 & \textbf{42,896} & 0.01 & 46,806 & 0.00 & 49,111 & 0.00 & \textbf{42,896} & 0.34 \\
miniYFJS04 & \textbf{24,485} & 0.00 & \textbf{24,485} & 0.00 & \textbf{24,485} & 0.00 & \textbf{24,485} & 0.00 & \textbf{24,485} & 0.00 & \textbf{24,485} & 0.00 & \textbf{24,485} & 0.31 \\
miniYFJS05 & \textbf{23,597} & 0.00 & \textbf{23,597} & 0.00 & \textbf{23,597} & 0.00 & \textbf{23,597} & 0.00 & \textbf{23,597} & 0.00 & \textbf{23,597} & 0.00 & \textbf{23,597} & 0.25 \\
miniYFJS06 & \textbf{28,655} & 0.01 & \textbf{28,655} & 0.01 & \textbf{28,655} & 0.00 & \textbf{28,655} & 0.00 & \textbf{28,655} & 0.01 & 28,750 & 0.01 & \textbf{28,655} & 0.32 \\
miniYFJS07 & \textbf{42,239} & 0.00 & \textbf{42,239} & 0.05 & \textbf{42,239} & 0.01 & \textbf{42,239} & 7.75 & \textbf{42,239} & 0.00 & 45,493 & 0.00 & \textbf{42,239} & 0.40 \\
miniYFJS08 & \textbf{31,471} & 0.02 & \textbf{31,471} & 0.15 & \textbf{31,471} & 0.00 & \textbf{31,471} & 22.11 & 31,597 & 0.00 & 32,669 & 0.00 & \textbf{31,471} & 0.42 \\
miniYFJS09 & \textbf{35,250} & 0.02 & \textbf{35,250} & 0.01 & \textbf{35,250} & 4.98 & \textbf{35,250} & 0.01 & 35,794 & 0.00 & 36,098 & 0.00 & \textbf{35,250} & 0.32 \\
miniYFJS10 & \textbf{26,145} & 0.01 & \textbf{26,145} & 0.00 & \textbf{26,145} & 5.41 & \textbf{26,145} & 0.00 & \textbf{26,145} & 0.00 & \textbf{26,145} & 0.00 & \textbf{26,145} & 0.32 \\
miniYFJS11 & \textbf{38,545} & 0.01 & \textbf{38,545} & 0.02 & \textbf{38,545} & 3.71 & \textbf{38,545} & 7.90 & \textbf{38,545} & 0.00 & 38,756 & 0.00 & \textbf{38,545} & 0.44 \\
miniYFJS12 & \textbf{27,895} & 0.15 & \textbf{27,895} & 0.09 & \textbf{27,895} & 48.02 & 28,281 & 0.01 & \textbf{27,895} & 0.10 & 30,851 & 0.00 & \textbf{27,895} & 7.80 \\
miniYFJS13 & \textbf{28,120} & 0.01 & \textbf{28,120} & 0.00 & \textbf{28,120} & 0.00 & \textbf{28,120} & 0.02 & \textbf{28,120} & 0.00 & \textbf{28,120} & 0.00 & \textbf{28,120} & 0.38 \\
miniYFJS14 & \textbf{29,682} & 0.02 & \textbf{29,682} & 0.00 & \textbf{29,682} & 0.00 & \textbf{29,682} & 0.00 & \textbf{29,682} & 0.00 & \textbf{29,682} & 0.00 & \textbf{29,682} & 0.37 \\
miniYFJS15 & \textbf{41,619} & 0.00 & \textbf{41,619} & 0.00 & \textbf{41,619} & 0.00 & \textbf{41,619} & 0.01 & \textbf{41,619} & 0.00 & 42,415 & 0.00 & \textbf{41,619} & 0.41 \\
miniYFJS16 & \textbf{31,280} & 0.02 & \textbf{31,280} & 0.02 & \textbf{31,280} & 0.00 & \textbf{31,280} & 0.00 & \textbf{31,280} & 0.05 & 32,366 & 0.00 & \textbf{31,280} & 0.48 \\
miniYFJS17 & \textbf{40,388} & 0.01 & \textbf{40,388} & 0.03 & \textbf{40,388} & 0.00 & \textbf{40,388} & 10.05 & \textbf{40,388} & 0.01 & 41,316 & 0.00 & \textbf{40,388} & 0.44 \\
miniYFJS18 & \textbf{26,297} & 0.00 & \textbf{26,297} & 0.01 & \textbf{26,297} & 0.00 & \textbf{26,297} & 0.01 & \textbf{26,297} & 0.01 & \textbf{26,297} & 0.00 & \textbf{26,297} & 0.38 \\
miniYFJS19 & \textbf{30,717} & 0.02 & \textbf{30,717} & 0.06 & \textbf{30,717} & 0.01 & \textbf{30,717} & 0.01 & \textbf{30,717} & 0.02 & 33,965 & 0.00 & \textbf{30,717} & 0.47 \\
miniYFJS20 & \textbf{28,832} & 0.03 & \textbf{28,832} & 0.02 & \textbf{28,832} & 0.04 & \textbf{28,832} & 0.03 & \textbf{28,832} & 0.03 & 31,304 & 0.00 & \textbf{28,832} & 0.38 \\
miniYFJS21 & \textbf{34,811} & 0.13 & \textbf{34,811} & 0.09 & \textbf{34,811} & 0.07 & \textbf{34,811} & 0.06 & 34,876 & 0.00 & \textbf{34,811} & 0.05 & \textbf{34,811} & 24.18 \\
miniYFJS22 & \textbf{31,702} & 0.01 & \textbf{31,702} & 0.03 & \textbf{31,702} & 0.03 & 31,747 & 0.17 & \textbf{31,702} & 0.03 & 31,946 & 0.00 & \textbf{31,702} & 0.54 \\
miniYFJS23 & \textbf{38,639} & 0.08 & \textbf{38,639} & 0.14 & \textbf{38,639} & 33.88 & \textbf{38,639} & 39.04 & \textbf{38,639} & 0.04 & 39,721 & 0.00 & \textbf{38,639} & 0.50 \\
miniYFJS24 & \textbf{28,884} & 0.05 & \textbf{28,884} & 0.05 & 30,200 & 21.84 & \textbf{28,884} & 46.64 & \textbf{28,884} & 0.06 & 30,568 & 0.00 & \textbf{28,884} & 0.53 \\
miniYFJS25 & \textbf{34,231} & 0.04 & \textbf{34,231} & 0.04 & \textbf{34,231} & 26.59 & \textbf{34,231} & 41.59 & \textbf{34,231} & 0.01 & 35,256 & 0.00 & \textbf{34,231} & 0.46 \\
miniYFJS26 & \textbf{47,519} & 0.01 & \textbf{47,519} & 0.01 & \textbf{47,519} & 1.48 & \textbf{47,519} & 5.00 & 49,202 & 0.04 & 48,066 & 0.00 & \textbf{47,519} & 0.50 \\
miniYFJS27 & \textbf{34,042} & 0.08 & \textbf{34,042} & 0.09 & \textbf{34,042} & 0.03 & \textbf{34,042} & 1.61 & \textbf{34,042} & 0.01 & 34,225 & 0.00 & \textbf{34,042} & 0.55 \\
miniYFJS28 & \textbf{32,080} & 0.04 & \textbf{32,080} & 0.04 & \textbf{32,080} & 69.04 & \textbf{32,080} & 43.21 & \textbf{32,080} & 0.02 & 32,852 & 0.00 & \textbf{32,080} & 0.60 \\
miniYFJS29 & \textbf{36,093} & 0.11 & \textbf{36,093} & 0.31 & \textbf{36,093} & 0.13 & \textbf{36,093} & 0.03 & 36,561 & 0.00 & 36,561 & 0.00 & \textbf{36,093} & 0.47 \\
miniYFJS30 & \textbf{31,888} & 0.53 & \textbf{31,888} & 0.42 & \textbf{31,888} & 59.15 & \textbf{31,888} & 46.43 & \textbf{31,888} & 0.26 & 33,118 & 0.00 & \textbf{31,888} & 1.15 \\ \hline
$\bar{C}_{\max}$ & \multicolumn{2}{c|}{32,807.80} & \multicolumn{2}{c|}{32,807.80} & \multicolumn{2}{c|}{32,851.67} & \multicolumn{2}{c|}{32,822.17} & \multicolumn{2}{c|}{33,034.33} & \multicolumn{2}{c|}{33,933.97} & \multicolumn{2}{c|}{32,807.80} \\
\#best & \multicolumn{2}{c|}{30} & \multicolumn{2}{c|}{30} & \multicolumn{2}{c|}{29} & \multicolumn{2}{c|}{28} & \multicolumn{2}{c|}{24} & \multicolumn{2}{c|}{7} & \multicolumn{2}{c|}{30} \\
gap(\%) & \multicolumn{2}{c|}{0.00} & \multicolumn{2}{c|}{0.00} & \multicolumn{2}{c|}{0.15} & \multicolumn{2}{c|}{0.05} & \multicolumn{2}{c|}{0.54} & \multicolumn{2}{c|}{3.36} & \multicolumn{2}{c|}{0.00} \\
$\overline{\mathrm{gap}}$(\%) & \multicolumn{2}{c|}{0.00} & \multicolumn{2}{c|}{0.00} & \multicolumn{2}{c|}{0.43} & \multicolumn{2}{c|}{0.17} & \multicolumn{2}{c|}{0.54} & \multicolumn{2}{c|}{3.36} & \multicolumn{2}{c|}{0.24} \\
\hline
\end{tabular}}
\end{center}
\caption{Results of applying the metaheuristics to the small-sized instances with learning rate $\alpha = 0.2$.}
\label{tab10}
\end{table}

\begin{table}[ht!]
\begin{center}
\resizebox{!}{0.45\textheight}{
\begin{tabular}{|c|cc|cc|cc|cc|cc|cc|cc|}
\hline
\multirow{2}{*}{instance} & \multicolumn{2}{c|}{ILS-RN} & \multicolumn{2}{c|}{ILS-CN} & \multicolumn{2}{c|}{GRASP-RN} & \multicolumn{2}{c|}{GRASP-CN} & \multicolumn{2}{c|}{TS-RN} & \multicolumn{2}{c|}{TS-CN} & \multicolumn{2}{c|}{SA} \\
\cline{2-15} 
& $C_{\max}$ & Time & $C_{\max}$ & Time & $C_{\max}$ & Time & $C_{\max}$ & Time & $C_{\max}$ & Time & $C_{\max}$ & Time & $C_{\max}$ & Time \\
\hline
\hline
miniDAFJS01 & \textbf{19,443} & 0.03 & \textbf{19,443} & 0.06 & \textbf{19,443} & 24.05 & \textbf{19,443} & 54.44 & 19,652 & 0.00 & 19,702 & 0.00 & \textbf{19,443} & 0.34 \\
miniDAFJS02 & \textbf{18,916} & 0.00 & \textbf{18,916} & 0.02 & \textbf{18,916} & 0.02 & \textbf{18,916} & 0.05 & \textbf{18,916} & 0.00 & 21,136 & 0.00 & \textbf{18,916} & 0.25 \\
miniDAFJS03 & \textbf{17,419} & 0.00 & \textbf{17,419} & 0.00 & \textbf{17,419} & 0.00 & \textbf{17,419} & 0.00 & \textbf{17,419} & 0.00 & 17,619 & 0.00 & \textbf{17,419} & 0.16 \\
miniDAFJS04 & \textbf{18,800} & 0.00 & \textbf{18,800} & 0.00 & \textbf{18,800} & 0.00 & \textbf{18,800} & 0.00 & \textbf{18,800} & 0.00 & 20,124 & 0.00 & \textbf{18,800} & 0.13 \\
miniDAFJS05 & \textbf{17,596} & 0.05 & \textbf{17,596} & 0.07 & \textbf{17,596} & 3.01 & \textbf{17,596} & 35.51 & 17,657 & 0.02 & 17,657 & 0.00 & \textbf{17,596} & 0.30 \\
miniDAFJS06 & \textbf{18,692} & 0.03 & \textbf{18,692} & 0.05 & \textbf{18,692} & 0.00 & \textbf{18,692} & 0.04 & \textbf{18,692} & 0.00 & 18,952 & 0.00 & \textbf{18,692} & 0.33 \\
miniDAFJS07 & \textbf{24,256} & 0.00 & \textbf{24,256} & 0.01 & \textbf{24,256} & 0.11 & \textbf{24,256} & 0.43 & \textbf{24,256} & 0.00 & 24,636 & 0.00 & \textbf{24,256} & 0.23 \\
miniDAFJS08 & \textbf{17,900} & 0.00 & \textbf{17,900} & 0.00 & \textbf{17,900} & 0.00 & \textbf{17,900} & 0.00 & \textbf{17,900} & 0.00 & \textbf{17,900} & 0.00 & \textbf{17,900} & 0.03 \\
miniDAFJS09 & \textbf{20,797} & 0.01 & \textbf{20,797} & 0.05 & \textbf{20,797} & 0.05 & \textbf{20,797} & 0.00 & \textbf{20,797} & 0.00 & 21,298 & 0.00 & \textbf{20,797} & 0.32 \\
miniDAFJS10 & \textbf{17,395} & 0.00 & \textbf{17,395} & 0.00 & \textbf{17,395} & 0.00 & \textbf{17,395} & 0.10 & \textbf{17,395} & 0.00 & 19,472 & 0.00 & \textbf{17,395} & 0.22 \\
miniDAFJS11 & \textbf{25,304} & 0.02 & \textbf{25,304} & 0.01 & \textbf{25,304} & 1.06 & \textbf{25,304} & 0.00 & 27,387 & 0.00 & 25,554 & 0.00 & \textbf{25,304} & 0.37 \\
miniDAFJS12 & \textbf{17,105} & 0.01 & \textbf{17,105} & 0.01 & \textbf{17,105} & 1.01 & \textbf{17,105} & 1.98 & \textbf{17,105} & 0.00 & 17,745 & 0.00 & \textbf{17,105} & 0.25 \\
miniDAFJS13 & \textbf{14,077} & 0.00 & \textbf{14,077} & 0.00 & \textbf{14,077} & 0.00 & \textbf{14,077} & 0.00 & \textbf{14,077} & 0.00 & \textbf{14,077} & 0.00 & \textbf{14,077} & 0.01 \\
miniDAFJS14 & \textbf{20,620} & 0.00 & \textbf{20,620} & 0.00 & \textbf{20,620} & 0.00 & \textbf{20,620} & 7.65 & \textbf{20,620} & 0.00 & 20,711 & 0.00 & \textbf{20,620} & 0.28 \\
miniDAFJS15 & \textbf{18,625} & 0.00 & \textbf{18,625} & 0.00 & \textbf{18,625} & 0.00 & \textbf{18,625} & 0.00 & \textbf{18,625} & 0.00 & \textbf{18,625} & 0.00 & \textbf{18,625} & 0.18 \\
miniDAFJS16 & \textbf{22,734} & 0.01 & \textbf{22,734} & 0.01 & \textbf{22,734} & 0.00 & \textbf{22,734} & 4.24 & 22,939 & 0.00 & 22,939 & 0.00 & \textbf{22,734} & 0.29 \\
miniDAFJS17 & \textbf{18,253} & 0.00 & \textbf{18,253} & 0.00 & \textbf{18,253} & 0.00 & \textbf{18,253} & 0.00 & \textbf{18,253} & 0.00 & \textbf{18,253} & 0.00 & \textbf{18,253} & 0.09 \\
miniDAFJS18 & \textbf{16,495} & 0.00 & \textbf{16,495} & 0.00 & \textbf{16,495} & 0.00 & \textbf{16,495} & 0.01 & \textbf{16,495} & 0.00 & 17,270 & 0.00 & \textbf{16,495} & 0.24 \\
miniDAFJS19 & \textbf{18,474} & 0.00 & \textbf{18,474} & 0.00 & \textbf{18,474} & 0.00 & \textbf{18,474} & 0.00 & \textbf{18,474} & 0.00 & 19,030 & 0.00 & \textbf{18,474} & 0.14 \\
miniDAFJS20 & \textbf{18,521} & 0.04 & \textbf{18,521} & 0.03 & \textbf{18,521} & 7.15 & \textbf{18,521} & 0.64 & 19,587 & 0.00 & 18,602 & 0.00 & \textbf{18,521} & 0.27 \\
miniDAFJS21$^+$ & \textbf{19,430} & 0.34 & \textbf{19,430} & 0.52 & \textbf{19,430} & 40.61 & 19,448 & 49.57 & 20,368 & 0.00 & 20,624 & 0.00 & 19,550 & 0.45 \\
miniDAFJS22 & \textbf{22,322} & 0.04 & \textbf{22,322} & 0.22 & \textbf{22,322} & 56.51 & \textbf{22,322} & 84.04 & 22,396 & 0.01 & 22,767 & 0.00 & \textbf{22,322} & 0.42 \\
miniDAFJS23$^+$ & \textbf{20,932} & 1.75 & \textbf{20,932} & 0.63 & 21,372 & 0.02 & 21,372 & 0.02 & 21,176 & 0.01 & 21,540 & 0.01 & 21,031 & 0.27 \\
miniDAFJS24 & \textbf{20,389} & 0.03 & \textbf{20,389} & 0.04 & \textbf{20,389} & 0.01 & \textbf{20,389} & 0.02 & 20,938 & 0.00 & \textbf{20,389} & 0.01 & \textbf{20,389} & 0.39 \\
miniDAFJS25 & \textbf{18,400} & 0.01 & \textbf{18,400} & 0.00 & \textbf{18,400} & 0.00 & \textbf{18,400} & 35.30 & \textbf{18,400} & 0.01 & 18,472 & 0.00 & \textbf{18,400} & 0.36 \\
miniDAFJS26 & \textbf{18,396} & 0.04 & \textbf{18,396} & 0.02 & \textbf{18,396} & 17.47 & \textbf{18,396} & 55.17 & 18,870 & 0.01 & 19,427 & 0.00 & \textbf{18,396} & 0.37 \\
miniDAFJS27 & \textbf{18,501} & 0.12 & \textbf{18,501} & 0.04 & \textbf{18,501} & 79.26 & \textbf{18,501} & 22.10 & \textbf{18,501} & 0.02 & 20,586 & 0.00 & \textbf{18,501} & 0.38 \\
miniDAFJS28 & \textbf{18,762} & 0.11 & \textbf{18,762} & 0.06 & \textbf{18,762} & 0.00 & \textbf{18,762} & 3.62 & \textbf{18,762} & 0.01 & 20,020 & 0.00 & 18,916 & 0.35 \\
miniDAFJS29 & \textbf{18,253} & 0.01 & \textbf{18,253} & 0.00 & \textbf{18,253} & 1.39 & \textbf{18,253} & 0.01 & \textbf{18,253} & 0.00 & 18,560 & 0.00 & \textbf{18,253} & 0.31 \\
miniDAFJS30$^+$ & \textbf{19,504} & 0.14 & \textbf{19,504} & 0.46 & \textbf{19,504} & 31.31 & 20,137 & 30.57 & \textbf{19,504} & 0.27 & 19,659 & 0.00 & 19,618 & 0.38 \\ \hline
$\bar{C}_{\max}$ & \multicolumn{2}{c|}{19,210.37} & \multicolumn{2}{c|}{19,210.37} & \multicolumn{2}{c|}{19,225.03} & \multicolumn{2}{c|}{19,246.73} & \multicolumn{2}{c|}{19,407.13} & \multicolumn{2}{c|}{19,778.20} & \multicolumn{2}{c|}{19,226.60} \\
\#best & \multicolumn{2}{c|}{30} & \multicolumn{2}{c|}{30} & \multicolumn{2}{c|}{29} & \multicolumn{2}{c|}{27} & \multicolumn{2}{c|}{20} & \multicolumn{2}{c|}{5} & \multicolumn{2}{c|}{26} \\
gap(\%) & \multicolumn{2}{c|}{0.00} & \multicolumn{2}{c|}{0.00} & \multicolumn{2}{c|}{0.07} & \multicolumn{2}{c|}{0.18} & \multicolumn{2}{c|}{0.93} & \multicolumn{2}{c|}{3.02} & \multicolumn{2}{c|}{0.08} \\
$\overline{\mathrm{gap}}$(\%) & \multicolumn{2}{c|}{0.00} & \multicolumn{2}{c|}{0.00} & \multicolumn{2}{c|}{0.24} & \multicolumn{2}{c|}{0.32} & \multicolumn{2}{c|}{0.93} & \multicolumn{2}{c|}{3.02} & \multicolumn{2}{c|}{0.61} \\
\hline
\multicolumn{15}{c}{}\\
\hline
\multirow{2}{*}{instance} & \multicolumn{2}{c|}{ILS--RN} & \multicolumn{2}{c|}{ILS--CN} & \multicolumn{2}{c|}{GRASP--RN} & \multicolumn{2}{c|}{GRASP--CN} & \multicolumn{2}{c|}{TS--RN} & \multicolumn{2}{c|}{TS--CN} & \multicolumn{2}{c|}{SA} \\
\cline{2-15} 
& $C_{\max}$ & time & $C_{\max}$ & time & $C_{\max}$ & time & $C_{\max}$ & time & $C_{\max}$ & time & $C_{\max}$ & time & $C_{\max}$ & time \\
\hline
\hline
miniYFJS01 & \textbf{31,008} & 0.00 & \textbf{31,008} & 0.00 & \textbf{31,008} & 0.00 & \textbf{31,008} & 13.25 & \textbf{31,008} & 0.00 & 32,506 & 0.00 & \textbf{31,008} & 0.25 \\
miniYFJS02 & \textbf{22,010} & 0.00 & \textbf{22,010} & 0.02 & \textbf{22,010} & 0.00 & \textbf{22,010} & 0.00 & \textbf{22,010} & 0.01 & 24,146 & 0.00 & \textbf{22,010} & 0.28 \\
miniYFJS03 & \textbf{38,935} & 0.00 & \textbf{38,935} & 0.00 & \textbf{38,935} & 0.00 & \textbf{38,935} & 0.00 & 41,699 & 0.00 & 42,339 & 0.00 & \textbf{38,935} & 0.34 \\
miniYFJS04 & \textbf{23,774} & 0.00 & \textbf{23,774} & 0.00 & \textbf{23,774} & 0.00 & \textbf{23,774} & 0.00 & \textbf{23,774} & 0.00 & 24,017 & 0.00 & \textbf{23,774} & 0.28 \\
miniYFJS05 & \textbf{22,843} & 0.00 & \textbf{22,843} & 0.00 & \textbf{22,843} & 0.00 & \textbf{22,843} & 0.00 & \textbf{22,843} & 0.00 & 23,236 & 0.00 & \textbf{22,843} & 0.24 \\
miniYFJS06 & \textbf{27,366} & 0.00 & \textbf{27,366} & 0.00 & \textbf{27,366} & 0.00 & \textbf{27,366} & 0.00 & \textbf{27,366} & 0.00 & \textbf{27,366} & 0.00 & \textbf{27,366} & 0.00 \\
miniYFJS07 & \textbf{38,932} & 0.00 & \textbf{38,932} & 0.06 & \textbf{38,932} & 0.24 & \textbf{38,932} & 13.11 & \textbf{38,932} & 0.01 & 41,487 & 0.00 & \textbf{38,932} & 0.39 \\
miniYFJS08 & \textbf{29,464} & 0.02 & \textbf{29,464} & 0.04 & \textbf{29,464} & 18.16 & \textbf{29,464} & 9.52 & 29,898 & 0.00 & 30,276 & 0.00 & \textbf{29,464} & 0.43 \\
miniYFJS09 & \textbf{33,763} & 0.01 & \textbf{33,763} & 0.01 & \textbf{33,763} & 0.00 & \textbf{33,763} & 0.00 & \textbf{33,763} & 0.04 & 34,357 & 0.00 & \textbf{33,763} & 0.30 \\
miniYFJS10 & \textbf{25,072} & 0.15 & \textbf{25,072} & 0.10 & \textbf{25,072} & 30.31 & \textbf{25,072} & 49.47 & \textbf{25,072} & 0.00 & 25,093 & 0.00 & \textbf{25,072} & 6.20 \\
miniYFJS11 & \textbf{36,307} & 0.02 & \textbf{36,307} & 0.05 & \textbf{36,307} & 2.55 & \textbf{36,307} & 9.76 & \textbf{36,307} & 0.01 & 39,637 & 0.00 & \textbf{36,307} & 0.44 \\
miniYFJS12 & \textbf{26,219} & 0.04 & \textbf{26,219} & 0.06 & 26,571 & 0.01 & \textbf{26,219} & 0.09 & \textbf{26,219} & 0.22 & 27,829 & 0.00 & \textbf{26,219} & 0.44 \\
miniYFJS13 & \textbf{25,619} & 0.00 & \textbf{25,619} & 0.01 & \textbf{25,619} & 0.00 & \textbf{25,619} & 0.02 & \textbf{25,619} & 0.00 & 25,881 & 0.00 & \textbf{25,619} & 0.36 \\
miniYFJS14 & \textbf{27,428} & 0.03 & \textbf{27,428} & 0.04 & \textbf{27,428} & 0.86 & \textbf{27,428} & 0.09 & \textbf{27,428} & 0.05 & \textbf{27,428} & 0.00 & \textbf{27,428} & 0.52 \\
miniYFJS15 & \textbf{38,256} & 0.01 & \textbf{38,256} & 0.03 & \textbf{38,256} & 0.00 & \textbf{38,256} & 7.91 & \textbf{38,256} & 0.01 & 39,294 & 0.00 & \textbf{38,256} & 0.50 \\
miniYFJS16 & \textbf{29,442} & 0.01 & \textbf{29,442} & 0.01 & \textbf{29,442} & 0.01 & \textbf{29,442} & 0.00 & \textbf{29,442} & 0.01 & 30,386 & 0.00 & \textbf{29,442} & 0.41 \\
miniYFJS17 & \textbf{37,465} & 0.02 & \textbf{37,465} & 0.02 & \textbf{37,465} & 0.00 & \textbf{37,465} & 0.10 & \textbf{37,465} & 0.12 & 38,739 & 0.00 & \textbf{37,465} & 0.39 \\
miniYFJS18 & \textbf{25,067} & 0.02 & \textbf{25,067} & 0.00 & \textbf{25,067} & 0.00 & \textbf{25,067} & 0.00 & \textbf{25,067} & 0.00 & \textbf{25,067} & 0.00 & \textbf{25,067} & 0.41 \\
miniYFJS19 & \textbf{29,207} & 0.01 & \textbf{29,207} & 0.02 & \textbf{29,207} & 0.06 & \textbf{29,207} & 0.29 & \textbf{29,207} & 0.04 & 30,218 & 0.00 & \textbf{29,207} & 0.45 \\
miniYFJS20 & \textbf{27,091} & 0.07 & \textbf{27,091} & 0.03 & \textbf{27,091} & 0.12 & \textbf{27,091} & 0.22 & \textbf{27,091} & 0.04 & 28,590 & 0.00 & \textbf{27,091} & 0.42 \\
miniYFJS21 & \textbf{32,166} & 0.30 & \textbf{32,166} & 0.16 & 32,238 & 0.06 & \textbf{32,166} & 0.04 & \textbf{32,166} & 0.01 & 33,027 & 0.00 & \textbf{32,166} & 0.59 \\
miniYFJS22 & \textbf{28,985} & 0.06 & \textbf{28,985} & 0.06 & \textbf{28,985} & 0.04 & \textbf{28,985} & 1.07 & \textbf{28,985} & 0.00 & 29,154 & 0.00 & \textbf{28,985} & 0.54 \\
miniYFJS23$^+$ & \textbf{35,441} & 0.14 & \textbf{35,441} & 0.18 & \textbf{35,441} & 15.36 & \textbf{35,441} & 33.00 & \textbf{35,441} & 0.20 & 35,961 & 0.00 & \textbf{35,441} & 0.61 \\
miniYFJS24 & \textbf{27,023} & 0.09 & \textbf{27,023} & 0.03 & 27,395 & 52.85 & \textbf{27,023} & 53.26 & \textbf{27,023} & 0.10 & 31,200 & 0.00 & \textbf{27,023} & 0.48 \\
miniYFJS25$^+$ & \textbf{32,346} & 0.07 & \textbf{32,346} & 0.07 & \textbf{32,346} & 0.38 & 32,465 & 0.55 & \textbf{32,346} & 0.02 & 32,513 & 0.00 & \textbf{32,346} & 3.21 \\
miniYFJS26 & \textbf{43,452} & 0.01 & \textbf{43,452} & 0.05 & \textbf{43,452} & 0.14 & \textbf{43,452} & 25.93 & 44,207 & 0.00 & 44,633 & 0.00 & \textbf{43,452} & 0.54 \\
miniYFJS27 & \textbf{31,571} & 0.00 & \textbf{31,571} & 0.01 & \textbf{31,571} & 0.01 & \textbf{31,571} & 0.49 & \textbf{31,571} & 0.01 & 31,877 & 0.00 & \textbf{31,571} & 0.51 \\
miniYFJS28 & \textbf{30,428} & 0.01 & \textbf{30,428} & 0.01 & \textbf{30,428} & 74.61 & \textbf{30,428} & 28.22 & \textbf{30,428} & 0.04 & \textbf{30,428} & 0.00 & \textbf{30,428} & 0.47 \\
miniYFJS29 & \textbf{32,826} & 0.26 & \textbf{32,826} & 0.36 & \textbf{32,826} & 0.31 & \textbf{32,826} & 11.23 & \textbf{32,826} & 0.45 & 33,455 & 0.00 & \textbf{32,826} & 0.48 \\
miniYFJS30 & \textbf{29,848} & 0.20 & \textbf{29,848} & 0.12 & \textbf{29,848} & 5.38 & \textbf{29,848} & 54.66 & \textbf{29,848} & 0.04 & 31,093 & 0.00 & \textbf{29,848} & 0.76 \\ \hline
$\bar{C}_{\max}$ & \multicolumn{2}{c|}{30,645.13} & \multicolumn{2}{c|}{30,645.13} & \multicolumn{2}{c|}{30,671.67} & \multicolumn{2}{c|}{30,649.10} & \multicolumn{2}{c|}{30,776.90} & \multicolumn{2}{c|}{31,707.77} & \multicolumn{2}{c|}{30,645.13} \\
\#best & \multicolumn{2}{c|}{30} & \multicolumn{2}{c|}{30} & \multicolumn{2}{c|}{27} & \multicolumn{2}{c|}{29} & \multicolumn{2}{c|}{27} & \multicolumn{2}{c|}{4} & \multicolumn{2}{c|}{30} \\
gap(\%) & \multicolumn{2}{c|}{0.00} & \multicolumn{2}{c|}{0.00} & \multicolumn{2}{c|}{0.10} & \multicolumn{2}{c|}{0.01} & \multicolumn{2}{c|}{0.34} & \multicolumn{2}{c|}{3.41} & \multicolumn{2}{c|}{0.00} \\
$\overline{\mathrm{gap}}$(\%) & \multicolumn{2}{c|}{0.00} & \multicolumn{2}{c|}{0.00} & \multicolumn{2}{c|}{0.28} & \multicolumn{2}{c|}{0.16} & \multicolumn{2}{c|}{0.34} & \multicolumn{2}{c|}{3.41} & \multicolumn{2}{c|}{0.12} \\ 
\hline
\end{tabular}}
\end{center}
\caption{Results of applying the metaheuristics to the small-sized instances with learning rate $\alpha = 0.3$.}
\label{tab11}
\end{table}

\begin{table}[ht!]
\begin{center}
\resizebox{\textwidth}{!}{
\begin{tabular}{|c|c|c|ccccccc|}
\hline
& & & ILS-RN & ILS-CN & GRASP-RN & GRASP-CN & TS-RN & TS-CN & SA \\
\hline
\hline
\multirow{6}{*}{\rotatebox{90}{$\alpha=0.1$}} & \multirow{3}{*}{DA-type} &
\#optimal                                     & 30   & 30   & 28   & 29   & 22   & 7    & 27   \\
& & gap to optimal (\%)                       & 0.00 & 0.00 & 0.14 & 0.13 & 0.95 & 2.61 & 0.10 \\
& & $\overline{\mathrm{gap}}$ to optimal (\%) & 0.00 & 0.00 & 0.18 & 0.17 & 0.95 & 2.61 & 0.35 \\
\cline{2-10}
& \multirow{3}{*}{Y-type} &
\#optimal                                     & 30   & 30   & 28   & 29   & 25   & 10   & 29   \\
& & gap to optimal (\%)                       & 0.00 & 0.00 & 0.03 & 0.02 & 0.67 & 2.92 & 0.02 \\
& & $\overline{\mathrm{gap}}$ to optimal (\%) & 0.00 & 0.00 & 0.39 & 0.06 & 0.67 & 2.92 & 0.18 \\
\hline
\multirow{6}{*}{\rotatebox{90}{$\alpha=0.2$}} & \multirow{3}{*}{DA-type} &
\#optimal                                     & 28   & 28   & 26   & 27   & 20   & 6    & 26   \\
& & gap to optimal (\%)                       & 0.00 & 0.00 & 0.04 & 0.03 & 1.13 & 2.63 & 0.02 \\
& & $\overline{\mathrm{gap}}$ to optimal (\%) & 0.00 & 0.00 & 0.13 & 0.16 & 1.13 & 2.63 & 0.31 \\
\cline{2-10}
& \multirow{3}{*}{Y-type} &
\#optimal                                     & 30   & 30   & 29   & 28   & 24   & 7    & 30   \\
& & gap to optimal (\%)                       & 0.00 & 0.00 & 0.15 & 0.05 & 0.54 & 3.36 & 0.00 \\
& & $\overline{\mathrm{gap}}$ to optimal (\%) & 0.00 & 0.00 & 0.43 & 0.17 & 0.54 & 3.36 & 0.24 \\
\hline
\multirow{6}{*}{\rotatebox{90}{$\alpha=0.3$}} & \multirow{3}{*}{DA-type} &
\#optimal                                     & 27   & 27   & 27   & 27   & 19   & 5    & 26   \\
& & gap to optimal (\%)                       & 0.00 & 0.00 & 0.00 & 0.00 & 0.81 & 2.99 & 0.03 \\
& & $\overline{\mathrm{gap}}$ to optimal (\%) & 0.00 & 0.00 & 0.13 & 0.11 & 0.81 & 2.99 & 0.48 \\
\cline{2-10}
& \multirow{3}{*}{Y-type} &
\#optimal                                     & 28   & 28   & 25   & 28   & 25   & 4    & 28   \\
& & gap to optimal (\%)                       & 0.00 & 0.00 & 0.11 & 0.00 & 0.37 & 3.58 & 0.00 \\
& & $\overline{\mathrm{gap}}$ to optimal (\%) & 0.00 & 0.00 & 0.29 & 0.14 & 0.37 & 3.58 & 0.13 \\
\hline
\hline
\multicolumn{2}{|c|}{}        & \#optimal                                 & 173  & 173  & 163  & 168  & 135  & 39   & 166  \\
\multicolumn{2}{|c|}{Summary} & gap to optimal (\%)                       & 0.00 & 0.00 & 0.08 & 0.04 & 0.74 & 3.01 & 0.03 \\
\multicolumn{2}{|c|}{}        & $\overline{\mathrm{gap}}$ to optimal (\%) & 0.00 & 0.00 & 0.26 & 0.13 & 0.74 & 3.01 & 0.28 \\
\hline
\end{tabular}}
\end{center}
\caption{Behavior of the seven metaheuristics on the 173 small-sized instances for which a proven optimal solution is known.}
\label{tab12}
\end{table}

\subsection{Experiments with classical instances without sequencing flexibility}

To the authors' knowledge, no other method has been developed so far that applies to the exact problem considered in this paper. For this reason, in order to be able to relate the developed methods with those existing in the literature, we decided to apply them to a problem that is a particular case of the problem considered. We are referring to the FJS with position-based learning effect but \textit{without} sequencing flexibility. A method that also applies to this problem and minimizes the makespan was developed in~\cite{TayebiAraghi2014}. The developed method is a hybrid metaheuristic that mixes genetic algorithms with variable neighborhood search with affinity function, thus named GVNSWAF. This method was chosen because, its careful and detailed description allowed us to reimplement it in C++, the same language in which the introduced methods were implemented. In~\cite{TayebiAraghi2014}, for lack of a better option, random instances were considered, and the Taguchi's robust design method was used to determine the best method's parameters. In our experiments, we considered the parameters determined in~\cite{TayebiAraghi2014}. In that respect, the comparison is fair, since neither method was specifically calibrated for the instances that were considered in this section. Our implementation of GVNSWAF is also available for download at \url{https://github.com/kennedy94/FJS}.

In the experiments of this section, we consider 35 instances of the classic FJS from the literature~\cite{Brandimarte1993,Fattahi2007}, without sequencing flexibility, to which we added the position-based learning effect. As we considered learning effect rates $\alpha \in \{ 0.1, 0.2, 0.3 \}$, we have in total 105 instances. On this set of instances, we compared the methods introduced in the present work with the method introduced in~\cite{TayebiAraghi2014}, which we denote GVNSWAF hereafter. Since the method has random components, it was run five times on each instance, as well as the methods introduced in the present work that contain random components as well. Tables~\ref{tab13}, \ref{tab14}, and~\ref{tab15} show the results. A summary of the results from these tables can be seen in Table~\ref{tab17}. Table~\ref{tab17} shows that the methods that stand out are ILS-RN and ILS-CN. It also shows that the ILS-RN, ILS-CN and SA methods have superior performance to the GVNSWAF method. The GVNSWAF method has a performance similar to the performance of GRASP-RN, GRASP-CN and TS-RN methods and only outperforms the TS-CN methods which already showed the worst performance among all the considered methods. The Wilcoxon test, details of which are shown in the Table~\ref{tabwil2} shows that, with a CPU time limit of 10 seconds, the performance of GVNSWAF is equal to the performance of TS-CN, and GVNSWAF is outperformed by all the other six methods. With the CPU time limit of 5 minutes, GVNSWAF outperforms only TS-CN and is equal to GRASP-RN, but is outperformed by all the other five methods. As a final experiment, we tried to solve with an exact commercial solver models of these classical instances to obtain proven optimal solutions. We were able to obtain proven optimal solutions for 47 instances, out of the total of 105 instances. Table~\ref{tab16} compares the solutions obtained by the eight methods in relation to the known optima on that set of 47 instances. Again, the ILS-RN and ILS-CN methods found all optimal solutions in all five runs of each method/instance pair. The GVNSWAF method also achieved the same result. Apart from TS-CN, all others found solutions with average gap less than 1\%.

In addition to the tests described above, we tested 168 additional instances of the classical FJS from the literature~\cite{Barnes1996,DauzrePrs1997,Hurink1994} which, considering $\alpha \in \{0.1, 0.2, 0.3 \}$, amounted to 504. The results essentially confirmed what we have already shown. For more details, see~\cite{teseKennedy}.

\begin{table}[ht!]
\begin{center}
\resizebox{\textwidth}{!}{
\begin{tabular}{|c|cc|cc|cc|cc|cc|cc|cc|cc|}
\hline
\multirow{2}{*}{instance} & \multicolumn{2}{c|}{ILS-RN} & \multicolumn{2}{c|}{ILS-CN} & \multicolumn{2}{c|}{GRASP-RN} & \multicolumn{2}{c|}{GRASP-CN} & \multicolumn{2}{c|}{TS-RN} & \multicolumn{2}{c|}{TS-CN} & \multicolumn{2}{c|}{SA} & \multicolumn{2}{c|}{GVNSWAF} \\
\cline{2-17} 
& $C_{\max}$ & Time & $C_{\max}$ & Time & $C_{\max}$ & Time & $C_{\max}$ & Time & $C_{\max}$ & Time & $C_{\max}$ & Time & $C_{\max}$ & Time & $C_{\max}$ & Time \\
\hline
\hline
mfjs01 & \textbf{45,306} & 0.01 & \textbf{45,306} & 0.00 & \textbf{45,306} & 0.50 & \textbf{45,306} & 0.02 & 46,264 & 0.00 & 46,264 & 0.00 & \textbf{45,306} & 0.25 & \textbf{45,306} & 0.01 \\
mfjs02 & \textbf{42,986} & 0.04 & \textbf{42,986} & 0.01 & \textbf{42,986} & 4.43 & \textbf{42,986} & 4.89 & \textbf{42,986} & 0.01 & \textbf{42,986} & 0.00 & \textbf{42,986} & 0.26 & \textbf{42,986} & 4.13 \\
mfjs03 & \textbf{45,331} & 0.03 & \textbf{45,331} & 0.02 & \textbf{45,331} & 0.90 & \textbf{45,331} & 0.68 & \textbf{45,331} & 0.01 & \textbf{45,331} & 0.00 & \textbf{45,331} & 0.74 & \textbf{45,331} & 0.08 \\
mfjs04 & \textbf{52,012} & 1.23 & \textbf{52,012} & 1.20 & 52,480 & 116.49 & \textbf{52,012} & 197.36 & \textbf{52,012} & 0.30 & 54,075 & 0.01 & 52,630 & 75.60 & \textbf{52,012} & 65.89 \\
mfjs05 & \textbf{47,630} & 0.23 & \textbf{47,630} & 0.12 & \textbf{47,630} & 133.82 & \textbf{47,630} & 0.22 & \textbf{47,630} & 0.02 & \textbf{47,630} & 0.02 & 49,988 & 0.47 & \textbf{47,630} & 0.34 \\
mfjs06 & \textbf{59,523} & 18.24 & \textbf{59,523} & 3.35 & \textbf{59,523} & 124.16 & \textbf{59,523} & 226.29 & \textbf{59,523} & 1.41 & 60,854 & 0.01 & 60,402 & 0.43 & \textbf{59,523} & 0.40 \\
mfjs07 & \textbf{80,877} & 76.04 & \textbf{80,877} & 36.62 & 82,438 & 24.41 & 81,453 & 28.70 & \textbf{80,877} & 14.62 & 82,686 & 0.06 & 81,371 & 3.15 & 81,364 & 40.51 \\
mfjs08 & 80,687 & 205.78 & \textbf{80,305} & 70.59 & 82,481 & 140.18 & 83,273 & 76.08 & \textbf{80,305} & 216.74 & 83,095 & 0.06 & 82,031 & 1.98 & 82,842 & 7.91 \\
mfjs09 & 96,922 & 28.84 & \textbf{96,236} & 110.23 & 100,332 & 90.60 & 99,159 & 160.79 & 98,028 & 93.66 & 99,576 & 0.58 & 97,358 & 4.42 & 98,417 & 76.61 \\
mfjs10 & 109,183 & 149.45 & \textbf{107,489} & 227.23 & 115,762 & 271.38 & 114,670 & 38.10 & 110,314 & 165.43 & 110,721 & 0.28 & 109,326 & 25.97 & 114,495 & 11.20 \\
\hline
$\bar{C}_{\max}$ & \multicolumn{2}{c|}{66,045.70} & \multicolumn{2}{c|}{65,769.50} & \multicolumn{2}{c|}{67,426.90} & \multicolumn{2}{c|}{67,134.30} & \multicolumn{2}{c|}{66,327.00} & \multicolumn{2}{c|}{67,321.80} & \multicolumn{2}{c|}{66,672.90} & \multicolumn{2}{c|}{66,990.60} \\
\#best & \multicolumn{2}{c|}{7} & \multicolumn{2}{c|}{10} & \multicolumn{2}{c|}{5} & \multicolumn{2}{c|}{6} & \multicolumn{2}{c|}{7} & \multicolumn{2}{c|}{3} & \multicolumn{2}{c|}{3} & \multicolumn{2}{c|}{6} \\
gap(\%) & \multicolumn{2}{c|}{0.28} & \multicolumn{2}{c|}{0.00} & \multicolumn{2}{c|}{1.75} & \multicolumn{2}{c|}{1.41} & \multicolumn{2}{c|}{0.66} & \multicolumn{2}{c|}{2.05} & \multicolumn{2}{c|}{1.33} & \multicolumn{2}{c|}{1.25} \\
$\overline{\mathrm{gap}}$(\%) & \multicolumn{2}{c|}{0.48} & \multicolumn{2}{c|}{0.19} & \multicolumn{2}{c|}{2.32} & \multicolumn{2}{c|}{1.86} & \multicolumn{2}{c|}{0.66} & \multicolumn{2}{c|}{2.05} & \multicolumn{2}{c|}{2.26} & \multicolumn{2}{c|}{2.22} \\
\hline
\multicolumn{17}{c}{} \\
\hline
\multirow{2}{*}{instance} & \multicolumn{2}{c|}{ILS-RN} & \multicolumn{2}{c|}{ILS-CN} & \multicolumn{2}{c|}{GRASP-RN} & \multicolumn{2}{c|}{GRASP-CN} & \multicolumn{2}{c|}{TS-RN} & \multicolumn{2}{c|}{TS-CN} & \multicolumn{2}{c|}{SA} & \multicolumn{2}{c|}{GVNSWAF} \\
\cline{2-17} 
& $C_{\max}$ & Time & $C_{\max}$ & Time & $C_{\max}$ & Time & $C_{\max}$ & Time & $C_{\max}$ & Time & $C_{\max}$ & Time & $C_{\max}$ & Time & $C_{\max}$ & Time \\
\hline
\hline
MK01 & \textbf{3,507} & 173.24 & \textbf{3,507} & 81.29 & 3,529 & 19.87 & 3,530 & 62.29 & 3,529 & 0.68 & 3,714 & 0.00 & 3,529 & 3.50 & 3,714 & 0.28 \\
MK02 & 2,347 & 38.47 & 2,327 & 176.39 & 2,348 & 147.70 & 2,404 & 203.40 & 2,332 & 17.24 & 2,372 & 0.09 & \textbf{2,305} & 84.67 & 2,348 & 31.99 \\
MK03 & \textbf{17,316} & 5.97 & \textbf{17,316} & 1.31 & \textbf{17,316} & 1.61 & \textbf{17,316} & 0.07 & \textbf{17,316} & 1.78 & 17,665 & 0.05 & \textbf{17,316} & 2.83 & \textbf{17,316} & 8.31 \\
MK04 & \textbf{5,173} & 112.01 & \textbf{5,173} & 49.49 & 5,255 & 189.35 & 5,187 & 197.41 & 5,665 & 0.79 & 5,792 & 0.05 & \textbf{5,173} & 11.75 & 5,725 & 3.24 \\
MK05 & 13,674 & 260.60 & 13,666 & 187.22 & 13,747 & 291.55 & 13,648 & 72.05 & 13,657 & 38.48 & 13,877 & 0.27 & \textbf{13,627} & 94.51 & 13,710 & 119.96 \\
MK06 & 5,045 & 183.79 & 5,040 & 114.41 & 5,360 & 129.50 & 5,429 & 59.20 & 5,055 & 287.61 & \textbf{4,982} & 23.00 & \textbf{4,982} & 184.62 & 5,380 & 61.45 \\
MK07 & 11,589 & 173.35 & 11,551 & 259.27 & 11,950 & 291.84 & 12,055 & 162.89 & 11,488 & 266.54 & 11,612 & 14.97 & \textbf{11,220} & 256.56 & 11,551 & 42.73 \\
MK08 & 39,002 & 263.25 & \textbf{38,999} & 245.98 & 39,213 & 198.56 & 39,181 & 135.52 & 39,660 & 1.29 & 39,307 & 0.51 & 39,177 & 49.61 & 39,026 & 265.08 \\
MK09 & 25,263 & 237.96 & \textbf{24,919} & 43.07 & 27,103 & 135.39 & 26,211 & 180.98 & 28,421 & 278.25 & 30,456 & 0.06 & 24,970 & 200.21 & 25,545 & 253.28 \\
MK10 & 16,913 & 153.92 & \textbf{16,802} & 218.45 & 19,141 & 204.30 & 19,186 & 202.84 & 16,855 & 255.05 & 18,605 & 3.03 & 17,060 & 171.62 & 18,359 & 131.33 \\
MK11 & 46,820 & 208.18 & 46,708 & 144.04 & 47,161 & 54.49 & 47,180 & 78.61 & \textbf{46,599} & 44.05 & 46,800 & 5.82 & 46,602 & 133.49 & 47,091 & 203.88 \\
MK12 & 39,877 & 128.22 & \textbf{39,869} & 32.84 & 39,955 & 190.20 & 39,911 & 120.55 & 42,093 & 0.56 & 42,232 & 0.05 & 39,949 & 63.10 & 39,907 & 104.38 \\
MK13 & 32,585 & 216.79 & 32,410 & 234.61 & 35,149 & 185.28 & 34,757 & 284.68 & 32,349 & 221.11 & \textbf{31,960} & 109.03 & 32,954 & 218.41 & 32,967 & 299.84 \\
MK14 & 52,376 & 209.56 & \textbf{52,349} & 165.80 & 53,269 & 1.04 & 52,514 & 212.09 & 56,617 & 6.21 & 56,143 & 0.90 & 52,531 & 271.72 & 52,406 & 109.93 \\
MK15 & 28,343 & 223.43 & \textbf{28,069} & 98.35 & 30,525 & 275.34 & 30,278 & 186.99 & 28,313 & 279.37 & 31,336 & 1.18 & 28,916 & 204.11 & 29,577 & 299.38 \\ \hline
$\bar{C}_{\max}$ & \multicolumn{2}{c|}{22,655.33} & \multicolumn{2}{c|}{22,580.33} & \multicolumn{2}{c|}{23,401.40} & \multicolumn{2}{c|}{23,252.47} & \multicolumn{2}{c|}{23,329.93} & \multicolumn{2}{c|}{23,790.20} & \multicolumn{2}{c|}{22,687.40} & \multicolumn{2}{c|}{22,974.80} \\
\#best & \multicolumn{2}{c|}{3} & \multicolumn{2}{c|}{9} & \multicolumn{2}{c|}{1} & \multicolumn{2}{c|}{1} & \multicolumn{2}{c|}{2} & \multicolumn{2}{c|}{2} & \multicolumn{2}{c|}{6} & \multicolumn{2}{c|}{1} \\
gap(\%) & \multicolumn{2}{c|}{0.82} & \multicolumn{2}{c|}{0.47} & \multicolumn{2}{c|}{4.28} & \multicolumn{2}{c|}{3.99} & \multicolumn{2}{c|}{3.15} & \multicolumn{2}{c|}{5.81} & \multicolumn{2}{c|}{0.63} & \multicolumn{2}{c|}{3.44} \\
$\overline{\mathrm{gap}}$(\%) & \multicolumn{2}{c|}{1.52} & \multicolumn{2}{c|}{0.78} & \multicolumn{2}{c|}{5.25} & \multicolumn{2}{c|}{4.69} & \multicolumn{2}{c|}{3.15} & \multicolumn{2}{c|}{5.81} & \multicolumn{2}{c|}{1.19} & \multicolumn{2}{c|}{5.44} \\
\hline
\multicolumn{17}{c}{} \\
\hline
\multirow{2}{*}{instance} & \multicolumn{2}{c|}{ILS-RN} & \multicolumn{2}{c|}{ILS-CN} & \multicolumn{2}{c|}{GRASP-RN} & \multicolumn{2}{c|}{GRASP-CN} & \multicolumn{2}{c|}{TS-RN} & \multicolumn{2}{c|}{TS-CN} & \multicolumn{2}{c|}{SA} & \multicolumn{2}{c|}{GVNSWAF} \\
\cline{2-17} 
& $C_{\max}$ & Time & $C_{\max}$ & Time & $C_{\max}$ & Time & $C_{\max}$ & Time & $C_{\max}$ & Time & $C_{\max}$ & Time & $C_{\max}$ & Time & $C_{\max}$ & Time \\
\hline
\hline
sfjs01 & \textbf{6,459} & 0.00 & \textbf{6,459} & 0.00 & \textbf{6,459} & 0.00 & \textbf{6,459} & 0.00 & \textbf{6,459} & 0.00 & \textbf{6,459} & 0.00 & \textbf{6,459} & 0.00 & \textbf{6,459} & 0.00 \\
sfjs02 & \textbf{10,271} & 0.00 & \textbf{10,271} & 0.00 & \textbf{10,271} & 0.00 & \textbf{10,271} & 0.00 & \textbf{10,271} & 0.00 & \textbf{10,271} & 0.00 & \textbf{10,271} & 0.00 & \textbf{10,271} & 0.00 \\
sfjs03 & \textbf{20,623} & 0.00 & \textbf{20,623} & 0.00 & \textbf{20,623} & 0.00 & \textbf{20,623} & 0.00 & \textbf{20,623} & 0.00 & 21,716 & 0.00 & \textbf{20,623} & 0.00 & \textbf{20,623} & 0.00 \\
sfjs04 & \textbf{33,429} & 0.00 & \textbf{33,429} & 0.00 & \textbf{33,429} & 0.00 & \textbf{33,429} & 0.00 & \textbf{33,429} & 0.00 & 34,483 & 0.00 & \textbf{33,429} & 0.00 & \textbf{33,429} & 0.00 \\
sfjs05 & \textbf{11,006} & 0.00 & \textbf{11,006} & 0.00 & \textbf{11,006} & 0.00 & \textbf{11,006} & 0.00 & \textbf{11,006} & 0.00 & 12,107 & 0.00 & \textbf{11,006} & 0.00 & \textbf{11,006} & 0.00 \\
sfjs06 & \textbf{29,926} & 0.00 & \textbf{29,926} & 0.00 & \textbf{29,926} & 0.00 & \textbf{29,926} & 0.00 & 31,835 & 0.00 & 32,057 & 0.00 & \textbf{29,926} & 0.00 & \textbf{29,926} & 0.00 \\
sfjs07 & \textbf{37,824} & 0.00 & \textbf{37,824} & 0.00 & \textbf{37,824} & 0.00 & \textbf{37,824} & 0.00 & \textbf{37,824} & 0.00 & \textbf{37,824} & 0.00 & \textbf{37,824} & 0.00 & \textbf{37,824} & 0.00 \\
sfjs08 & \textbf{23,842} & 0.00 & \textbf{23,842} & 0.00 & \textbf{23,842} & 0.00 & \textbf{23,842} & 0.00 & \textbf{23,842} & 0.00 & \textbf{23,842} & 0.00 & \textbf{23,842} & 0.01 & \textbf{23,842} & 0.00 \\
sfjs09 & \textbf{19,406} & 0.00 & \textbf{19,406} & 0.00 & \textbf{19,406} & 0.00 & \textbf{19,406} & 0.00 & \textbf{19,406} & 0.00 & \textbf{19,406} & 0.00 & \textbf{19,406} & 0.01 & \textbf{19,406} & 0.00 \\
sfjs10 & \textbf{49,368} & 0.00 & \textbf{49,368} & 0.00 & \textbf{49,368} & 0.00 & \textbf{49,368} & 0.00 & \textbf{49,368} & 0.00 & \textbf{49,368} & 0.00 & \textbf{49,368} & 0.06 & \textbf{49,368} & 0.00 \\ \hline
$\bar{C}_{\max}$ & \multicolumn{2}{c|}{24,215.40} & \multicolumn{2}{c|}{24,215.40} & \multicolumn{2}{c|}{24,215.40} & \multicolumn{2}{c|}{24,215.40} & \multicolumn{2}{c|}{24,406.30} & \multicolumn{2}{c|}{24,753.30} & \multicolumn{2}{c|}{24,215.40} & \multicolumn{2}{c|}{24,215.40} \\
\#best & \multicolumn{2}{c|}{10} & \multicolumn{2}{c|}{10} & \multicolumn{2}{c|}{10} & \multicolumn{2}{c|}{10} & \multicolumn{2}{c|}{9} & \multicolumn{2}{c|}{6} & \multicolumn{2}{c|}{10} & \multicolumn{2}{c|}{10} \\
gap(\%) & \multicolumn{2}{c|}{0.00} & \multicolumn{2}{c|}{0.00} & \multicolumn{2}{c|}{0.00} & \multicolumn{2}{c|}{0.00} & \multicolumn{2}{c|}{0.64} & \multicolumn{2}{c|}{2.56} & \multicolumn{2}{c|}{0.00} & \multicolumn{2}{c|}{0.00} \\
$\overline{\mathrm{gap}}$(\%) & \multicolumn{2}{c|}{0.00} & \multicolumn{2}{c|}{0.00} & \multicolumn{2}{c|}{0.00} & \multicolumn{2}{c|}{0.00} & \multicolumn{2}{c|}{0.64} & \multicolumn{2}{c|}{2.56} & \multicolumn{2}{c|}{0.00} & \multicolumn{2}{c|}{0.00} \\
\hline
\end{tabular}}
\end{center}
\caption{Results of applying the metaheuristics and the method introduced in~\cite{TayebiAraghi2014} to classical instances of the FSJ with learning effect and without sequencing flexibility, with learning effect rate $\alpha = 0.1$.}
\label{tab13}
\end{table}

\begin{table}[ht!]
\begin{center}
\resizebox{\textwidth}{!}{
\begin{tabular}{|c|cc|cc|cc|cc|cc|cc|cc|cc|}
\hline
\multirow{2}{*}{instance} & \multicolumn{2}{c|}{ILS-RN} & \multicolumn{2}{c|}{ILS-CN} & \multicolumn{2}{c|}{GRASP-RN} & \multicolumn{2}{c|}{GRASP-CN} & \multicolumn{2}{c|}{TS-RN} & \multicolumn{2}{c|}{TS-CN} & \multicolumn{2}{c|}{SA} & \multicolumn{2}{c|}{GVNSWAF} \\
\cline{2-17} 
& $C_{\max}$ & Time & $C_{\max}$ & Time & $C_{\max}$ & Time & $C_{\max}$ & Time & $C_{\max}$ & Time & $C_{\max}$ & Time & $C_{\max}$ & Time & $C_{\max}$ & Time \\
\hline
\hline
mfjs01 & \textbf{43,208} & 0.01 & \textbf{43,208} & 0.01 & \textbf{43,208} & 2.14 & \textbf{43,208} & 1.04 & \textbf{43,208} & 0.01 & 44,880 & 0.00 & \textbf{43,208} & 0.33 & \textbf{43,208} & 0.01 \\
mfjs02 & \textbf{41,273} & 0.07 & \textbf{41,273} & 0.08 & \textbf{41,273} & 5.09 & \textbf{41,273} & 0.00 & \textbf{41,273} & 0.01 & 45,208 & 0.00 & \textbf{41,273} & 0.30 & \textbf{41,273} & 47.06 \\
mfjs03 & \textbf{43,412} & 0.11 & \textbf{43,412} & 0.12 & \textbf{43,412} & 0.31 & \textbf{43,412} & 0.16 & \textbf{43,412} & 0.05 & \textbf{43,412} & 0.00 & \textbf{43,412} & 0.39 & \textbf{43,412} & 4.23 \\
mfjs04 & \textbf{47,717} & 2.80 & \textbf{47,717} & 0.90 & 49,024 & 78.72 & 49,024 & 0.19 & \textbf{47,717} & 0.74 & \textbf{47,717} & 0.09 & 49,921 & 0.50 & \textbf{47,717} & 1.16 \\
mfjs05 & \textbf{45,670} & 0.13 & \textbf{45,670} & 0.06 & \textbf{45,670} & 0.09 & \textbf{45,670} & 0.04 & \textbf{45,670} & 0.18 & \textbf{45,670} & 0.00 & \textbf{45,670} & 0.44 & \textbf{45,670} & 1.88 \\
mfjs06 & \textbf{56,743} & 3.25 & \textbf{56,743} & 6.66 & \textbf{56,743} & 135.74 & 57,005 & 281.92 & 57,093 & 4.14 & 57,889 & 0.00 & 57,055 & 178.01 & \textbf{56,743} & 3.24 \\
mfjs07 & \textbf{73,865} & 11.61 & \textbf{73,865} & 6.71 & 75,640 & 262.56 & 74,821 & 50.66 & \textbf{73,865} & 0.60 & 77,533 & 0.01 & 74,646 & 70.58 & \textbf{73,865} & 0.73 \\
mfjs08 & 74,561 & 165.80 & \textbf{74,359} & 246.39 & 76,148 & 275.69 & 76,815 & 3.04 & 75,880 & 7.06 & 78,362 & 0.17 & 75,032 & 1.54 & 75,032 & 269.43 \\
mfjs09 & 88,457 & 180.68 & \textbf{87,676} & 9.06 & 92,100 & 219.33 & 90,970 & 63.40 & 89,991 & 137.38 & 90,853 & 0.54 & \textbf{87,676} & 57.99 & 91,006 & 7.48 \\
mfjs10 & 97,980 & 129.30 & 97,929 & 206.55 & 100,069 & 164.40 & 102,337 & 214.69 & \textbf{97,780} & 90.39 & 98,459 & 0.31 & 98,058 & 158.98 & 99,630 & 38.54 \\ \hline
$\bar{C}_{\max}$ & \multicolumn{2}{c|}{61,288.60} & \multicolumn{2}{c|}{61,185.20} & \multicolumn{2}{c|}{62,328.70} & \multicolumn{2}{c|}{62,453.50} & \multicolumn{2}{c|}{61,588.90} & \multicolumn{2}{c|}{62,998.30} & \multicolumn{2}{c|}{61,595.10} & \multicolumn{2}{c|}{61,755.60} \\
\#best & \multicolumn{2}{c|}{7} & \multicolumn{2}{c|}{9} & \multicolumn{2}{c|}{5} & \multicolumn{2}{c|}{4} & \multicolumn{2}{c|}{7} & \multicolumn{2}{c|}{3} & \multicolumn{2}{c|}{5} & \multicolumn{2}{c|}{7} \\
gap(\%) & \multicolumn{2}{c|}{0.14} & \multicolumn{2}{c|}{0.02} & \multicolumn{2}{c|}{1.49} & \multicolumn{2}{c|}{1.62} & \multicolumn{2}{c|}{0.53} & \multicolumn{2}{c|}{3.01} & \multicolumn{2}{c|}{0.74} & \multicolumn{2}{c|}{0.66} \\
$\overline{\mathrm{gap}}$(\%) & \multicolumn{2}{c|}{0.21} & \multicolumn{2}{c|}{0.10} & \multicolumn{2}{c|}{2.14} & \multicolumn{2}{c|}{1.72} & \multicolumn{2}{c|}{0.53} & \multicolumn{2}{c|}{3.01} & \multicolumn{2}{c|}{2.15} & \multicolumn{2}{c|}{2.38} \\
\hline
\multicolumn{17}{c}{} \\
\hline 
\multirow{2}{*}{instance} & \multicolumn{2}{c|}{ILS-RN} & \multicolumn{2}{c|}{ILS-CN} & \multicolumn{2}{c|}{GRASP-RN} & \multicolumn{2}{c|}{GRASP-CN} & \multicolumn{2}{c|}{TS-RN} & \multicolumn{2}{c|}{TS-CN} & \multicolumn{2}{c|}{SA} & \multicolumn{2}{c|}{GVNSWAF} \\
\cline{2-17} 
& $C_{\max}$ & Time & $C_{\max}$ & Time & $C_{\max}$ & Time & $C_{\max}$ & Time & $C_{\max}$ & Time & $C_{\max}$ & Time & $C_{\max}$ & Time & $C_{\max}$ & Time \\
\hline
\hline
MK01 & \textbf{3,041} & 65.31 & \textbf{3,041} & 85.94 & 3,047 & 64.27 & 3,048 & 141.02 & 3,052 & 8.72 & 3,179 & 0.01 & \textbf{3,041} & 26.52 & 3,201 & 0.72 \\
MK02 & \textbf{2,118} & 28.65 & \textbf{2,118} & 37.89 & 2,120 & 92.83 & 2,129 & 97.06 & \textbf{2,118} & 6.64 & 2,134 & 0.13 & \textbf{2,118} & 8.97 & 2,129 & 9.10 \\
MK03 & \textbf{14,781} & 2.17 & \textbf{14,781} & 0.33 & \textbf{14,781} & 1.71 & \textbf{14,781} & 0.16 & \textbf{14,781} & 2.68 & 14,908 & 0.09 & \textbf{14,781} & 3.08 & \textbf{14,781} & 2.68 \\
MK04 & \textbf{4,483} & 53.69 & \textbf{4,483} & 94.48 & 4,484 & 184.38 & 4,484 & 192.52 & 4,738 & 9.04 & 4,983 & 0.03 & \textbf{4,483} & 5.74 & 4,853 & 1.24 \\
MK05 & 10,885 & 160.49 & 10,892 & 266.12 & 10,945 & 189.24 & 10,894 & 220.20 & 10,875 & 61.39 & 10,902 & 3.61 & \textbf{10,818} & 272.66 & 10,940 & 210.49 \\
MK06 & 4,377 & 149.74 & 4,243 & 83.62 & 4,499 & 266.52 & 4,592 & 58.84 & 4,958 & 30.18 & 4,255 & 20.77 & \textbf{4,134} & 289.16 & 4,591 & 32.71 \\
MK07 & 9,469 & 222.13 & \textbf{9,393} & 23.84 & 9,642 & 153.59 & 9,728 & 87.87 & 9,427 & 21.45 & 9,857 & 0.14 & \textbf{9,393} & 22.45 & 9,531 & 20.83 \\
MK08 & 29,784 & 217.05 & \textbf{29,758} & 98.74 & 30,184 & 181.23 & 29,841 & 159.05 & 30,962 & 4.87 & 31,816 & 0.20 & 29,868 & 242.44 & 29,895 & 290.79 \\
MK09 & 20,641 & 295.19 & \textbf{20,359} & 206.55 & 21,454 & 148.53 & 21,176 & 13.41 & 23,186 & 55.94 & 23,806 & 0.52 & 20,435 & 127.97 & 20,939 & 284.25 \\
MK10 & \textbf{13,648} & 275.95 & 13,694 & 185.13 & 15,879 & 145.53 & 15,179 & 81.71 & 14,079 & 231.51 & 13,770 & 16.78 & 14,246 & 217.07 & 14,974 & 145.47 \\
MK11 & 35,747 & 185.40 & \textbf{35,680} & 212.31 & 36,110 & 49.28 & 36,238 & 181.46 & 36,093 & 1.58 & 35,813 & 16.15 & 35,756 & 25.72 & 35,893 & 284.24 \\
MK12 & 31,539 & 218.99 & \textbf{31,512} & 181.46 & 31,654 & 169.42 & 31,592 & 258.89 & 32,605 & 2.07 & 33,358 & 0.23 & 31,625 & 11.56 & 31,586 & 144.52 \\
MK13 & 26,428 & 290.21 & 26,238 & 292.55 & 28,472 & 129.83 & 27,598 & 58.96 & \textbf{25,764} & 197.92 & 26,690 & 8.26 & 26,471 & 45.26 & 27,285 & 122.26 \\
MK14 & 39,864 & 299.22 & \textbf{39,829} & 148.40 & 40,464 & 277.76 & 40,047 & 263.88 & 41,968 & 7.49 & 43,123 & 0.23 & 40,030 & 70.85 & 39,891 & 222.06 \\
MK15 & 23,815 & 300.52 & \textbf{23,442} & 141.00 & 24,740 & 285.86 & 24,859 & 46.53 & 23,472 & 146.44 & 25,612 & 0.71 & 23,862 & 148.48 & 24,466 & 245.22 \\
\hline
$\bar{C}_{\max}$ & \multicolumn{2}{c|}{18,041.33} & \multicolumn{2}{c|}{17,964.20} & \multicolumn{2}{c|}{18,565.00} & \multicolumn{2}{c|}{18,412.40} & \multicolumn{2}{c|}{18,538.53} & \multicolumn{2}{c|}{18,947.07} & \multicolumn{2}{c|}{18,070.73} & \multicolumn{2}{c|}{18,330.33} \\
\#best & \multicolumn{2}{c|}{5} & \multicolumn{2}{c|}{11} & \multicolumn{2}{c|}{1} & \multicolumn{2}{c|}{1} & \multicolumn{2}{c|}{3} & \multicolumn{2}{c|}{0} & \multicolumn{2}{c|}{7} & \multicolumn{2}{c|}{1} \\
gap(\%) & \multicolumn{2}{c|}{0.89} & \multicolumn{2}{c|}{0.37} & \multicolumn{2}{c|}{3.69} & \multicolumn{2}{c|}{3.14} & \multicolumn{2}{c|}{3.87} & \multicolumn{2}{c|}{5.20} & \multicolumn{2}{c|}{0.72} & \multicolumn{2}{c|}{3.46} \\
$\overline{\mathrm{gap}}$(\%) & \multicolumn{2}{c|}{1.58} & \multicolumn{2}{c|}{0.63} & \multicolumn{2}{c|}{4.53} & \multicolumn{2}{c|}{3.88} & \multicolumn{2}{c|}{3.87} & \multicolumn{2}{c|}{5.20} & \multicolumn{2}{c|}{1.18} & \multicolumn{2}{c|}{5.03} \\
\hline
\multicolumn{17}{c}{} \\
\hline
\multirow{2}{*}{instance} & \multicolumn{2}{c|}{ILS-RN} & \multicolumn{2}{c|}{ILS-CN} & \multicolumn{2}{c|}{GRASP-RN} & \multicolumn{2}{c|}{GRASP-CN} & \multicolumn{2}{c|}{TS-RN} & \multicolumn{2}{c|}{TS-CN} & \multicolumn{2}{c|}{SA} & \multicolumn{2}{c|}{GVNSWAF} \\
\cline{2-17} 
& $C_{\max}$ & Time & $C_{\max}$ & Time & $C_{\max}$ & Time & $C_{\max}$ & Time & $C_{\max}$ & Time & $C_{\max}$ & Time & $C_{\max}$ & Time & $C_{\max}$ & Time \\
\hline
\hline
sfjs01 & \textbf{6,328} & 0.00 & \textbf{6,328} & 0.00 & \textbf{6,328} & 0.00 & \textbf{6,328} & 0.00 & \textbf{6,328} & 0.00 & \textbf{6,328} & 0.00 & \textbf{6,328} & 0.00 & \textbf{6,328} & 0.00 \\
sfjs02 & \textbf{9,872} & 0.00 & \textbf{9,872} & 0.00 & \textbf{9,872} & 0.00 & \textbf{9,872} & 0.00 & \textbf{9,872} & 0.00 & \textbf{9,872} & 0.00 & \textbf{9,872} & 0.00 & \textbf{9,872} & 0.00 \\
sfjs03 & \textbf{19,281} & 0.00 & \textbf{19,281} & 0.00 & \textbf{19,281} & 0.00 & \textbf{19,281} & 0.00 & \textbf{19,281} & 0.00 & 20,027 & 0.00 & \textbf{19,281} & 0.00 & \textbf{19,281} & 0.00 \\
sfjs04 & \textbf{31,553} & 0.00 & \textbf{31,553} & 0.00 & \textbf{31,553} & 0.00 & \textbf{31,553} & 0.00 & 32,472 & 0.00 & 32,472 & 0.00 & \textbf{31,553} & 0.00 & \textbf{31,553} & 0.00 \\
sfjs05 & \textbf{10,198} & 0.00 & \textbf{10,198} & 0.00 & \textbf{10,198} & 0.00 & \textbf{10,198} & 0.00 & \textbf{10,198} & 0.00 & 11,209 & 0.00 & \textbf{10,198} & 0.00 & \textbf{10,198} & 0.00 \\
sfjs06 & \textbf{28,024} & 0.00 & \textbf{28,024} & 0.00 & \textbf{28,024} & 0.00 & \textbf{28,024} & 0.00 & \textbf{28,024} & 0.00 & \textbf{28,024} & 0.00 & \textbf{28,024} & 0.00 & \textbf{28,024} & 0.00 \\
sfjs07 & \textbf{36,075} & 0.00 & \textbf{36,075} & 0.00 & \textbf{36,075} & 0.00 & \textbf{36,075} & 0.00 & \textbf{36,075} & 0.00 & \textbf{36,075} & 0.00 & \textbf{36,075} & 0.00 & \textbf{36,075} & 0.00 \\
sfjs08 & \textbf{22,515} & 0.00 & \textbf{22,515} & 0.00 & \textbf{22,515} & 0.00 & \textbf{22,515} & 0.00 & \textbf{22,515} & 0.00 & \textbf{22,515} & 0.00 & \textbf{22,515} & 0.03 & \textbf{22,515} & 0.01 \\
sfjs09 & \textbf{17,552} & 0.00 & \textbf{17,552} & 0.00 & \textbf{17,552} & 0.00 & \textbf{17,552} & 0.00 & \textbf{17,552} & 0.00 & \textbf{17,552} & 0.00 & \textbf{17,552} & 0.01 & \textbf{17,552} & 0.00 \\
sfjs10 & \textbf{47,323} & 0.00 & \textbf{47,323} & 0.00 & \textbf{47,323} & 0.00 & \textbf{47,323} & 0.00 & \textbf{47,323} & 0.00 & \textbf{47,323} & 0.00 & \textbf{47,323} & 0.09 & \textbf{47,323} & 0.00 \\
\hline
$\bar{C}_{\max}$ & \multicolumn{2}{c|}{22,872.10} & \multicolumn{2}{c|}{22,872.10} & \multicolumn{2}{c|}{22,872.10} & \multicolumn{2}{c|}{22,872.10} & \multicolumn{2}{c|}{22,964.00} & \multicolumn{2}{c|}{23,139.70} & \multicolumn{2}{c|}{22,872.10} & \multicolumn{2}{c|}{22,872.10} \\
\#best & \multicolumn{2}{c|}{10} & \multicolumn{2}{c|}{10} & \multicolumn{2}{c|}{10} & \multicolumn{2}{c|}{10} & \multicolumn{2}{c|}{9} & \multicolumn{2}{c|}{7} & \multicolumn{2}{c|}{10} & \multicolumn{2}{c|}{10} \\
gap(\%) & \multicolumn{2}{c|}{0.00} & \multicolumn{2}{c|}{0.00} & \multicolumn{2}{c|}{0.00} & \multicolumn{2}{c|}{0.00} & \multicolumn{2}{c|}{0.29} & \multicolumn{2}{c|}{1.67} & \multicolumn{2}{c|}{0.00} & \multicolumn{2}{c|}{0.00} \\
$\overline{\mathrm{gap}}$(\%) & \multicolumn{2}{c|}{0.00} & \multicolumn{2}{c|}{0.00} & \multicolumn{2}{c|}{0.00} & \multicolumn{2}{c|}{0.00} & \multicolumn{2}{c|}{0.29} & \multicolumn{2}{c|}{1.67} & \multicolumn{2}{c|}{0.00} & \multicolumn{2}{c|}{0.00} \\
\hline
\end{tabular}}
\end{center}
\caption{Results of applying the metaheuristics and the method introduced in~\cite{TayebiAraghi2014} to classical instances of the FSJ with learning effect and without sequencing flexibility, with learning effect rate $\alpha = 0.2$.}
\label{tab14}
\end{table}

\begin{table}[ht!]
\begin{center}
\resizebox{\textwidth}{!}{
\begin{tabular}{|c|cc|cc|cc|cc|cc|cc|cc|cc|}
\hline
\multirow{2}{*}{instance} & \multicolumn{2}{c|}{ILS-RN} & \multicolumn{2}{c|}{ILS-CN} & \multicolumn{2}{c|}{GRASP-RN} & \multicolumn{2}{c|}{GRASP-CN} & \multicolumn{2}{c|}{TS-RN} & \multicolumn{2}{c|}{TS-CN} & \multicolumn{2}{c|}{SA} & \multicolumn{2}{c|}{GVNSWAF} \\
\cline{2-17} 
& $C_{\max}$ & Time & $C_{\max}$ & Time & $C_{\max}$ & Time & $C_{\max}$ & Time & $C_{\max}$ & Time & $C_{\max}$ & Time & $C_{\max}$ & Time & $C_{\max}$ & Time \\
\hline
\hline
mfjs01 & \textbf{40,508} & 0.02 & \textbf{40,508} & 0.03 & \textbf{40,508} & 3.74 & \textbf{40,508} & 1.20 & 42,562 & 0.00 & 41,785 & 0.00 & \textbf{40,508} & 0.37 & \textbf{40,508} & 0.18 \\
mfjs02 & \textbf{38,996} & 0.05 & \textbf{38,996} & 0.05 & \textbf{38,996} & 2.80 & \textbf{38,996} & 1.60 & \textbf{38,996} & 0.00 & 39,834 & 0.00 & \textbf{38,996} & 0.29 & \textbf{38,996} & 4.26 \\
mfjs03 & \textbf{41,318} & 0.06 & \textbf{41,318} & 0.04 & \textbf{41,318} & 72.50 & \textbf{41,318} & 14.85 & \textbf{41,318} & 0.00 & 44,254 & 0.00 & \textbf{41,318} & 0.26 & \textbf{41,318} & 7.82 \\
mfjs04 & \textbf{44,869} & 2.16 & \textbf{44,869} & 0.60 & \textbf{44,869} & 79.94 & 46,048 & 61.15 & 46,048 & 0.03 & 46,558 & 0.01 & 46,048 & 25.11 & \textbf{44,869} & 101.21 \\
mfjs05 & \textbf{44,376} & 0.21 & \textbf{44,376} & 0.06 & \textbf{44,376} & 90.04 & 44,738 & 20.87 & \textbf{44,376} & 0.33 & 44,738 & 0.00 & \textbf{44,376} & 28.30 & \textbf{44,376} & 0.29 \\
mfjs06 & \textbf{53,618} & 2.31 & \textbf{53,618} & 10.19 & 53,760 & 84.79 & 53,760 & 37.32 & 53,760 & 0.39 & 57,329 & 0.00 & \textbf{53,618} & 184.69 & \textbf{53,618} & 184.69 \\
mfjs07 & \textbf{69,086} & 9.26 & \textbf{69,086} & 4.61 & 70,577 & 192.72 & 70,523 & 80.50 & \textbf{69,086} & 12.88 & 74,941 & 0.01 & \textbf{69,086} & 0.84 & 70,583 & 0.02 \\
mfjs08 & \textbf{68,053} & 152.28 & \textbf{68,053} & 188.36 & 71,413 & 4.05 & 70,333 & 154.93 & 69,845 & 12.97 & 70,333 & 0.26 & 68,362 & 11.09 & \textbf{68,053} & 74.29 \\
mfjs09 & \textbf{79,947} & 3.22 & 80,277 & 249.67 & 83,975 & 92.83 & 82,272 & 203.95 & 80,983 & 27.37 & 81,939 & 0.83 & 80,897 & 1.69 & 81,411 & 40.30 \\
mfjs10 & 89,188 & 217.74 & \textbf{88,515} & 195.67 & 91,930 & 114.34 & 91,477 & 158.34 & 89,644 & 41.99 & 90,178 & 1.03 & 90,419 & 7.67 & 90,144 & 117.97 \\
\hline
$\bar{C}_{\max}$ & \multicolumn{2}{c|}{56,995.90} & \multicolumn{2}{c|}{56,961.60} & \multicolumn{2}{c|}{58,172.20} & \multicolumn{2}{c|}{57,997.30} & \multicolumn{2}{c|}{57,661.80} & \multicolumn{2}{c|}{59,188.90} & \multicolumn{2}{c|}{57,362.80} & \multicolumn{2}{c|}{57,387.60} \\
\#best & \multicolumn{2}{c|}{9} & \multicolumn{2}{c|}{9} & \multicolumn{2}{c|}{5} & \multicolumn{2}{c|}{3} & \multicolumn{2}{c|}{4} & \multicolumn{2}{c|}{0} & \multicolumn{2}{c|}{6} & \multicolumn{2}{c|}{7} \\
gap(\%) & \multicolumn{2}{c|}{0.08} & \multicolumn{2}{c|}{0.04} & \multicolumn{2}{c|}{1.63} & \multicolumn{2}{c|}{1.54} & \multicolumn{2}{c|}{1.32} & \multicolumn{2}{c|}{4.01} & \multicolumn{2}{c|}{0.64} & \multicolumn{2}{c|}{0.58} \\
$\overline{\mathrm{gap}}$(\%) & \multicolumn{2}{c|}{0.21} & \multicolumn{2}{c|}{0.16} & \multicolumn{2}{c|}{2.09} & \multicolumn{2}{c|}{1.97} & \multicolumn{2}{c|}{1.32} & \multicolumn{2}{c|}{4.01} & \multicolumn{2}{c|}{1.89} & \multicolumn{2}{c|}{2.64} \\
\hline
\multicolumn{17}{c}{} \\
\hline
\multirow{2}{*}{instance} & \multicolumn{2}{c|}{ILS-RN} & \multicolumn{2}{c|}{ILS-CN} & \multicolumn{2}{c|}{GRASP-RN} & \multicolumn{2}{c|}{GRASP-CN} & \multicolumn{2}{c|}{TS-RN} & \multicolumn{2}{c|}{TS-CN} & \multicolumn{2}{c|}{SA} & \multicolumn{2}{c|}{GVNSWAF} \\
\cline{2-17} 
& $C_{\max}$ & Time & $C_{\max}$ & Time & $C_{\max}$ & Time & $C_{\max}$ & Time & $C_{\max}$ & Time & $C_{\max}$ & Time & $C_{\max}$ & Time & $C_{\max}$ & Time \\
\hline
\hline
MK01 & \textbf{2,618} & 111.30 & \textbf{2,618} & 79.93 & 2,669 & 26.88 & 2,668 & 40.52 & 2,669 & 38.63 & 2,793 & 0.00 & \textbf{2,618} & 10.59 & 2,714 & 1.38 \\
MK02 & \textbf{1,920} & 7.79 & \textbf{1,920} & 3.27 & \textbf{1,920} & 80.62 & 1,921 & 89.79 & \textbf{1,920} & 1.45 & 1,932 & 0.03 & \textbf{1,920} & 1.90 & 1,921 & 1.74 \\
MK03 & \textbf{12,439} & 3.89 & \textbf{12,439} & 9.44 & \textbf{12,439} & 3.27 & \textbf{12,439} & 2.08 & \textbf{12,439} & 11.02 & 12,669 & 0.03 & \textbf{12,439} & 3.99 & 12,481 & 2.88 \\
MK04 & 3,877 & 296.61 & 3,885 & 298.13 & 3,904 & 117.29 & 3,917 & 237.10 & 3,917 & 6.97 & 4,282 & 0.08 & \textbf{3,872} & 253.43 & 4,246 & 4.82 \\
MK05 & 8,693 & 256.96 & 8,692 & 246.56 & 8,736 & 276.98 & 8,716 & 163.93 & 8,708 & 56.12 & 8,688 & 3.02 & 8,678 & 63.11 & \textbf{8,675} & 266.29 \\
MK06 & 3,609 & 168.58 & 3,611 & 232.13 & 3,830 & 266.93 & 3,784 & 120.00 & 3,866 & 188.73 & \textbf{3,556} & 40.62 & \textbf{3,556} & 148.27 & 4,183 & 267.07 \\
MK07 & 7,913 & 46.15 & 7,914 & 224.48 & 8,055 & 103.05 & 8,015 & 216.96 & 7,938 & 19.54 & 8,036 & 0.17 & \textbf{7,879} & 54.92 & 7,938 & 7.96 \\
MK08 & 23,265 & 117.12 & \textbf{23,240} & 198.23 & 23,413 & 30.22 & 23,442 & 246.91 & 23,813 & 12.18 & 24,437 & 0.21 & 23,333 & 44.37 & 23,385 & 172.23 \\
MK09 & 16,983 & 219.10 & \textbf{16,761} & 121.54 & 17,571 & 206.91 & 17,332 & 260.19 & 18,955 & 47.76 & 19,508 & 0.22 & 16,833 & 61.28 & 17,214 & 48.30 \\
MK10 & 11,576 & 263.74 & \textbf{11,172} & 130.85 & 12,995 & 278.38 & 12,732 & 234.83 & 11,828 & 297.77 & 12,504 & 2.51 & 11,625 & 155.35 & 12,039 & 220.25 \\
MK11 & 27,604 & 208.94 & \textbf{27,542} & 125.25 & 27,786 & 235.45 & 27,990 & 57.88 & 27,692 & 81.59 & 27,644 & 10.97 & 27,694 & 74.08 & 27,757 & 275.70 \\
MK12 & 25,138 & 187.88 & \textbf{25,097} & 127.56 & 25,301 & 30.19 & 25,197 & 72.47 & 25,882 & 2.39 & 26,101 & 0.28 & 25,207 & 44.11 & 25,188 & 268.75 \\
MK13 & 21,633 & 143.11 & 21,150 & 227.99 & 22,728 & 141.57 & 22,682 & 15.97 & \textbf{21,020} & 133.81 & 21,773 & 3.22 & 21,466 & 47.61 & 22,182 & 214.19 \\
MK14 & 30,597 & 134.76 & \textbf{30,550} & 94.94 & 31,067 & 283.03 & 30,815 & 124.47 & 31,674 & 131.73 & 31,223 & 1.77 & 30,778 & 88.62 & 30,672 & 280.78 \\
MK15 & 19,882 & 260.19 & \textbf{19,775} & 196.63 & 20,792 & 10.26 & 20,570 & 79.55 & 20,477 & 266.85 & 24,081 & 0.19 & 19,900 & 269.25 & 21,345 & 56.19 \\
\hline
$\bar{C}_{\max}$ & \multicolumn{2}{c|}{14,516.47} & \multicolumn{2}{c|}{14,424.40} & \multicolumn{2}{c|}{14,880.40} & \multicolumn{2}{c|}{14,814.67} & \multicolumn{2}{c|}{14,853.20} & \multicolumn{2}{c|}{15,281.80} & \multicolumn{2}{c|}{14,519.87} & \multicolumn{2}{c|}{14,796.00} \\
\#best & \multicolumn{2}{c|}{3} & \multicolumn{2}{c|}{10} & \multicolumn{2}{c|}{2} & \multicolumn{2}{c|}{1} & \multicolumn{2}{c|}{3} & \multicolumn{2}{c|}{1} & \multicolumn{2}{c|}{6} & \multicolumn{2}{c|}{1} \\
gap(\%) & \multicolumn{2}{c|}{0.75} & \multicolumn{2}{c|}{0.21} & \multicolumn{2}{c|}{3.46} & \multicolumn{2}{c|}{2.99} & \multicolumn{2}{c|}{3.02} & \multicolumn{2}{c|}{5.82} & \multicolumn{2}{c|}{0.63} & \multicolumn{2}{c|}{3.88} \\
$\overline{\mathrm{gap}}$(\%) & \multicolumn{2}{c|}{1.50} & \multicolumn{2}{c|}{0.49} & \multicolumn{2}{c|}{4.23} & \multicolumn{2}{c|}{3.54} & \multicolumn{2}{c|}{3.02} & \multicolumn{2}{c|}{5.82} & \multicolumn{2}{c|}{0.89} & \multicolumn{2}{c|}{5.50} \\
\hline
\multicolumn{17}{c}{} \\
\hline
\multirow{2}{*}{instance} & \multicolumn{2}{c|}{ILS-RN} & \multicolumn{2}{c|}{ILS-CN} & \multicolumn{2}{c|}{GRASP-RN} & \multicolumn{2}{c|}{GRASP-CN} & \multicolumn{2}{c|}{TS-RN} & \multicolumn{2}{c|}{TS-CN} & \multicolumn{2}{c|}{SA} & \multicolumn{2}{c|}{GVNSWAF} \\
\cline{2-17} 
& $C_{\max}$ & Time & $C_{\max}$ & Time & $C_{\max}$ & Time & $C_{\max}$ & Time & $C_{\max}$ & Time & $C_{\max}$ & Time & $C_{\max}$ & Time & $C_{\max}$ & Time \\
\hline
\hline
sfjs01 & \textbf{6,206} & 0.00 & \textbf{6,206} & 0.00 & \textbf{6,206} & 0.00 & \textbf{6,206} & 0.00 & \textbf{6,206} & 0.00 & \textbf{6,206} & 0.00 & \textbf{6,206} & 0.00 & \textbf{6,206} & 0.00 \\
sfjs02 & \textbf{9,498} & 0.00 & \textbf{9,498} & 0.00 & \textbf{9,498} & 0.00 & \textbf{9,498} & 0.00 & \textbf{9,498} & 0.00 & \textbf{9,498} & 0.00 & \textbf{9,498} & 0.00 & \textbf{9,498} & 0.00 \\
sfjs03 & \textbf{18,062} & 0.00 & \textbf{18,062} & 0.00 & \textbf{18,062} & 0.00 & \textbf{18,062} & 0.00 & \textbf{18,062} & 0.00 & \textbf{18,062} & 0.00 & \textbf{18,062} & 0.00 & \textbf{18,062} & 0.00 \\
sfjs04 & \textbf{29,852} & 0.00 & \textbf{29,852} & 0.00 & \textbf{29,852} & 0.00 & \textbf{29,852} & 0.00 & 30,648 & 0.00 & 30,648 & 0.00 & \textbf{29,852} & 0.00 & \textbf{29,852} & 0.00 \\
sfjs05 & \textbf{9,465} & 0.00 & \textbf{9,465} & 0.00 & \textbf{9,465} & 0.00 & \textbf{9,465} & 0.00 & \textbf{9,465} & 0.00 & 10,288 & 0.00 & \textbf{9,465} & 0.00 & \textbf{9,465} & 0.00 \\
sfjs06 & \textbf{26,281} & 0.00 & \textbf{26,281} & 0.00 & \textbf{26,281} & 0.00 & \textbf{26,281} & 0.00 & \textbf{26,281} & 0.00 & \textbf{26,281} & 0.00 & \textbf{26,281} & 0.00 & \textbf{26,281} & 0.00 \\
sfjs07 & \textbf{34,443} & 0.00 & \textbf{34,443} & 0.00 & \textbf{34,443} & 0.00 & \textbf{34,443} & 0.00 & \textbf{34,443} & 0.00 & \textbf{34,443} & 0.00 & \textbf{34,443} & 0.00 & \textbf{34,443} & 0.00 \\
sfjs08 & \textbf{21,309} & 0.00 & \textbf{21,309} & 0.00 & \textbf{21,309} & 0.00 & \textbf{21,309} & 0.00 & 21,715 & 0.00 & \textbf{21,309} & 0.00 & \textbf{21,309} & 0.02 & \textbf{21,309} & 0.00 \\
sfjs09 & \textbf{15,973} & 0.00 & \textbf{15,973} & 0.00 & \textbf{15,973} & 0.00 & \textbf{15,973} & 0.00 & \textbf{15,973} & 0.00 & \textbf{15,973} & 0.00 & \textbf{15,973} & 0.01 & \textbf{15,973} & 0.00 \\
sfjs10 & \textbf{45,450} & 0.00 & \textbf{45,450} & 0.00 & \textbf{45,450} & 0.00 & \textbf{45,450} & 0.00 & \textbf{45,450} & 0.00 & \textbf{45,450} & 0.00 & \textbf{45,450} & 0.01 & \textbf{45,450} & 0.00 \\ \hline
$\bar{C}_{\max}$ & \multicolumn{2}{c|}{21,653.90} & \multicolumn{2}{c|}{21,653.90} & \multicolumn{2}{c|}{21,653.90} & \multicolumn{2}{c|}{21,653.90} & \multicolumn{2}{c|}{21,774.10} & \multicolumn{2}{c|}{21,815.80} & \multicolumn{2}{c|}{21,653.90} & \multicolumn{2}{c|}{21,653.90} \\
\#best & \multicolumn{2}{c|}{10} & \multicolumn{2}{c|}{10} & \multicolumn{2}{c|}{10} & \multicolumn{2}{c|}{10} & \multicolumn{2}{c|}{8} & \multicolumn{2}{c|}{8} & \multicolumn{2}{c|}{10} & \multicolumn{2}{c|}{10} \\
gap(\%) & \multicolumn{2}{c|}{0.00} & \multicolumn{2}{c|}{0.00} & \multicolumn{2}{c|}{0.00} & \multicolumn{2}{c|}{0.00} & \multicolumn{2}{c|}{0.46} & \multicolumn{2}{c|}{1.14} & \multicolumn{2}{c|}{0.00} & \multicolumn{2}{c|}{0.00} \\
$\overline{\mathrm{gap}}$(\%) & \multicolumn{2}{c|}{0.00} & \multicolumn{2}{c|}{0.00} & \multicolumn{2}{c|}{0.00} & \multicolumn{2}{c|}{0.00} & \multicolumn{2}{c|}{0.46} & \multicolumn{2}{c|}{1.14} & \multicolumn{2}{c|}{0.00} & \multicolumn{2}{c|}{0.00} \\
\hline
\end{tabular}}
\end{center}
\caption{Results of applying the metaheuristics and the method introduced in~\cite{TayebiAraghi2014} to classical instances of the FSJ with learning effect and without sequencing flexibility, with learning effect rate $\alpha = 0.3$.}
\label{tab15}
\end{table}

\begin{table}[ht!]
\begin{center}
\resizebox{!}{0.45\textheight}{
\begin{tabular}{|ccc|cccccccc|}
\hline
Instance & $\alpha$ & $C_{\max}^{\star}$ & ILS-RN & ILS-CN & GRASP-RN & GRASP-CN & TS-RN & TS-CN & SA & GVNSWAF \\
\hline
\hline
mfjs01 & 0.10 & 45,306 & 45,306 & 45,306 & 45,306 & 45,306 & 46,264 & 46,264 & 45,306 & 45,306 \\
mfjs02 & 0.10 & 42,986 & 42,986 & 42,986 & 42,986 & 42,986 & 42,986 & 42,986 & 42,986 & 42,986 \\
mfjs03 & 0.10 & 45,331 & 45,331 & 45,331 & 45,331 & 45,331 & 45,331 & 45,331 & 45,331 & 45,331 \\
mfjs04 & 0.10 & 52,012 & 52,012 & 52,012 & 52,480 & 52,012 & 52,012 & 54,075 & 52,630 & 52,012 \\
mfjs05 & 0.10 & 47,630 & 47,630 & 47,630 & 47,630 & 47,630 & 47,630 & 47,630 & 49,988 & 47,630 \\
mfjs06 & 0.10 & 59,523 & 59,523 & 59,523 & 59,523 & 59,523 & 59,523 & 60,854 & 60,402 & 59,523 \\
sfjs01 & 0.10 &  6,459 &  6,459 &  6,459 &  6,459 &  6,459 &  6,459 &  6,459 &  6,459 &  6,459 \\
sfjs02 & 0.10 & 10,271 & 10,271 & 10,271 & 10,271 & 10,271 & 10,271 & 10,271 & 10,271 & 10,271 \\
sfjs03 & 0.10 & 20,623 & 20,623 & 20,623 & 20,623 & 20,623 & 20,623 & 21,716 & 20,623 & 20,623 \\
sfjs04 & 0.10 & 33,429 & 33,429 & 33,429 & 33,429 & 33,429 & 33,429 & 34,483 & 33,429 & 33,429 \\
sfjs05 & 0.10 & 11,006 & 11,006 & 11,006 & 11,006 & 11,006 & 11,006 & 12,107 & 11,006 & 11,006 \\
sfjs06 & 0.10 & 29,926 & 29,926 & 29,926 & 29,926 & 29,926 & 31,835 & 32,057 & 29,926 & 29,926 \\
sfjs07 & 0.10 & 37,824 & 37,824 & 37,824 & 37,824 & 37,824 & 37,824 & 37,824 & 37,824 & 37,824 \\
sfjs08 & 0.10 & 23,842 & 23,842 & 23,842 & 23,842 & 23,842 & 23,842 & 23,842 & 23,842 & 23,842 \\
sfjs09 & 0.10 & 19,406 & 19,406 & 19,406 & 19,406 & 19,406 & 19,406 & 19,406 & 19,406 & 19,406 \\
sfjs10 & 0.10 & 49,368 & 49,368 & 49,368 & 49,368 & 49,368 & 49,368 & 49,368 & 49,368 & 49,368 \\
\hline
\multicolumn{3}{|c|}{\#optimal} & 16 & 16 & 15 & 16 & 14 & 9 & 13 & 16 \\
\multicolumn{3}{|c|}{gap to optimal (\%)} & 0.00 & 0.00 & 0.06 & 0.00 & 0.53 & 2.12 & 0.48 & 0.00 \\
\multicolumn{3}{|c|}{$\overline{\mathrm{gap}}$ to optimal (\%)} & 0.00 & 0.00 & 0.32 & 0.18 & 0.53 & 2.12 & 0.69 & 0.00 \\
\hline
\multicolumn{11}{c}{} \\
\hline
Instance & $\alpha$ & $C_{\max}^{\star}$ & ILS-RN & ILS-CN & GRASP-RN & GRASP-CN & TS-RN & TS-CN & SA & GVNSWAF \\
\hline
mfjs01 & 0.20 & 43,208 & 43,208 & 43,208 & 43,208 & 43,208 & 43,208 & 44,880 & 43,208 & 43,208 \\
mfjs02 & 0.20 & 41,273 & 41,273 & 41,273 & 41,273 & 41,273 & 41,273 & 45,208 & 41,273 & 41,273 \\
mfjs03 & 0.20 & 43,412 & 43,412 & 43,412 & 43,412 & 43,412 & 43,412 & 43,412 & 43,412 & 43,412 \\
mfjs04 & 0.20 & 47,717 & 47,717 & 47,717 & 49,024 & 49,024 & 47,717 & 47,717 & 49,921 & 47,717 \\
mfjs05 & 0.20 & 45,670 & 45,670 & 45,670 & 45,670 & 45,670 & 45,670 & 45,670 & 45,670 & 45,670 \\
mfjs06 & 0.20 & 56,743 & 56,743 & 56,743 & 56,743 & 57,005 & 57,093 & 57,889 & 57,055 & 56,743 \\
sfjs01 & 0.20 &  6,328 &  6,328 &  6,328 &  6,328 &  6,328 &  6,328 &  6,328 &  6,328 &  6,328 \\
sfjs02 & 0.20 &  9,872 &  9,872 &  9,872 &  9,872 &  9,872 &  9,872 &  9,872 &  9,872 &  9,872 \\
sfjs03 & 0.20 & 19,281 & 19,281 & 19,281 & 19,281 & 19,281 & 19,281 & 20,027 & 19,281 & 19,281 \\
sfjs04 & 0.20 & 31,553 & 31,553 & 31,553 & 31,553 & 31,553 & 32,472 & 32,472 & 31,553 & 31,553 \\
sfjs05 & 0.20 & 10,198 & 10,198 & 10,198 & 10,198 & 10,198 & 10,198 & 11,209 & 10,198 & 10,198 \\
sfjs06 & 0.20 & 28,024 & 28,024 & 28,024 & 28,024 & 28,024 & 28,024 & 28,024 & 28,024 & 28,024 \\
sfjs07 & 0.20 & 36,075 & 36,075 & 36,075 & 36,075 & 36,075 & 36,075 & 36,075 & 36,075 & 36,075 \\
sfjs08 & 0.20 & 22,515 & 22,515 & 22,515 & 22,515 & 22,515 & 22,515 & 22,515 & 22,515 & 22,515 \\
sfjs09 & 0.20 & 17,552 & 17,552 & 17,552 & 17,552 & 17,552 & 17,552 & 17,552 & 17,552 & 17,552 \\
sfjs10 & 0.20 & 47,323 & 47,323 & 47,323 & 47,323 & 47,323 & 47,323 & 47,323 & 47,323 & 47,323 \\
\hline
\multicolumn{3}{|c|}{\#optimal} & 16 & 16 & 15 & 14 & 14 & 10 & 14 & 16 \\
\multicolumn{3}{|c|}{gap to optimal (\%)} & 0.00 & 0.00 & 0.17 & 0.20 & 0.22 & 2.01 & 0.32 & 0.00 \\
\multicolumn{3}{|c|}{$\overline{\mathrm{gap}}$ to optimal (\%)} & 0.00 & 0.00 & 0.28 & 0.21 & 0.22 & 2.01 & 0.74 & 0.00 \\
\hline
\multicolumn{11}{c}{} \\
\hline
Instance & $\alpha$ & $C_{\max}^{\star}$ & ILS-RN & ILS-CN & GRASP-RN & GRASP-CN & TS-RN & TS-CN & SA & GVNSWAF \\
\hline
mfjs01 & 0.30 & 40,508 & 40,508 & 40,508 & 40,508 & 40,508 & 42,562 & 41,785 & 40,508 & 40,508 \\
mfjs02 & 0.30 & 38,996 & 38,996 & 38,996 & 38,996 & 38,996 & 38,996 & 39,834 & 38,996 & 38,996 \\
mfjs03 & 0.30 & 41,318 & 41,318 & 41,318 & 41,318 & 41,318 & 41,318 & 44,254 & 41,318 & 41,318 \\
mfjs04 & 0.30 & 44,869 & 44,869 & 44,869 & 44,869 & 46,048 & 46,048 & 46,558 & 46,048 & 44,869 \\
mfjs05 & 0.30 & 44,376 & 44,376 & 44,376 & 44,376 & 44,738 & 44,376 & 44,738 & 44,376 & 44,376 \\
sfjs01 & 0.30 &  6,206 &  6,206 &  6,206 &  6,206 &  6,206 &  6,206 &  6,206 &  6,206 &  6,206 \\
sfjs02 & 0.30 &  9,498 &  9,498 &  9,498 &  9,498 &  9,498 &  9,498 &  9,498 &  9,498 &  9,498 \\
sfjs03 & 0.30 & 18,062 & 18,062 & 18,062 & 18,062 & 18,062 & 18,062 & 18,062 & 18,062 & 18,062 \\
sfjs04 & 0.30 & 29,852 & 29,852 & 29,852 & 29,852 & 29,852 & 30,648 & 30,648 & 29,852 & 29,852 \\
sfjs05 & 0.30 &  9,465 &  9,465 &  9,465 &  9,465 &  9,465 &  9,465 & 10,288 &  9,465 &  9,465 \\
sfjs06 & 0.30 & 26,281 & 26,281 & 26,281 & 26,281 & 26,281 & 26,281 & 26,281 & 26,281 & 26,281 \\
sfjs07 & 0.30 & 34,443 & 34,443 & 34,443 & 34,443 & 34,443 & 34,443 & 34,443 & 34,443 & 34,443 \\
sfjs08 & 0.30 & 21,309 & 21,309 & 21,309 & 21,309 & 21,309 & 21,715 & 21,309 & 21,309 & 21,309 \\
sfjs09 & 0.30 & 15,973 & 15,973 & 15,973 & 15,973 & 15,973 & 15,973 & 15,973 & 15,973 & 15,973 \\
sfjs10 & 0.30 & 45,450 & 45,450 & 45,450 & 45,450 & 45,450 & 45,450 & 45,450 & 45,450 & 45,450 \\
\hline
\multicolumn{3}{|c|}{\#optimal} & 15 & 15 & 15 & 13 & 11 & 8 & 14 & 15 \\
\multicolumn{3}{|c|}{gap to optimal (\%)} & 0.00 & 0.00 & 0.00 & 0.23 & 0.82 & 1.89 & 0.18 & 0.00 \\
\multicolumn{3}{|c|}{$\overline{\mathrm{gap}}$ to optimal (\%)} & 0.00 & 0.00 & 0.18 & 0.26 & 0.82 & 1.89 & 0.33 & 0.00 \\
\hline
\multicolumn{11}{c}{} \\
\hline
\multicolumn{3}{|c|}{\#optimal} & 47 & 47 & 45 & 43 & 39 & 27 & 41 & 47 \\
\multicolumn{3}{|c|}{gap to optimal (\%)} & 0.00 & 0.00 & 0.08 & 0.14 & 0.52 & 2.01 & 0.33 & 0.00 \\
\multicolumn{3}{|c|}{$\overline{\mathrm{gap}}$ to optimal (\%)} & 0.00 & 0.00 & 0.26 & 0.21 & 0.52 & 2.01 & 0.59 & 0.00 \\
\hline
\end{tabular}}
\end{center}
\caption{Performance of the metaheuristics and GVNSWAF in the 47 classical instances for which proved optimal solutions are known.}
\label{tab16}
\end{table}

\begin{table}
\centering
\begin{tabular}{|c|ccc|ccc|} 
\hline
\multirow{2}{*}{comparison} & \multicolumn{3}{c|}{10 seconds of CPU time limit} & \multicolumn{3}{c|}{5 minutes of CPU time limit} \\ 
\cline{2-7}
& $R^+$ & $R^-$ & $p$-value 
& $R^+$ & $R^-$ & $p$-value \\
\hline
ILS-RN   versus GVNSWAF & 2,262 &   153 & 0.000 & 2,004 &    12 & 0.000 \\
ILS-CN   versus GVNSWAF & 2,415 &     0 & 0.000 & 2,016 &     0 & 0.000 \\
GRASP-RN versus GVNSWAF & 1,792 &   623 & 0.000 & 1,016 & 1,064 & 0.872 \\
GRASP-CN versus GVNSWAF & 2,193 &   222 & 0.000 & 1,327 &   689 & 0.029 \\
TS-RN    versus GVNSWAF & 2,132 &   718 & 0.000 & 1,628 &   787 & 0.012 \\
TS-CN    versus GVNSWAF & 1,991 & 1,412 & 0.181 & 1,081 & 2,079 & 0.015 \\
SA       versus GVNSWAF & 2,187 &   228 & 0.000 & 1,953 &   532 & 0.000 \\
\hline
\end{tabular}
\caption{Details of the Wilcoxon test comparing each introduced method versus GVNSWAF, to accept or reject the null hypothesis ``the difference between the two methods follows a symmetrical distribution around zero''.}
\label{tabwil2}
\end{table}

\begin{table}[ht!]
\begin{center}
\resizebox{\textwidth}{!}{
\begin{tabular}{|c|c|cccccccc|}
\hline
\begin{tabular}{c}instance \\ group \end{tabular} & metric & ILS-RN & ILS-CN & GRASP-RN & GRASP-CN & TS-RN & TS-CN & SA & GVNSWAF \\
\hline
\hline
\multirow{4}{*}{large-sized}
& $\bar{C}_{\max}$               & 48,458.93 & 48,319.72 & 50,206.70 & 50,327.07 & 48,125.00 & 49,446.14 & 47,907.31 & -- \\
& \#best                        & 59 & 58 & 9 & 7 & 55 & 10 & 71 & -- \\
& gap(\%)                       & 1.80 & 1.55 & 4.96 & 5.10 & 1.02 & 3.39 & 0.66 & -- \\
& $\overline{\mathrm{gap}}$(\%) & 2.33 & 1.95 & 5.80 & 5.73 & 1.02 & 3.39 & 1.51 & -- \\
\hline
\multirow{4}{*}{small-sized}
& $\bar{C}_{\max}$               & 26,756.48 & 26,756.48 & 26,783.73 & 26,778.00 & 26,977.45 & 27,570.54 & 26,768.89 & -- \\
& \#best                        & 180 & 180 & 168 & 170 & 140 & 39 & 169 & -- \\
& gap(\%)                       & 0.00 & 0.00 & 0.11 & 0.09 & 0.75 & 2.99 & 0.05 & -- \\
& $\overline{\mathrm{gap}}$(\%) & 0.00 & 0.00 & 0.30 & 0.19 & 0.75 & 2.99 & 0.32 & -- \\
\hline
\multirow{4}{*}{classical}
& $\bar{C}_{\max}$               & 31,989.65 & 31,915.34 & 32,565.66 & 32,480.55 & 32,362.34 & 32,975.85 & 32,122.78 & 32,288.28 \\
& \#best                        & 64 & 78 & 39 & 46 & 30 & 53 & 63 & 52 \\
& gap(\%)                       & 0.40 & 0.15 & 2.10 & 1.88 & 1.81 & 3.78 & 0.54 & 1.78 \\
& $\overline{\mathrm{gap}}$(\%) & 0.74 & 0.31 & 2.63 & 2.26 & 1.81 & 3.78 & 1.06 & 2.97 \\
\hline
\end{tabular}}
\end{center}
\caption{Summary of the performance of the seven methods analyzed on the 150 large-sized instances, the 180 small-sized instances, and the 105 classical instances from the literature, the latter without sequencing flexibility.}
\label{tab17}
\end{table}

%               LocalSearch     ECT             EST
% large-sized	55,100.85	62,469.98	63,439.83
% 0	0	0
% 16.18%	33.01%	32.36%
% 16.18%	33.01%	32.36%
% small-sized	28,541.37	32,357.42	32,183.25
% 35	8	8
% 6.38%	19.17%	19.19%
% 6.38%	19.17%	19.19%
% classical	36,614.04	43,377.30	40,774.73
% 15	3	3
% 13.33%	34.26%	26.73%
% 13.33%	34.26%	26.73%

\section{Conclusions}

In this paper we introduced a local search and four trajectory metaheuristics for the flexible job shop scheduling problem with sequencing flexibility and position-based learning effect. For the local search, we showed that, in the presence of learning effect, the classical approach of considering reallocations of operations in the critical path only fails to consider potentially better neighbors than the current solution. Consequently, we proposed a new neighborhood reduction that does not eliminate potentially better neighbors than the current solution and cuts the neighborhood by approximately 50\%. We further proposed a neighborhood cutoff that reduces the neighborhood size significantly (by about an order of magnitude) and finds solutions that are at most 1\% worse. The introduced local search and/or neighborhoods were used in the development of four trajectory metaheuristics. We performed extensive numerical experiments and showed that variants of ILS, TS, and SA stands out for its effectiveness and efficiency. As a whole, we build a test suite that can be used in the development of future work. The methods introduced and the solutions found are freely available.

As future work, we intend to consider different learning effects, which do not depend only on the position of the operation in the machine to which it was attributed. We also intend to consider objective functions that take into account the energy consumption (green scheduling).\\

\noindent
\textbf{Conflict of interest statement:} On behalf of all authors, the
corresponding author states that there is no conflict of interest.\\

\noindent
\textbf{Data availability:} The datasets generated during and/or
analysed during the current study are available in the GitHub
repository, \url{https://github.com/kennedy94/FJS}.

\bibliographystyle{plain}

\bibliography{arbro2024}

\end{document}